\documentclass[11pt]{article}

\usepackage[margin=1.10in]{geometry}
\usepackage{amsmath,amssymb,amsthm,mathtools,amscd}
\usepackage{enumitem,microtype,graphicx,subfigure}
\microtypesetup{expansion=false}
\usepackage{comment,upgreek,eucal,mathrsfs}
\usepackage{tikz,pgfplots,refcount}
\usetikzlibrary{matrix,calc,arrows.meta}
\pgfplotsset{compat=1.18}
\usepackage[all]{xy}
\usepackage[dvipsnames]{xcolor}
\usepackage[colorlinks=true,linkcolor=MidnightBlue,citecolor=MidnightBlue,urlcolor=MidnightBlue]{hyperref}

\hypersetup{pdftitle={Closed-String Mirror Reconstruction in Dimension Four},
 pdfauthor={Yoel Groman},pdfsubject={Closed-string mirror reconstruction from relative symplectic cohomology}}

\newcommand{\kk}{\Bbbk}
\newcommand{\PV}{\operatorname{PV}}
\newcommand{\Der}{\operatorname{Der}}

\newcommand{\Log}{\operatorname{Log}}
\newcommand{\diver}{\operatorname{div}}

\newcommand{\bC}{{\mathbb C}}
\newcommand{\bP}{{\mathbb P}}

\newcommand{\bR}{{\mathbb R}}
\newcommand{\bZ}{{\mathbb Z}}
\newcommand{\bN}{{\mathbb N}}

\newcommand{\bT}{{\mathbb T}}

\newcommand{\cA}{\mathcal A}

\newcommand{\cF}{\mathcal F}

\newcommand{\cO}{\mathcal O}
\newcommand{\cP}{\mathcal P}

\DeclareMathOperator{\val}{val}

\DeclareMathOperator{\id}{id}

\newtheorem{theorem}{Theorem}[section]

\newtheorem{proposition}[theorem]{Proposition}

\newtheorem{lemma}[theorem]{Lemma}

\newtheorem{corollary}[theorem]{Corollary}

\theoremstyle{definition}
\newtheorem{definition}[theorem]{Definition}
\newtheorem{df}[theorem]{Definition}

\theoremstyle{remark}
\newtheorem{remark}[theorem]{Remark}
\newtheorem{rem}[theorem]{Remark}

\title{Closed-String Mirror Reconstruction in Dimension Four}
\author{Yoel Groman}
\date{30 November 2026}
\begin{document}
\maketitle
\begin{abstract}
    We consider Maslov-zero Lagrangian torus fibrations on closed or geometrically finite-type symplectic Calabi--Yau four-manifolds. We relate the relative symplectic cohomology sheaf in all degrees to the analytic geometry of a non-Archimedean mirror. Over the regular locus we identify the relative symplectic cohomologies, together with their restriction maps and BV structures, with analytic polyvector fields equipped with the divergence operator associated with the logarithmic volume form. Instanton corrections are incorporated implicitly, without computing them explicitly, and are associated with the appearance of holomorphic planes under SFT-type neck stretching for non-generic almost complex structures. Under additional assumptions on the singular locus, the mirror extends to a Gorenstein Berkovich analytic space carrying a continuous proper surjection onto the base, modeled over the regular locus by the logarithm map. The space carries a nowhere-vanishing canonical form, and the relative symplectic cohomology sheaf embeds into the sheaf of reflexive polyvector fields. Under these assumptions, this embedding is an isomorphism in degrees zero and two and has finite-length cokernel in degree one.
\end{abstract}

\tableofcontents
\section{Introduction}\label{sec:introduction}

Consider a Maslov-zero Lagrangian torus fibration $\pi:M\to B$ on a
symplectic Calabi--Yau four-manifold. The Strominger--Yau--Zaslow
picture of mirror symmetry \cite{SYZ96} suggests seeking a mirror
whose geometry reflects the Floer theory of $M$. In the
non-Archimedean formulation developed by Kontsevich and Soibelman
\cite{KontsevichSoibelmanTorus}, one seeks such a mirror as an
analytic space over a valued field. The Arnold--Liouville integral
affine structure on the regular locus of $B$ provides local models in
dual analytic tori, while their gluing is expected to incorporate
instanton corrections.

In the K3 case, Kontsevich and Soibelman solved an associated lifting
problem \cite[\S 7.3 and Theorem~5]{koso}. Starting from a sphere with
standard integral affine singularities and a compatible $K$-affine
structure satisfying a local monodromy condition, they constructed a
Berkovich analytic K3 surface realizing these data. Their construction
modifies the semi-flat gluing by analytic automorphisms that play the
role of instanton corrections. From the symplectic point of view, one
would like to construct the mirror directly from a given Lagrangian
torus fibration, with a precise relation between its analytic geometry
and the Floer-theoretic invariants of the original manifold.

In this paper we address this problem using \emph{relative symplectic
cohomology}. To a compact subset $P\subset B$ one associates the
relative symplectic cohomology of its inverse image
\[
 K_P=\pi^{-1}(P),\qquad \cF^*(P)=SH_M^*(K_P;\Lambda).
\]
These groups carry a product and a BV operator. Inclusions of
compact sets induce restriction maps preserving these structures. The
first task is to show that, for suitable domains in the base, the
degree-zero algebras are normal affinoid algebras and the restriction
homomorphisms correspond to inclusions of affinoid domains. We then
seek an interpretation of the higher-degree groups in terms of the
analytic geometry of the resulting space.

This approach was developed in our earlier joint work with Varolg\"unes
\cite[Theorems~1.1 and~1.4, Remark~1.5]{GromanVarolgunes2022} for nodal
Lagrangian torus fibrations on symplectic cluster manifolds, where
reconstruction was established through computations in local models.
Here we obtain reconstruction without prescribing a particular local
model for the singular fibers or computing the local wall-crossing
transformations. The comparison also extends beyond functions: over
the regular locus we identify relative symplectic cohomology in all
degrees with analytic polyvector fields, and across singular fibers we
relate it to reflexive polyvector fields on the reconstructed space.
These comparisons respect the BV structures, with the Floer BV
operator corresponding to divergence determined by the mirror's volume
form.

\subsection{Reconstruction over the regular locus}\label{subsec:intro-regular}

Let $(M^4,\omega)$ be closed or geometrically of finite type,
with $c_1(M)=0$. Let $\pi:M\to B$ be a proper map to a
surface, restricting over a connected open set $B_{\rm reg}$
to a fibration by Maslov-zero Lagrangian tori.
We fix an orientation of $B_{\rm reg}$ and assume that the
fibration-induced trivialization of $K_M$ extends over $M$.
A sufficient condition for this extension is given in
Section~\ref{subsec:background-grading}.
Our coefficients are the universal Novikov field $\Lambda$
over a field $\kk$ of characteristic zero, with
$|T^a|=e^{-a}$.

We work with rigid analytic spaces over $\Lambda$
\cite{Bosch14}. For the rigid analytic spaces appearing below, we
denote the Berkovich realization by $Y^{\rm Berk}$. In particular,
if $Y=\operatorname{Sp}A$ is affinoid, then
$Y^{\rm Berk}=\mathcal M(A)$ is the space of bounded
multiplicative seminorms on $A$ extending the absolute
value of $\Lambda$. We use Berkovich realizations when
considering continuous maps to the base $B$; see also
\cite{Conrad2008,TemkinBerkovich}.

Consider the Arnold--Liouville integral affine structure
on $B_{\rm reg}$. By a regular convex polytope we mean a
compact convex polytope with rational facet directions
contained in a regular integral affine chart. We allow
lower-dimensional polytopes.

Let $L$ be a regular fiber and put
$M_{\mathbb Z}=H_1(L;\mathbb Z)$. Choosing a basis $e_1,e_2$
of $M_{\mathbb Z}$ and an origin gives local action
coordinates on the base. The dual analytic torus
\[
 \mathbb T_L^\vee=
 (\operatorname{Spec}\Lambda[M_{\mathbb Z}])^{\rm an}
\]
has a valuation map on its Berkovich realization given by
\[
 \Log(x)=
 \bigl(-\log|z^{e_1}(x)|,-\log|z^{e_2}(x)|\bigr).
\]
For a convex polytope $Q$ in these coordinates the
candidate mirror is the affinoid torus domain over $Q$,
which we denote by $\Log^{-1}(Q)$.

Writing $z_i=z^{e_i}$, the logarithmic volume form
$\Omega_0=d\log z_1\wedge d\log z_2$ gives the
polyvector fields on this domain a BV structure with
operator $\diver_{\Omega_0}$.

As in \cite{Paper1}, we keep track of the
tropical data of the Floer and analytic algebras.
The corresponding tropical markings, which record
the lattice labels and coefficients of leading terms,
are described in Section~\ref{subsec:regular-leading-terms}.

\begin{theorem}[Reconstruction over the regular locus]
\label{thm:intro-regular-reconstruction}
For every regular convex polytope $Q$, the algebra
$A_Q=\cF^0(Q)$ is a strict affinoid algebra.
The affinoids $Y_Q=\operatorname{Sp}A_Q$ glue by Floer
restriction to a smooth rigid analytic surface
$Y_{\rm reg}$ equipped with a continuous proper surjection
\[
 p_{\rm reg}:Y_{\rm reg}^{\rm Berk}\longrightarrow B_{\rm reg},
 \qquad p_{\rm reg}^{-1}(Q)=Y_Q^{\rm Berk}.
\]
The surface carries a nowhere-vanishing two-form
$\Omega_{\rm reg}$. There are isomorphisms in all degrees
\[
 \cF^q(Q)\simeq
 \Gamma(Y_Q,\wedge^q T_{Y_{\rm reg}}),
\]
compatible with restriction and with the BV structures,
where the BV operator on polyvector fields is
$\diver_{\Omega_{\rm reg}}$.

Fix an action--angle chart over a neighborhood of a
regular convex polytope $P$. For all regular convex
polytopes $Q\subset P$ there are isomorphisms of BV algebras
\[
 \cF^*(Q)\simeq
 \Gamma(\Log^{-1}(Q),\wedge^*T_{\mathbb T_L^\vee}),
\]
compatible with restriction and respecting the tropical
markings, with BV operator $\diver_{\Omega_0}$ on the right.
The induced analytic identifications intertwine
$p_{\rm reg}$ with $\Log$ on Berkovich realizations and
identify $\Omega_{\rm reg}$ with $\Omega_0$ in oriented
coordinates.

The space, its projection and its volume form are
independent of the auxiliary choices, for a fixed
orientation of $B_{\rm reg}$. Finite regular polyhedral
covers pull back to admissible affinoid covers.
\end{theorem}

Related constructions using Lagrangian Floer theory were
developed by Tu \cite{TuReconstruction} and Yuan
\cite{YuanFamilyFloer}, under their respective Floer-theoretic
hypotheses. An essential feature of these constructions is that the
mirror is not simply the dual torus fibration determined by the
integral affine structure. Rather, the gluing maps between the local
analytic tori incorporate corrections coming from holomorphic discs
with boundary on the Lagrangian fibers. The resulting space is
referred to as the instanton-corrected mirror.

Our construction is based on closed string Floer theory, and
the instanton corrections enter in a somewhat different way.
The main ingredient in the construction of the local charts
in Section~\ref{sec:regular} is a neck-stretching argument.
Choosing an exterior almost complex structure $J$ for which the
relevant holomorphic planes are absent, and a compatible
Hamiltonian ray, we obtain a comparison with the Floer theory of the
cotangent bundle of the torus. This gives the local identification
with the algebra of polyvector fields on a domain in the dual
analytic torus as a tropically marked BV algebra.

In passing from one such ray to another, however, the continuation
map need not agree with the identification coming from the
integral affine structure alone. The transition maps obtained from
these continuation maps thus incorporate the instanton corrections
to the semi-flat mirror.
The fact that they preserve the tropical marking ensures that
the corrections are compatible with the projection to the
integral affine base. The fact that they are maps of BV algebras
imposes a further constraint on the corrected transition maps
and is used in establishing the existence of the global volume
form $\Omega_{\rm reg}$.

\subsection{Extension across singular values}\label{subsec:intro-singular}

Suppose now that $B\setminus B_{\rm reg}$ consists of
isolated points. We consider three requirements on a
polygonal disk $P$ surrounding a single singular value $s$.

\begin{itemize}[leftmargin=2.2em,itemsep=0.35em,topsep=0.5em]
\item[$(G_\partial)$]
The boundary of $P$ lies in the regular locus and is
locally convex and Delzant. There is a conical affine
structure near $\partial P$ with outward-pointing Euler
field. We require at least three labeled edges with
distinct adjacent labeled pairs.
\item[$(G_{\rm rad})$]
The conical structure extends over $P\setminus\{s\}$,
and its inward flow gives nested homothetic polygons
forming a neighborhood basis of $s$.
\item[$(T)_j$]
For a regular fiber $L$ near $\partial P$, ordinary
restriction
\[
 H^j(K_P;\kk)\longrightarrow H^j(L;\kk)
\]
is injective.
\end{itemize}

The precise definitions are given in
Section~\ref{sec:detection}.

These conditions are satisfied in a number of familiar
examples. Throughout, we assume the standing hypotheses
of Section~\ref{subsec:intro-regular}.

\begin{itemize}
\item Almost toric fibrations with one focus--focus point
on each singular fiber satisfy all three requirements,
with $(T)_j$ for $j=0,1,2$. For fibers with several
focus--focus points, $(T)_2$ fails.
\item Hyperk\"ahler rotations of isotrivial elliptic
K3 surfaces whose singular fibers are all cuspidal
satisfy all the requirements. For the elliptic
fibration on a Kummer K3 with four fibers of type
$I_0^*$, the conditions $(G_\partial)$, $(G_{\rm rad})$,
$(T)_0$ and $(T)_1$ hold, but $(T)_2$ fails.
\item There are analogous noncompact examples obtained
from rational elliptic surfaces by removing a fiber
at infinity, as well as from the strongly parabolic
$\mathrm{SL}(2,\mathbb C)$ Hitchin system on the
four-punctured sphere. In these examples the same
weaker conditions hold, and $(T)_2$ holds for the
cuspidal fibers.
\end{itemize}

The verification of these examples is given in
Section~\ref{subsec:reconstruction-examples}.

Let $\cP$ consist of the regular convex polytopes and the
admissible singular polygons satisfying $(G_\partial)$,
$(G_{\rm rad})$, $(T)_0$ and $(T)_1$.
We equip $\cP$, with inclusions as morphisms, with the
Grothendieck topology generated by finite set-theoretic
coverings by its members. Intersections admit finite
refinements by members of $\cP$, as explained in
Section~\ref{sec:radial}.
Under the hypotheses of the following theorem, the interiors
of the two-dimensional members of $\cP$ form a neighborhood
basis in $B$.

On a normal analytic surface $Y$ we denote the sheaf
of reflexive polyvector fields by
$\PV_Y^{[q]}=(\wedge^q T_Y)^{**}$, taken to be zero
in negative degrees. A trivialization $\Omega$ of
the canonical sheaf determines a BV operator on
reflexive polyvector fields by considering divergence
with respect to $\Omega$ on the smooth locus.

\begin{theorem}[Global reconstruction]\label{thm:main}
Suppose every singular value admits a polygon satisfying
$(G_\partial)$, $(G_{\rm rad})$ and $(T)_j$ for $j=0,1,2$.

There is a normal Gorenstein rigid analytic surface $Y$
with a continuous proper surjection
\[
 p:Y^{\rm Berk}\longrightarrow B.
\]
Over $B_{\rm reg}$ this is the analytic mirror and
projection of Theorem~\ref{thm:intro-regular-reconstruction}.
The space $Y^{\rm Berk}$ is without boundary.
Every point of $B$ has a neighborhood basis of open sets
$U$ for which $p^{-1}(U)$ is Stein.
The form $\Omega_{\rm reg}$ extends uniquely to a
trivialization $\Omega$ of the canonical sheaf of $Y$.

For each $Q\in\cP$ there is an affinoid domain
$Y_Q=\operatorname{Sp}A_Q$ in $Y$ with
\[
 p^{-1}(Q)=Y_Q^{\rm Berk},\qquad A_Q=\cF^0(Q).
\]
These identifications are compatible with restriction,
and coverings in $\cP$ pull back to admissible
affinoid covers.

There are injective homomorphisms of graded BV algebras
\[
 \Sigma_Q:\cF^*(Q)\longrightarrow
       \Gamma(Y_Q,\PV_Y^{[*]}),
\]
compatible with restriction. The BV operator on the
right is $\diver_\Omega$ on the smooth locus.
Each $\Sigma_Q$ is an isomorphism in all degrees other
than one. In degree one its cokernel is a finite-length
$A_Q$-module supported on the analytic singular locus
of $Y_Q$. In particular, when $Y_Q$ is smooth the
comparison is an isomorphism of BV algebras.

On regular polytopes the comparison agrees with
Theorem~\ref{thm:intro-regular-reconstruction}.
\end{theorem}

The maps $\Sigma_Q$ are given in higher degrees by
iterated Gerstenhaber brackets with functions.
Their construction is described in
Section~\ref{sec:nonsmooth}.
In particular, the degree-one comparison identifies
the reflexive hull of $\cF^1(Q)$ with the analytic
tangent module. We leave open whether the comparison
itself is surjective.

Some of the conclusions of the theorem hold under
weaker hypotheses. Under $(G_\partial)$, $(T)_0$
and $(T)_1$ alone, the degree-zero Floer algebra of a
single singular polygon is a normal affinoid algebra,
and its spectrum is Gorenstein with trivial canonical
sheaf. Under $(G_{\rm rad})$, the projection also
extends over the singular value.

The role of $(T)_2$ is to ensure injectivity of the
degree-two comparison. Under the other hypotheses,
it is equivalent to $b_2(K_P;\kk)=1$.
Without $(T)_2$ the degree-two comparison is still
surjective, and its possible kernel vanishes upon
localization to the smooth locus.
See Theorem~\ref{thm:polygon-higher-degree}
and Section~\ref{sec:radial}.

\subsection{Idea of the proof}

\paragraph{The regular locus.}
The starting point for the proof over the regular locus is the
comparison between intrinsic and ambient symplectic cohomology
developed in \cite{GromanBoundaryDepth}. Given a Liouville domain
$D$ embedded in a symplectic manifold $M$, we can think of
$SH_M^*(D)$ as a deformation of $SH_{\widehat D}^*(D)$, where
$\widehat D$ is the Liouville completion of $D$.
The advantage of this point of view in the present setting is
that the intrinsic theory is well understood. Indeed, the
Viterbo isomorphism, in the BV form proved by Abouzaid
\cite{AbouzaidViterbo}, gives
\[
 SH^*(T^*\bT^2;\kk)
 \simeq H_{2-*}(\mathcal L\bT^2;\kk)
 \simeq \kk[z_1^{\pm1},z_2^{\pm1}]
             \otimes H^*(\bT^2;\kk).
\]
In the supported version over a convex polytope $Q$ in the
base of $T^*\bT^2$, the Laurent polynomial algebra is replaced
by the algebra $A(Q)$ of convergent Laurent series on
$\Log^{-1}(Q)$. The resulting BV algebra is
\[
 SH_{\widehat D}^*(K_Q;\Lambda)
 \simeq A(Q)\otimes_\kk H^*(\bT^2;\kk)
 \simeq \Gamma(\Log^{-1}(Q),\PV^*),
\]
with the BV operator on the right given by logarithmic
divergence.

Section~1.5 of \cite{GromanBoundaryDepth}, as part of a collaboration
with Kyler Siegel, outlines a conjectural description of the ambient
deformation in terms of the Cieliebak--Latschev Maurer--Cartan
element \cite{CieliebakLatschev} obtained by neck stretching.
The relevant geometric objects for our argument are
holomorphic planes in the complement of $D$ with a
negative end asymptotic to a Reeb orbit on $\partial D$.
The case of a Weinstein neighborhood of a Maslov-zero
torus in real dimension four is particularly favorable:
these planes have Fredholm index $-1$ and can thus be
excluded by a generic choice of cylindrical almost
complex structure $J$ on the complement. In this case
we can establish locality directly, without developing
the general Maurer--Cartan description of the deformation.

The technical means for doing so is provided by the relative
filtration of \cite{GromanBoundaryDepth}. Working in a finite
ordinary Novikov window $[a,b)$ on the nonnegative relative
lattice, we can use strictly $S$-shaped Hamiltonians, which
have a flat exterior plateau, despite their degeneracy there.
Taking the plateau sufficiently high, and in the noncompact
case sufficiently wide, excludes exterior generators and
breaking through them from the calculation in the given
window. On the plateau the Floer equation
then reduces to the holomorphic curve equation, and a
neck-stretching argument shows that an escaping curve would
give rise to one of the planes described above. This gives
the collapse of the locality spectral sequence and, by
applying the same argument to the curves defining the Floer
operations, an isomorphism of BV algebras
\[
 \Phi_{J,Q}:SH_M^*(K_Q;\Lambda)
       \xrightarrow{\sim} SH_{\widehat D}^*(K_Q;\Lambda).
\]
The dependence of this identification on the choice of
exterior almost complex structure accounts for the corrections
to the transition maps between the local mirror charts.

\paragraph{The singular locus.}
The starting point is the calculation of local Floer
cohomology in \cite{Paper1}. Associated with an admissible
polygon $P$ is its integral tropicalization $P_{trop}(\bZ)$,
obtained by gluing the monoids
\[
 C_{i,i+1}=\bN\xi_i+\bN\xi_{i+1}
\]
along their common rays, where $\xi_i$ are the primitive
outward conormals to the edges of $P$. After rounding
the polygon, the nonzero integral points label families
of Reeb orbits on the boundary. In sufficiently short
Novikov windows $I$, and at bounded Reeb period, the
Floer cohomology is accordingly described by
\[
 \left(H^*(K_P;\kk)\oplus
 \bigoplus_m H^*(L_{b_m};\kk)\right)
 \otimes_\kk\Lambda_I^{\log},
\]
where $m$ ranges over the nonzero labels up to the period cutoff
and $b_m$ is a base point under the corresponding family. Restriction to a vertex neighborhood
retains precisely the Reeb summands belonging to its
cone. On the ordinary cohomology summand it is ordinary
restriction. Together with $(T)_j$ and the regular locality
theorem, this description leads to torsion-freeness of
the ordinary open Floer lattice in degree $j$ and the
norm formula
\[
 \|a\|_P=\max_{v\in V(P)}\|a|_v\|_v.
\]

The significance of torsion-freeness is that it allows
us to lift classes from short windows. Indeed, writing
$C_t$ for the subcomplex of elements of norm less than
$e^t$, the obstruction to lifting a class in
$H^j(C_t/T^\lambda C_t)$ to $H^j(C_t)$ lies in the
$T^\lambda$-torsion of $H^{j+1}(C_t)$.
Thus, although a priori we know very little about the
analytic functions on the mirror over the punctured
polygon, we can construct for each edge $e_i$ a function
$z_i\in A_P$ whose leading term along the edge is the
normalized Laurent character of exponent $-\xi_i$.
The monomials
\[
 z_i^az_{i+1}^b,\qquad
 a\xi_i+b\xi_{i+1}\in P_{trop}(\bZ),
\]
are then shown to form an orthonormal basis of $A_P$.
This establishes affinoidness and the extension property
\[
 A_P\xrightarrow{\sim}\Gamma(U_\partial,\cO),
\]
where $U_\partial$ is the regular mirror over $\partial P$.
A related basis indexed by tropical integral points was constructed
by Pascaleff for the degree-zero symplectic cohomology of affine
log Calabi--Yau surfaces with maximal boundary
\cite[Theorem~1.2]{PascaleffLogCY}.

The same edge functions allow us to identify restrictions
associated with small inward displacements of the edges
with Weierstrass and Laurent localizations. Under the
radial conicity assumption we iterate this construction
to extend the projection over the singular value $s$,
with
\[
 p^{-1}(s)=\bigcap_{0<r\leq1}X_{P_r}^{\rm Berk}.
\]

Once the degree-zero reconstruction is established,
much of the higher-degree comparison is a matter of
commutative algebra. Reflexive polyvectors on the normal
surface $X_P=\operatorname{Sp}A_P$ are described by the
modules
\[
 \operatorname{Hom}_{A_P}
 (\wedge^q\Omega^{1,\mathrm{an}}_{A_P/\Lambda},A_P).
\]
The extension property for functions thus gives an
extension property for polyvectors and leads to the
trivialization $\Omega$ of the canonical sheaf.
A further application of the lifting argument, together
with Poincar\'e--Lefschetz duality, produces a Floer class
$s_P\in\cF^2(P)$ representing the dual bivector.
Its brackets $\{s_P,f\}$ give Hamiltonian vector fields,
which generate the tangent sheaf on the smooth locus.
Using the BV identities we show that the possible
defects of the comparison are confined to the analytic
singularities of the mirror.

\subsection{Organization of the paper}

Section~\ref{sec:floer-background} recalls the ordinary and relative filtrations on
relative symplectic cohomology and the local Floer groups that enter
the computation. In Section~\ref{sec:collapse} we prove that the
locality spectral sequence collapses near a regular torus.
Section~\ref{sec:regular} compares the Floer operations with those
of the cotangent bundle of the torus and constructs the mirror over
the regular locus.

The remaining sections concern a polygon surrounding a singular
value. Section~\ref{sec:detection} gives conditions under which a
Floer class is detected by its restrictions near the vertices, and
verifies the geometric and topological hypotheses in the examples of
Section~\ref{subsec:reconstruction-examples}. In
Section~\ref{sec:affinoid} we use this to construct the normal
affinoid surface associated to the polygon. Section~\ref{sec:localization} identifies the affinoid
domains obtained by moving its edges slightly inward.
Section~\ref{sec:nonsmooth} compares the higher Floer groups with
reflexive polyvector fields on the surface. Finally,
Section~\ref{sec:radial} glues the singular charts to the regular
mirror and extends the projection over the singular values.

\subsection{Acknowledgments}

The author was supported by the ISF (grant no. 2445/20
and no. 3605/24).

The questions addressed in this paper are part of a broader circle of ideas that Umut Varolg\"unes and I have been pursuing, both jointly and independently, since 2018, and which have also been the subject of our ongoing collaboration with Mohammed Abouzaid. I am grateful to both for the many stimulating discussions over the years.

The basic idea of the reconstruction near a singular fiber appeared, for the case of a smooth mirror, in Section~7 of the first version of my paper on local Floer cohomology
\cite{GromanIndicatorV1}. The idea of the neck-stretching argument also goes back several years. In developing the present paper, I made extensive use of OpenAI's ChatGPT, particularly in working out the commutative algebra used to treat singular mirrors in Section~\ref{sec:nonsmooth} and to prove normality in Section~\ref{sec:affinoid}. The argument in Section~\ref{sec:radial} for extending the fibration beyond the regime of shallow localization was developed largely by ChatGPT in response to my prompts. I used ChatGPT in the same way to work out the verification of the hypotheses in the examples of Section~\ref{subsec:reconstruction-examples}. I also used ChatGPT extensively in correcting and refining the neck-stretching argument, checking the other proofs, and revising the exposition. I take responsibility for the mathematical content of the paper.

\section{Background in Floer theory}\label{sec:floer-background}
Let $K$ be a compact subset of a Liouville domain $D\subset M$.
To compare the Floer theory with support in $K$ in the two ambient
manifolds $M$ and $\widehat D$, we use the relative filtration of
\cite{GromanBoundaryDepth}. Sections~\ref{sec:collapse}
and~\ref{sec:regular} combine this comparison with neck stretching
near a regular torus. For a polygon surrounding a singular value,
we use the local Floer cohomology of indicator functions
\cite{Paper1} to compute restriction to regular vertex caps in
short Novikov windows. Section~\ref{sec:detection} explains how
these calculations detect classes in completed cohomology.

When $M$ is noncompact, we assume it is \emph{geometrically of finite type}
in the sense of \cite[Section~4.2]{GromanVarolgunes2022}. Thus there
are a geometrically bounded compatible almost complex structure $J_0$
and a proper bounded-below exhaustion $f_0$ whose critical set is
compact, such that $X_{f_0}$ and $\nabla X_{f_0}$ are uniformly bounded
in the metric of $J_0$. We fix $\epsilon>0$ sufficiently small and
use Hamiltonians equal to $\epsilon f_0+C$ outside a compact set.
These tails are dissipative and have no one-periodic orbits outside
a sufficiently large compact set
\cite[Section~2.2]{GromanBoundaryDepth}.

Since $c_1(M)=0$, we use the classical transversality
framework for Hamiltonian Floer theory and its operations,
as in \cite[Section~3]{AbouzaidGromanVarolgunes24}.

\subsection{Floer complexes and symplectic cohomology with supports}
\label{subsec:background-floer}
Recall that $\Lambda$ is the universal Novikov field over $\kk$. We use
its ordinary valuation $\nu$ and the corresponding norm
\[
 \nu\Bigl(\sum_j a_jT^{\lambda_j}\Bigr)
 =\inf\{\lambda_j:a_j\ne0\},\qquad |x|=e^{-\nu(x)}.
\]
Let $H:S^1\times M\to\mathbb R$ be an admissible nondegenerate
Hamiltonian, with convention $\omega(X_H,\cdot)=dH$. Its Floer complex
$CF^*(H;\Lambda)$ is the completed direct sum of the orientation
lines $o_\gamma\otimes\Lambda$ of its one-periodic orbits, graded by
their Conley--Zehnder indices. The differential is given by
\[
 d a_{\gamma_-}
 =\sum_{\gamma_+}\;\sum_{u\in\mathcal M_0(\gamma_-,\gamma_+)}
       \epsilon(u)T^{E_{\mathrm{top}}(u)}a_{\gamma_+},
\]
where $u$ runs through the rigid solutions of
\[
 \partial_su+J_{s,t}(u)
            \bigl(\partial_tu-X_{H_{s,t}}(u)\bigr)=0.
\]
The notation suppresses the orientation-line maps. For an ordinary
Floer cylinder, topological energy agrees with geometric energy and
is nonnegative. For an admissible monotone continuation datum,
$\partial_sH_s\geq0$, the same construction gives a continuation map
from the complex at the negative end to that at the positive end,
with $E_{\mathrm{top}}\geq E_{\mathrm{geo}}\geq0$. These maps and the
Floer operations are continuous in the ordinary Novikov norm.

For a compact subset $K\subset M$, symplectic cohomology with support
in $K$ is obtained by approximating the indicator Hamiltonian
\[
 H_K(x)=\begin{cases}0,&x\in K,\\+\infty,&x\notin K.\end{cases}
\]
Choose a cofinal increasing sequence of admissible Hamiltonians
approaching this indicator, together with monotone continuation maps.
The completed mapping telescope of their Floer complexes is denoted
$SC_M^*(K;\Lambda)$, and
\[
 SH_M^*(K;\Lambda)=H^*(SC_M^*(K;\Lambda)).
\]
For finite-type $M$, the approximation uses the prescribed small
exhaustion profile outside a large compact region. See \cite[Section~2]{GromanBoundaryDepth}. If $K_1\subset K_2$, monotone continuation induces the
restriction map
\[
 SH_M^*(K_2;\Lambda)\longrightarrow SH_M^*(K_1;\Lambda).
\]
We call this ambient theory with supports \emph{relative symplectic
cohomology}. If $K$ lies in a compact Liouville domain $D\subset M$
with primitive $\theta$ near $D$, one can also count trajectories in
its completion $\widehat D$. This gives the intrinsic supported complex
$SC^*_{\widehat D,\theta}(K;\kk)$.

\subsection{\texorpdfstring{$S$}{S}-shaped Hamiltonians and the relative filtration}
\label{subsec:background-filtrations}

Fix $K\subset D$ as above and write $V(D)$ for a neighborhood on
which $\theta$ is defined. For the $S$-shaped Hamiltonians of
\cite[Definition~3.5]{GromanBoundaryDepth}, the relative filtration
$\val_{M,\theta}$ compares the ambient and intrinsic complexes.
We recall the features used below.

An $S$-shaped Hamiltonian rises from its values near $K$ to a high,
nearly constant value on the complement of $D$. Write
$\theta=e^\rho\alpha$ in a Liouville collar and let $D_t$ denote the
domain obtained by adjoining the collar up to $\rho=t$. On a cutting
annulus $A=\{\delta'\leq\rho\leq\delta\}$, with
$\delta/2<\delta'<\delta$ sufficiently small, we require $H=h(e^\rho)$ with
$e^\rho h'(e^\rho)=\sigma$, where $\sigma>0$ is smaller than the
least Reeb period of $\alpha$. The almost complex structure is of
contact type on $A$. Outside $D_\delta$ the Hamiltonian is autonomous
and dynamically small, so that its one-periodic orbits are constant,
and its exterior critical values satisfy
\[
 H|_{\rho=\delta}\leq H(p)\leq H|_{\rho=\delta'}+\eta_{\mathrm{end}},
\]
where $\eta_{\mathrm{end}}$ is small. There are no one-periodic
orbits in the cutting annulus; we call the orbits on its inner side
local. The precise smallness and critical-value conditions
are those of \cite[Definition~3.5]{GromanBoundaryDepth}. In the
noncompact case the Hamiltonian has a small dissipative exhaustion
profile at infinity.

For continuation and its homotopies we use a common cutting annulus.
In addition to $\partial_sH_s\geq0$ we require
\begin{equation}\label{eq:rc-normalized-monotonicity}
 c(s)=H_s|_{\rho=\delta},\qquad
 \partial_s(H_s-c(s))\geq0
 \quad\text{on }M\setminus D_{\delta'},
\end{equation}
with the derivative bounds on the exterior from \cite{GromanBoundaryDepth}.
For a cofinal sequence we may keep $H_s-c(s)$ fixed there, adding
constants on the exterior as we increase the interior slopes
\cite[Proposition~3.9]{GromanBoundaryDepth}. In the noncompact case
we also fix the coefficient of the exhaustion at infinity.

For a local periodic orbit $\gamma$ set
\[
 A_{H,\theta}(\gamma)
       =-\int_0^1H(t,\gamma(t))\,dt-\int_\gamma\theta.
\]
For an exterior constant orbit $p$ put $A_{H,\theta}(p)=-H(p)$.
If $a_\gamma$ generates its orientation line, the ordinary and
relative valuations are
\[
 \nu(T^s a_\gamma)=s,\qquad
 \operatorname{val}_\theta(T^s a_\gamma)
                         =s+A_{H,\theta}(\gamma).
\]
In the action-normalized basis
$\widetilde a_\gamma=T^{-A_{H,\theta}(\gamma)}a_\gamma$ this becomes
\begin{equation}\label{eq:rc-two-valuations}
 \operatorname{val}_\theta(T^q\widetilde a_\gamma)=q,
 \qquad
 \nu(T^q\widetilde a_\gamma)=q-A_{H,\theta}(\gamma).
\end{equation}
On a sum, each valuation is the infimum of its values on the
nonzero coordinates. For a Floer or continuation cylinder $u$ from
$\gamma_-$ to $\gamma_+$, set
\[
 E_\theta(u)=E_{\mathrm{top}}(u)
       +A_{H_+,\theta}(\gamma_+)-A_{H_-,\theta}(\gamma_-).
\]
In the normalized basis the contribution of $u$ has weight
$T^{E_\theta(u)}$. By the relative-energy estimates of
\cite[Lemmas~3.11--3.12 and Corollary~3.16]{GromanBoundaryDepth},
$E_\theta(u)\geq0$. Thus the differential, admissible continuation
maps and their homotopies preserve the relative filtration as well
as the ordinary Novikov filtration.

In particular, writing $\Lambda_0=\Lambda_{\geq0}$, we obtain a
Floer complex over the Novikov ring,
\[
 CF^+_\theta(H)=\bigoplus_\gamma\Lambda_0\widetilde o_\gamma,
\]
where $\widetilde o_\gamma=T^{-A_{H,\theta}(\gamma)}o_\gamma$.
For a cofinal $S$-shaped acceleration ray $(H_i)$ let
$SC^+_{M,\theta}(K)$ be the closure of the telescope of these complexes
in the ordinary Novikov completion $SC_M^*(K;\Lambda)$. Equivalently,
it is the closed subcomplex of nonnegative relative valuation.
The comparison theorem of \cite[Lemma~4.8 and
Corollary~4.9]{GromanBoundaryDepth} identifies these lattices by
$\Lambda_0$-linear chain-homotopy equivalences preserving the relative
filtration and the ordinary filtration (as do their inverse maps
and homotopies). Here and below, in the noncompact case this
comparison uses a common admissible coefficient at infinity.

We shall also use truncations in ordinary Novikov depth. For any of
the relative lattices $C^+$ above and $a<b$, write
\[
 G^aC^+=\{x\in C^+:\nu(x)\geq a\},\qquad
 C^+_{[a,b)}=G^aC^+/G^bC^+.
\]
Write $HF^+_{\theta,[a,b)}(H)$ and $SH^+_{M,\theta,[a,b)}(K)$ for
the cohomologies of the Hamiltonian and supported complexes in this
window. The relative filtration on a window is the filtration induced
by intersection and quotient. Multiplication by $T^\lambda$ for
$\lambda\geq0$ gives a map from the window $[a,b)$ to the window
$[a+\lambda,b+\lambda)$.

In a finite ordinary window there is no difference between the
algebraic and completed telescopes. Indeed, below ordinary depth
$b$ an element of the completion has only finitely many nonzero
orbit coordinates, and discarding the other coordinates preserves
all lower bounds on relative valuation. The inclusion of the
algebraic telescope therefore gives
\begin{equation}\label{eq:rc-completed-window}
 SC^+_{M,\theta,[a,b)}(K)
 \cong\operatorname{Tel}_i CF^+_{\theta,[a,b)}(H_i),\qquad
 SH^+_{M,\theta,[a,b)}(K)
 \cong\varinjlim_i HF^+_{\theta,[a,b)}(H_i).
\end{equation}
The second isomorphism uses exactness of filtered direct limits.
These identifications and the comparison of acceleration data commute
with multiplication by $T^\lambda$ and the corresponding shifts of
windows.

The reason for introducing ordinary windows is that they bound the
energy of the trajectories. If $\nu_-$ and $\nu_+$ are the ordinary
depths of the input and output of a contribution, then
\begin{equation}\label{eq:rc-energy-depth}
              \nu_+-\nu_-=E_{\mathrm{top}}(u).
\end{equation}
In particular, a contribution surviving in $[a,b)$ satisfies
\begin{equation}\label{eq:rc-window-energy}
 0\leq E_{\mathrm{geo}}(u)\leq E_{\mathrm{top}}(u)<b-a.
\end{equation}
Thus fixing the width of the ordinary window gives the energy bound
for neck stretching in Sections~\ref{sec:collapse} and~\ref{sec:regular}.

\paragraph{Window conventions.}
In \cite{GromanBoundaryDepth} ordinary windows are indexed by
valuation, whereas in \cite{Paper1} they are indexed by
logarithmic norm. Thus, writing
\[
 G^aC=\{x:\nu(x)\geq a\},\qquad
 C_t=\{x:\|x\|<e^t\},
\]
the two conventions are respectively
\[
 C^\nu_{[a,b)}=G^aC/G^bC,\qquad
 C^{\log}_{[r,s)}=C_s/C_r.
\]
The corresponding coefficient modules are
\[
 \Lambda^\nu_{[a,b)}
 =\Lambda_{\geq a}/\Lambda_{\geq b},\qquad
 \Lambda^{\log}_{[r,s)}
 =T^{-s}\Lambda_{>0}/T^{-r}\Lambda_{>0}.
\]
Since $C_t=\bigcup_{a>-t}G^aC$, we have
\[
 C^{\log}_{[r,s)}
 =\varinjlim_{\eta\downarrow0}
 C^\nu_{[-s+\eta,-r+\eta)}.
\]
We use the valuation convention in the comparison between
ambient and intrinsic Floer theory, and the log-norm
convention in the application of the indicator-function
theory. The passage between the two is justified by
exactness of filtered direct limits.

\subsection{The locality spectral sequence and quantitative deformation}
\label{subsec:background-deformation}

Let $D\subset M$ be a compact Liouville domain, with $\theta$
defined on a neighborhood of $D$, and let $K\subset D$ be compact.
At the lowest relative energies the ambient complex recovers the
intrinsic one. When the intrinsic boundary depth is small enough,
one can also transfer the ambient differential to intrinsic
cohomology. These are Theorems~A and~B of \cite{GromanBoundaryDepth},
respectively.

\begin{theorem}[Locality and quantitative deformation]
\label{thm:deformation-package}
\begin{enumerate}[label=\textup{(\Alph*)}]
\item For the relative lattice $SC^+_{M,\theta}(K)$ there is a geometric
constant $\hbar>0$ and a quasi-isomorphism
\[
 SC^+_{M,\theta}(K)/T^\hbar SC^+_{M,\theta}(K)
 \simeq
 SC^*_{\widehat D,\theta}(K;\kk)
             \otimes_\kk\Lambda^\nu_{[0,\hbar)}.
\]
After extending scalars to $\Lambda$, the
relative lattice is dense in the ambient complex, whose ordinary-norm
completion it recovers. The reduction is natural under restriction
of supports and identifies the leading product and BV operator with
the intrinsic ones. In particular, the $T^\eta$-adic filtration of
the relative lattice, for $0<\eta\leq\hbar$ chosen compatibly with
the gap, gives the one-sided locality spectral sequence
\begin{equation}\label{eq:background-locality-ss}
 E_1^{p,q}
 =SH^{p+q}_{\widehat D,\theta}(K;\kk)
          \otimes_\kk\Lambda^\nu_{[p\eta,(p+1)\eta)},
       \qquad p\geq0.
\end{equation}
The spectral sequence need not converge to the cohomology of the
completed ambient complex.

\item For this part take $\kk=\mathbb Q$, or take a convex
cotangent-torus support as in the coefficient extension below.
Let $\beta(K)$ denote the boundary depth of the intrinsic
supported complex, namely the maximal loss of its intrinsic
filtration needed to solve boundaries. If $\beta(K)<\hbar$, there
is a bounded relative-filtered special deformation retract from the
ambient complex onto
\[
 \bigl(SH^*_{\widehat D,\theta}(K;\kk)
                     \widehat\otimes_\kk\Lambda,
                       D_{\mathrm{def}}\bigr).
\]
The transferred differential raises relative valuation by at least
$\hbar$ and has ordinary operator norm strictly smaller than
$e^{-\hbar}$. The comparison induces an isometry on cohomology.
Its inclusion, projection and homotopy preserve the relative
filtration and are bounded in the ordinary norm.

The gap admits a positive lower bound determined by the geometry of
a fixed neighborhood of $D$, independently of $K$ and of the
exterior geometry. Its local bound scales linearly under Liouville
rescaling. In particular, convex cotangent-torus domains have
intrinsic boundary depth zero.

\end{enumerate}
\end{theorem}

We index the spectral sequence by
$E_0=\operatorname{gr}_F$ and $E_1=H(E_0,d_0)$. In the
boundary-depth range of part~\textup{(B)}, collapse at $E_1$ is
equivalent to $D_{\mathrm{def}}=0$.

Theorem~B of \cite{GromanBoundaryDepth} is stated over $\mathbb Q$.
To use it over $\kk$ we need the maps of the deformation retract to
extend to the completed Novikov modules with their norm bounds.

\begin{lemma}[Bounded coefficient extension]
\label{lem:bounded-column-extension}
Let $E$ and $F$ be completed free Novikov modules with specified
orthonormal $c_0$-bases over $\Lambda_{\mathbb Q}$. Every bounded
$\Lambda_{\mathbb Q}$-linear operator $A:E\to F$ extends uniquely
to the corresponding completed free modules over $\Lambda_\kk$,
with the same operator norm. Finite collections of such operators
retain their algebraic identities and norm bounds under extension.
When the relative filtrations are specified by weights on the two
bases, the corresponding relative-filtration bounds are preserved as well. The same assertion holds
for bounded multilinear maps between completed free modules.
\end{lemma}
\begin{proof}
Write $A$ as a matrix in the specified bases. Each column lies in
$c_0$ and the norms of the columns have a common upper bound.
These conditions define an operator on finite-support vectors over
$\Lambda_\kk$ with the same norm bound, and hence a unique
continuous extension to $c_0$. Testing on basis vectors shows that
the norm is unchanged. The asserted identities hold on the dense
finite-support submodules, hence everywhere by continuity. For a
bounded multilinear map, use the same argument on tuples of basis
vectors, followed by multilinearity and continuity.

The same argument applies to a weighted basis with $\|e_i\|=w_i>0$.
A matrix defines a bounded map of norm at most $C$ precisely when
its $i$th column lies in the target space and has norm at most $Cw_i$.
Indeed, the ultrametric inequality gives this bound on finite sums,
which are dense among the sums $\sum_i a_i e_i$ with $|a_i|w_i\to0$.
\end{proof}

For the convex cotangent-torus supports used here, the source of the
retract has the telescope orbit basis and its target has the weighted
Laurent basis indexed by winding classes and a basis of fiber
cohomology. Both filtrations on these modules are given by weights
on the basis vectors. The intrinsic Laurent calculation is defined
over the ground field
\cite[Section~5]{GromanVarolgunes2022}. Thus the coefficient extension
in Lemma~\ref{lem:bounded-column-extension} identifies the two
modules over $\Lambda_\kk$ and extends the retract identities and
bounds. In particular we may use part~\textup{(B)} in
Section~\ref{sec:collapse} over our given field $\kk$.

\subsection{Local Floer cohomology of indicator functions}
\label{subsec:background-local-indicator}

Now let $K\subset M$ be a compact domain with smooth contact-type boundary. We will use
the following results of \cite[Theorems~1--4]{Paper1} to compute
restriction from a rounded polygon to its vertex caps.

To each isolated component $\gamma$ of the closed Reeb orbits of
$\partial K$, \cite{Paper1} assigns an \emph{unweighted} local Floer
cohomology group $SH^*_{\mathrm{uw}}(\gamma)$. It is computed in an
arbitrarily small isolating neighborhood of the component in loop
space.

Let $\Sigma_\tau$, $\tau\in[0,1]$, be a smooth family of contact-type
hypersurfaces. An \emph{isolating isotopy} of Reeb components is a
compact family in $LM\times[0,1]$ whose fibers $\gamma_\tau$ are
isolated path-connected components of the closed Reeb orbits of
$\Sigma_\tau$, isolated in a common neighborhood as in
\cite[Definitions~4.2 and~4.5]{Paper1}. The topology and dimension
of $\gamma_\tau$ need not remain constant. By
\cite[Theorem~1]{Paper1}, such an isotopy $\Gamma$ induces an
isomorphism
\[
 \iota_\Gamma:SH^*_{\mathrm{uw}}(\gamma_0)
 \longrightarrow SH^*_{\mathrm{uw}}(\gamma_1),
\]
which we call \emph{local Floer transport}. For short isotopies it is
computed by local continuation, and in general by subdivision and
composition. It depends only on the homotopy class of the isolating
isotopy relative to its endpoints.

We call the groups contributed by the Reeb components the \emph{Reeb
blocks}. The \emph{critical block} $H^*(K;\kk)$ comes from the constant
orbits. We refer to both kinds of summands as \emph{local blocks}.
For a positive contact-type primitive $\theta$ near $\partial K$,
we parametrize the orbits in the Hamiltonian direction, opposite to the
positive Reeb flow. Write $\mathcal T_\theta(\gamma)=-\int_\gamma\theta$
for the period and put
\[
 RH^{\leq A}(K;\theta)=H^*(K;\kk)\oplus
   \bigoplus_{\mathcal T_\theta(\gamma)\leq A}SH^*_{\mathrm{uw}}(\gamma;\kk),
\]
where the sum runs over Reeb components.

For an inclusion $K_1\subset K_2$ of smooth contact-type domains,
a source Reeb component is \emph{strictly common} if the boundaries
agree on a neighborhood of it, and \emph{non-common} if its image
is disjoint from $K_1$. Fix a contact-type primitive $\theta_1$
near the target boundary. A component which is not strictly common
is \emph{weakly common} if near the component the source boundary
can be moved monotonically inward to the target boundary through
graphs over a transverse slice to the Liouville flow of $\theta_1$,
carrying an isolating isotopy $\gamma_\tau$ to a target component
with $-\int_{\gamma_\tau}\theta_1$ independent of $\tau$.
Here the primitive is defined on the neighborhood of the deformation.
An inclusion is \emph{admissible} if every source Reeb component is
strictly common, weakly common or non-common. For the precise local
conditions and choices of isotopies see
\cite[Definitions~5.12--5.14]{Paper1}.

The period filtration is different from the Novikov filtration.
In sufficiently short Novikov windows the concave generators in the
finite-slope approximations span a subcomplex. On the quotient the
differential and continuation maps do not increase the actions of
the inner generators. In the cofinal limit these actions approach
the Reeb periods (and zero for the critical orbits), giving the
period filtration. We take cutoffs outside the period spectrum.

\begin{theorem}[Indicator-function short-window theory]
\label{thm:background-indicator-input}
Assume that the relevant Reeb components are isolated and that the
local orientation and transversality hypotheses of \cite{Paper1}
hold.
\begin{enumerate}[label=\textup{(\roman*)}]
\item There is a width $\hbar>0$, independent of the period cutoff,
such that for every ordinary Novikov window $I$ of width less than
$\hbar$, the supported Floer cohomology has an exhaustive increasing
period filtration $SH^{\leq A}_{M,I}(K)$. For each finite cutoff $A$
there is $0<\delta(A)\leq\hbar$ such that for $|I|<\delta(A)$ this
filtered part splits off as a direct summand, identified with
\[
 RH^{\leq A}(K;\theta)\otimes_\kk\Lambda^{\log}_I.
\]
\item For an admissible inclusion $K_1\subset K_2$, choose a target
period cutoff containing the target components corresponding to the
source common components below the source cutoff. In sufficiently
short windows the restriction is ordinary cohomological restriction
on the critical block, the identity on a strictly common Reeb block,
and $\iota_\Gamma\otimes\operatorname{id}$ on a weakly common block,
where $\Gamma$ is its chosen isolating isotopy. It vanishes on
non-common blocks. The primitives near the two boundaries need not
agree. This is \cite[Theorem~3]{Paper1}.
\item Consider a nested inward family of contact-type domains, with
one primitive $\theta$ positive on every boundary. For an isolating
isotopy of Reeb components, fix a finite cutoff bounding their
periods throughout the family. In a sufficiently short window the
restriction on this block is local Floer transport multiplied by
$T^{\Delta A}$, where
\[
 \Delta A=\mathcal T_\theta(\gamma_0)-\mathcal T_\theta(\gamma_1).
\]
See \cite[Theorem~4]{Paper1}.
\end{enumerate}
\end{theorem}

For a Morse--Bott family the unweighted group is the cohomology of
the orbit family with its transverse Conley--Zehnder shift and
orientation local system \cite[Theorem~3.28]{Paper1}. Along a
partially shared polygon edge the source and cap have different
orbit families. Section~\ref{subsec:rounded-cap-restriction}
constructs the isolating deformation required to compare them.

\subsection{Grading}\label{subsec:background-grading}

Let $\pi:M^{2n}\to B^n$ be a Lagrangian torus fibration,
and let $\Sigma\subset M$ denote its critical locus. Write
$M^\circ=M\setminus\Sigma$ and $V=\ker d\pi|_{M^\circ}$.

The vertical Lagrangian distribution $V$ determines a canonical
homotopy class of trivializations of $K_M^{\otimes 2}$
over $M^\circ$. If $V$ is orientable, it also determines
a canonical homotopy class of trivializations of $K_M$.
We call the latter the \emph{fibration-induced grading trivialization}.

We assume that the fibration-induced grading trivialization
is defined and extends, up to homotopy, over $M$.
The following lemma gives a sufficient condition.

\begin{lemma}[Extension of the grading trivialization]
Suppose that $c_1(M)=0$ and that the critical locus $\Sigma$
is a closed Whitney-stratified subset of real codimension
at least three.

Then the canonical homotopy class of trivializations of
$K_M^{\otimes 2}$ over $M^\circ$ extends uniquely to a
homotopy class of trivializations over $M$.

If $V$ is orientable (in particular, if the regular base
is orientable), the same holds for the canonical homotopy
class of trivializations of $K_M$.
\end{lemma}

\begin{proof}
For a compatible $J$, an orthonormal frame of the vertical Lagrangian
distribution determines the unit complex volume form taking value
one on that frame. A change of frame changes this form only by the
orientation sign, so its square is canonical. If $V$ is orientable,
an orientation defines the form itself; reversing it multiplies the
form by $-1$ and leaves its homotopy class unchanged. These classes
are independent of $J$.

General position moves loops and homotopies between them off the
Whitney-stratified set $\Sigma$, since its real codimension is at least
three. Thus restriction induces an isomorphism
\[
 H^1(M;\mathbb Z)\xrightarrow{\;\cong\;}H^1(M^\circ;\mathbb Z).
\]
Since $c_1(M)=0$, both $K_M$ and its square are trivial. Their
homotopy classes of trivializations form torsors over $H^1$, so
restriction is a bijection on these classes. The fibration-induced
classes therefore extend uniquely up to homotopy.
\end{proof}

\section{Collapse of the locality spectral sequence}
\label{sec:collapse}

Let $(M^4,\omega)$ be closed or geometrically of finite type, with
$c_1(M)=0$, and fix the fibration-induced grading trivialization of
$K_M$ from 
Section~\ref{subsec:background-grading}. Let $L\subset M$ be a
Maslov-zero Lagrangian torus with respect to this grading, and
choose a Weinstein domain $D=L\times P$ for $P$ sufficiently
small, smooth, and strictly convex. Write $V(D)$ for a slightly
larger Weinstein neighborhood and $\theta$ for its action--angle
primitive.

The aim of this section is to prove the following theorem.

\begin{theorem}[Regular collapse]\label{thm:dimension-four-collapse}
The locality spectral sequence of a sufficiently small convex
neighborhood of $L$ collapses at $E_1$.
\end{theorem}

For $P$ convex, the intrinsic boundary depth of $D$ is zero by a standard computation. See \cite[Cor. 5.7]{GromanBoundaryDepth}. By
Theorem~\ref{thm:deformation-package}\textup{(B)}, collapse therefore
amounts to the vanishing of the transferred differential
$D_{\mathrm{def}}$. Write the deformation retract as
\begin{equation}\label{eq:rc-retract}
 (C,d)\ \underset{\iota}{\overset{p}{\rightleftarrows}}\
 (\mathsf H,D_{\mathrm{def}}),\qquad
 \mathsf H=SH^*_{\widehat D,\theta}(D;\kk)
                    \widehat\otimes_{\kk}\Lambda,
\end{equation}
where $C$ is a completed ambient telescope. The maps $p,\iota$ are
bounded in the ordinary norm, preserve the relative filtration, and
satisfy $p\iota=\operatorname{id}$, $d\iota=\iota D_{\mathrm{def}}$,
and $pd=D_{\mathrm{def}}p$. The transferred differential raises
relative valuation by at least $\hbar>0$.

If $D_{\mathrm{def}}$ is nonzero, multiplication by some $T^\lambda$,
$\lambda>0$, kills a nonzero class in the cohomology of the relative
lattice. We will show that this failure of injectivity is visible in a finite ordinary
window for a single Hamiltonian. We will reach a contradiction by showing that for appropriately chosen
Floer data neck stretching confines every contributing
cylinder to $V(D)$. There its relative energy vanishes, and
multiplication by $T^\lambda$ into the shifted window is injective.

\subsection{Strictly \texorpdfstring{$S$}{S}-shaped Hamiltonians}

We need the Floer equation on the exterior of the bend to be
holomorphic. Choose a protected cutting annulus
$A=\{\eta'\leq\rho\leq\eta\}$ as in
Section~\ref{subsec:background-filtrations}. We call an $S$-shaped
Hamiltonian \emph{strictly $S$-shaped} if its positive radial slope
tapers to zero beyond $A$, before a level $\delta>\eta$, and it is
spatially constant on the ensuing exterior plateau. Take the taper
with slope below the least Reeb period, so it has no nonconstant
one-periodic orbits. In the compact case the plateau extends over
$M\setminus D_\delta$. In the noncompact case it ends in the tail
$H=C(H)+q_R$, where $q_R=\chi_R\circ f_0$, $R$ is above all
critical values of $f_0$, and $\chi_R$ vanishes on $(-\infty,R]$
and is eventually linear with a small positive slope.

The local generators are the orbits on the inner side of $A$.
By \emph{regular strict data} we mean that these orbits are
nondegenerate and that the moduli spaces with local asymptotics
used in the window are regular. The constant orbits on the plateau
are exterior, including those in $V(D)$. The distinction between
the cutting level $\eta$ and the outer bound $\delta$ allows us
to perturb these constant orbits while keeping the cutting annulus
fixed.

\begin{lemma}[Wide exterior plateaux]
\label{lem:finite-type-flat-exterior}
In the geometrically finite-type case, strictly $S$-shaped data can
be chosen with arbitrarily wide flat exterior plateaux. Given
$E>0$, the plateau can be chosen so that every Floer or admissible continuation cylinder
with a local asymptotic and geometric energy at most $E$ stays in
a fixed compact set before the exhaustion tail. This set is
independent of further enlargement of the plateau, and the choice
is uniform for compact families with fixed exterior data.
\end{lemma}

\begin{proof}
Choose $\chi_R$ with positive derivative on $(R,\infty)$ and with
first and second derivatives sufficiently small. The bounds on
$X_{f_0}$ and $\nabla X_{f_0}$ give the small-period and dissipativity
bounds for $q_R$; the tail has no one-periodic orbits where $\chi_R'>0$.
The eventual slope is the prescribed sufficiently small exhaustion
coefficient. Keep the normalized exterior profile fixed along a ray,
as in \cite[Proposition~3.9]{GromanBoundaryDepth}.

Let $W_R$ be the distance from a fixed compact neighborhood containing
the cutting collar and the part of the exterior where $J$ differs from
$J_0$ to the beginning of the exhaustion tail. Since $\|df_0\|$ is
uniformly bounded and $f_0$ is proper, $W_R\to\infty$. On the intervening
flat region the equation is $J_0$-holomorphic. By geometric boundedness
there are constants $r,c>0$ such that crossing a ball of radius $r$
centered on the curve, with no boundary inside the ball, costs at least
$cr^2$ of area. A curve reaching the tail crosses this region and hence
has energy at least
\[
 cr^2\left\lfloor\frac{W_R-2r}{3r}\right\rfloor.
\]
Choose a compact neighborhood wide enough that crossing it already
costs more than $E$, and put the tail beyond it. The same neighborhood
works when the plateau is enlarged. The argument also applies to
continuation on a common plateau, where spatially constant additions
to the Hamiltonian leave the equation unchanged.
\end{proof}

\subsection{Floer theory in a finite ordinary window}

The flat exterior consists of degenerate constant orbits
\cite[Remark~3.6]{GromanBoundaryDepth}. To be able to disregard these we will show that raising the plateau puts the exterior
generators beyond the upper endpoint of the window and that an energy
estimate excludes breaking through them. For a window $[a,b)$ choose
the Hamiltonian high enough on the cutting annulus that
\begin{equation}\label{eq:rc-high-cut}
 H|_{\rho=r}-\sigma>b\qquad(\eta'\leq r\leq\eta).
\end{equation}
This can be achieved by increasing the plateau height: the difference
between the cutting values and the plateau height is fixed along
an acceleration ray. Every exterior constant $p$ has ordinary Novikov valuation
at least $H(p)$ in the relative lattice, hence is absent below $b$.
We use the normalized local orbit lines to form
$C^{+,\mathrm{str}}_\theta(H)=\bigoplus_{\gamma\ \mathrm{local}}
\Lambda_0\widetilde o_\gamma$ and its ordinary window
$C^{\mathrm{str}}_{[a,b)}(H)$ as in
Section~\ref{subsec:background-filtrations}. We write
$HF^{+,\mathrm{str}}_{\theta,[a,b)}(H)$ for its Floer cohomology.

We use the following estimate both to exclude breaking through the
omitted orbits and later in the neck-stretching argument.

\begin{lemma}[The two ends lie in one interior component]
\label{lem:rc-cut-component}
Let $u$ be a Floer cylinder or an admissible monotone continuation
cylinder with both asymptotic orbits in the interior of $D_r$, where
$r$ is a regular cutting level for $u$. Assume that near $\partial D_r$
the Hamiltonian has the form $H_s=h_s(e^\rho)$, $J$ is of contact type, and
$\theta(X_{H_s})=-\sigma_r$ on $\partial D_r$, with $\sigma_r$ independent
of $s$. If the input is $T^q\widetilde a_{\gamma_-}$, $q\geq0$, the
output has ordinary valuation less than $b$, and
\[
 H_-|_{\rho=r}-\sigma_r>b,
\]
then both ends belong to the same connected component of $u^{-1}(D_r)$.
\end{lemma}

\begin{proof}
Suppose the component $S$ containing the negative end does not contain
the positive end. Write $c(s)=H_s|_{\rho=r}$ and orient the cut as the
boundary of $S$. Since $\omega=d\theta$ on $D_r$, the topological
energy is the integral of $d(u^*\theta+H_s(u)\,dt)$ over $S$.
Stokes' theorem therefore gives
\[
 E_{\mathrm{top}}(u|_S)
 =A_{H_-,\theta}(\gamma_-)
   +\int_{\mathrm{cut}}(u^*\theta+c(s)\,dt).
\]
Here and below one truncates the negative end before applying Stokes
and then takes the limit.

Let $n$ be an outward normal to $S$ and $\tau=jn$ the oriented
boundary tangent. The Floer equation gives
\[
 du(\tau)-X_{H_s}dt(\tau)
 =J\bigl(du(n)-X_{H_s}dt(n)\bigr).
\]
On the cut, $\theta\circ J$ is a positive multiple of $d\rho$.
Since $d\rho(X_{H_s})=0$ and $d(\rho\circ u)(n)\geq0$, applying
$\theta$ gives $u^*\theta(\tau)+\sigma_r dt(\tau)\geq0$. Thus
\[
 \begin{aligned}
 E_{\mathrm{top}}(u|_S)
 &\geq A_{H_-,\theta}(\gamma_-)
       +\int_{\mathrm{cut}}(c(s)-\sigma_r)\,dt\\
 &= A_{H_-,\theta}(\gamma_-)+c(-\infty)-\sigma_r
       +\int_S c'(s)\,ds\wedge dt\\
 &\geq A_{H_-,\theta}(\gamma_-)+c(-\infty)-\sigma_r.
 \end{aligned}
\]
The equality is Stokes' theorem and the last inequality uses
monotonicity. Topological energy is nonnegative on the complement of
$S$, so this is also a lower bound for $E_{\mathrm{top}}(u)$.
On the other hand, the input has ordinary valuation
$q-A_{H_-,\theta}(\gamma_-)$, so by \eqref{eq:rc-energy-depth},
\[
 E_{\mathrm{top}}(u)<b+A_{H_-,\theta}(\gamma_-)-q.
\]
Combining the bounds gives $c(-\infty)-\sigma_r<b-q\leq b$,
contrary to the hypothesis.
\end{proof}

On the protected annulus we have $\sigma_r=\sigma$, so
\eqref{eq:rc-high-cut} gives the hypothesis of the lemma. Breaking
through an exterior constant would produce slices converging to that
constant and hence lying entirely outside a cutting level strictly
inside the annulus (chosen regular for the sequence of cylinders).
Such a slice separates the two ends, contradicting
Lemma~\ref{lem:rc-cut-component}. The argument is uniform for compact
families with a common cutting annulus and the same bound.

In the noncompact case, the $C^0$ bound needed for compactness comes
from the dissipative exhaustion data, or, for strict cylinders on a
sufficiently wide common plateau, from
Lemma~\ref{lem:finite-type-flat-exterior}.

Consequently the ordinary compactification of the zero- and
one-dimensional moduli spaces with local asymptotics suffices to
define the differential, admissible continuation maps and their
homotopies in the window. To see that the counts preserve the
relative filtration, apply the cut-domain calculation of
\cite[Lemmas~3.11--3.12]{GromanBoundaryDepth} to a level in the
protected annulus. With both ends local there are no exterior
end terms, and \eqref{eq:rc-normalized-monotonicity} gives
$E_\theta(u)\geq E_{\mathrm{geo}}(u|_{u^{-1}(M\setminus D_r)})\geq0$.
In particular, a trajectory of relative energy zero stays on the
local side of the cut. Indeed, on the exterior its $s$-derivative
would vanish, which is incompatible with entering and leaving that
region. Thus the strict counts define a theory over $\Lambda_0$, with
continuation identities and compatibility with Novikov shifts.
All intermediate terms remain in the window because ordinary depth
cannot decrease and relative valuation stays nonnegative. No
nondegeneracy of the omitted exterior orbits is needed.

We now compare the strict window complexes with those of ordinary
$S$-shaped data. This comparison allows us to use strict Hamiltonians
without changing the supported theory.

\begin{lemma}[Comparison with regular $S$-shaped data]
\label{lem:rc-regular-companion}
Fix an ordinary window $[a,b)$ and regular strictly $S$-shaped data
$(H,J)$ satisfying \eqref{eq:rc-high-cut}. There is a nondegenerate
$S$-shaped Hamiltonian $H^{\mathrm{reg}}$, agreeing with $H$ on
$D_\eta$, for which monotone continuation induces an isomorphism
of complexes
\[
 \Phi_H:C^{\mathrm{str}}_{[a,b)}(H)
 \xrightarrow{\ \sim\ }
 CF^+_{\theta,[a,b)}(H^{\mathrm{reg}}).
\]
Both $\Phi_H$ and its inverse preserve the ordinary and relative
filtrations. The comparison commutes with the ordinary window maps and with
multiplication by Novikov monomials whenever the same data satisfy
the cutting condition for all the windows involved.

For a cofinal strict acceleration ray $(H_i)$ with a common cutting
annulus and fixed normalized exterior profile, these isomorphisms
can be chosen to commute on cohomology with structural continuation.
In particular, for a ray computing the theory with support in $D$,
\begin{equation}\label{eq:rc-strict-supported-window}
 SH^+_{M,\theta,[a,b)}(D)
   \cong\varinjlim_i HF^{+,\mathrm{str}}_{\theta,[a,b)}(H_i).
\end{equation}
The same holds for a compact support $K\subset D$ using its indicator
Hamiltonians. In the finite-type case the comparison of supported
windows can also be made between sufficiently small exhaustion
coefficients by using strict data with a sufficiently wide common
flat plateau.
\end{lemma}

\begin{proof}
Choose $H^{\mathrm{reg}}=H+h$, where $h\geq0$ is a sufficiently
small autonomous perturbation supported outside $D_\eta$, making
the Hamiltonian Morse near the entire flat critical set, including
its boundary. Choose it to preserve the smallness, end-value and
exhaustion conditions. The local generators and their actions are
unchanged, and the exterior constants of both Hamiltonians are
above the window. Thus the two window complexes have the same
underlying module.

For the continuation use $H_s=H+\beta(s)h$ with $\beta$ increasing
from zero to one. Keep the almost complex structure independent of
$s$ on $D_\eta$, and vary it generically on the exterior, retaining
the geometric bounds at infinity. A trajectory contained in $D_\eta$
satisfies the regular stationary equation. Any other trajectory has
a point outside $D_\eta$ at which $\partial_su\ne0$, so
domain-dependent variations there give transversality. Sphere
attachment is excluded in dimensions zero and one by the marked-point
transversality argument of
\cite[Section~3.3]{AbouzaidGromanVarolgunes24},
using exactness on the fixed local region. Since $h$
vanishes on the cutting annulus, both ordinary monotonicity and
\eqref{eq:rc-normalized-monotonicity} hold. Stretching $\beta$ in
the $s$-direction gives the required derivative bounds.
Lemma~\ref{lem:rc-cut-component} excludes breaking through the
omitted orbits, as above, so continuation defines a chain map
preserving both filtrations.

At relative energy zero, the preceding cut-domain energy estimate
confines the trajectory to the local side of the cut. There the
continuation equation is independent of $s$. Every nonconstant
solution then has an $\mathbb R$-translation symmetry and cannot
be rigid for the fixed continuation equation. The constant cylinders
induce the identity on the local orientation lines. Consequently,
in the action-normalized basis,
\[
                  \Phi_H=\id+P,
\]
where every nonzero term of $P$ has positive relative exponent.
There are only finitely many local generators at this fixed
Hamiltonian. By the Novikov condition there is therefore an
$\epsilon>0$ such that every term of $P$ surviving in the window
raises the relative exponent by at least $\epsilon$. On the other
hand, a term $T^q\widetilde a_\gamma$ in the window has
\[
                 0\leq q<b+A_{H,\theta}(\gamma).
\]
Thus $P^N=0$ in the window as soon as
$N\epsilon>b+\max_\gamma A_{H,\theta}(\gamma)$ (the zero module
case is immediate). The inverse of the continuation map is the
finite sum
\[
                 \Phi_H^{-1}=\sum_{j=0}^{N-1}(-P)^j.
\]
It is a chain map and preserves both filtrations, since $P$ does.
The counts defining $\Phi_H$ are independent of the endpoints of
the ordinary window; with common data the maps therefore commute
with truncations, inclusions of lower bounds and multiplication
by Novikov monomials into shifted windows.

For the assertion about a ray, add the same sufficiently small
exterior perturbation at every level. The resulting regular
companions form a cofinal $S$-shaped acceleration ray as in
\cite[Proposition~3.9]{GromanBoundaryDepth}. Between consecutive
levels there are two ways of continuing from the strict Hamiltonian
at the earlier level to the regular one at the later level.
The two concatenations are joined, after reparametrization, by
convex interpolation of the Hamiltonians. This preserves ordinary
and normalized exterior monotonicity and the common cutting
annulus; stretching the parameters gives the derivative bounds.
The resulting continuation homotopy is defined in the window by
the same compactness argument as above. Hence the continuation
square commutes on cohomology. Taking direct limits and using
\eqref{eq:rc-completed-window} proves
\eqref{eq:rc-strict-supported-window}. Nothing in this comparison
depends on the compact support except the cofinality of the
chosen Hamiltonians for its indicator.

In the noncompact case we keep the exhaustion coefficient fixed
along each ray. To compare two sufficiently small coefficients,
choose strict rays with the respective tails but identical data
up to the end of a common flat plateau. For contributions in
$[a,b)$, the geometric energy is bounded by $b-a$. By
Lemma~\ref{lem:finite-type-flat-exterior}, with the plateau wide
enough for both rays, the differential and structural continuation
trajectories cannot reach either tail. Thus the strict window
complexes and continuation maps agree. Applying the comparison
with regular companions to both rays identifies the supported
windows for the two coefficients. The argument also works for a
finite collection of shifted windows by making the plateau wide
and high enough for all of them. One may then use the regular
comparison theorem at the smaller common coefficient.
\end{proof}

\subsection{Detection of the deformation}
We now test whether $D_{\mathrm{def}}$ vanishes. The window in
which multiplication by $T^\lambda$ detects a nonzero differential
must be fixed before we choose the exterior data for neck stretching.
\begin{proposition}[Finite-window detection]\label{prop:rc-strict-detection}
If $D_{\mathrm{def}}\ne0$, there are $a<b$ and $\lambda>0$ such
that for a regular strictly $S$-shaped Hamiltonian $H$ the map
\[
 T^\lambda:HF^{+,\mathrm{str}}_{\theta,[a,b)}(H)
 \longrightarrow HF^{+,\mathrm{str}}_{\theta,[a+\lambda,b+\lambda)}(H)
\]
is not injective. The cutting collar may be chosen arbitrarily small,
the plateau arbitrarily high, and a sphere-free admissible exterior
almost complex structure prescribed, after $a,b,\lambda$ are fixed.
This is a choice of auxiliary data, not an additional hypothesis on
$M$. In the noncompact case the flat plateau may also be required to
be arbitrarily wide.
\end{proposition}

\begin{proof}
Let $\mathsf H^+$ denote the nonnegative relative lattice in
$\mathsf H$. Since finite Laurent sums are dense and the transferred
differential is bounded, choose a normalized Laurent generator
$x\in\mathsf H^+$ with $D_{\mathrm{def}}x\ne0$,
put $\lambda=\operatorname{val}_\theta(D_{\mathrm{def}}x)$ and
$y=T^{-\lambda}D_{\mathrm{def}}x$. Then $D_{\mathrm{def}}y=0$ and
$\operatorname{val}_\theta(y)=0$, whereas every boundary in
$\mathsf H^+$ has relative valuation at least $\hbar>0$. Thus $[y]$
is a nonzero torsion class killed by $T^\lambda$. In $C^+$ put
$v=\iota y$, $z=\iota x$, so that $dv=0$ and $dz=T^\lambda v$.

This torsion is detected in two finite ordinary windows. Choose $a$
so that $v\in G^aC^+$ and $z\in G^{a+\lambda}C^+$. Choose $b$ large
enough that $v$ remains nonzero in cohomology modulo $G^bC^+$:
indeed, choose a monomial in a Laurent-basis coefficient of $y$
with relative exponent less than $\hbar$ and ordinary depth $w$. If
$\nu(pu)\geq\nu(u)-L$, take $b>\max\{a,w+L\}$. The equality
$v=du+e$ with $u\in G^aC^+$ and $e\in G^bC^+$ would imply
$y=D_{\mathrm{def}}pu+pe$, which is impossible at the chosen
coefficient and depth. Hence multiplication by $T^\lambda$ on
$SH^+_{M,\theta,[a,b)}(D)$ has a nonzero kernel in the shifted
window, since $dz=T^\lambda v$.

Now choose a cofinal strict ray with the prescribed exterior geometry
and cutting values above $b+\lambda+\sigma$. By
Lemma~\ref{lem:rc-regular-companion} and exactness of filtered direct
limits, the noninjectivity
is visible at some finite stage. More explicitly, a nonzero class in
the kernel is represented at some stage; its image becomes zero at a
later stage, while the class itself remains nonzero there. The
Hamiltonian at that stage has the claimed properties.
\end{proof}
\subsection{Neck stretching for strict data}

Fix the windows supplied by Proposition~\ref{prop:rc-strict-detection}.
If cylinders of bounded energy escape $V(D)$ as the outer bend
approaches $\partial D$, their exterior portions have a holomorphic
limit. The argument below extracts from it a plane with one negative
Reeb end. Dimension four enters in excluding that plane.

We choose the convex set $P$ with everywhere positive curvature. Put
$Y=\partial D$ and $\alpha=\theta|_Y$, so that
$\theta=e^\rho\alpha$ on the exterior collar. Regard $W=M\setminus D$
as a manifold with a negative cylindrical end, with coordinate
$\ell=\chi(\rho)\to-\infty$ as $\rho\downarrow0$. Fix an exterior
$\omega$-compatible almost complex structure $J$ which is translation
invariant in $\ell$ near this end, satisfies $J\partial_\ell=R_\alpha$,
and restricts to a $d\alpha$-compatible structure on $\ker\alpha$. Extend $J$
smoothly across each finite neck, with all modifications supported
in a collar shrinking to $\partial D$ together with the bend. Thus
the extensions agree with $J$ on every fixed compact subset of $W$
once the bend is sufficiently small.
In the geometrically finite-type case we retain $J_0$ outside a compact
set. The wide plateaux of Lemma~\ref{lem:finite-type-flat-exterior}
then keep the bounded-energy cylinders considered here away from the
exhaustion tail, where
the Hamiltonian would cease to be spatially constant.

\begin{lemma}[Exterior planes]\label{cor:no-exterior-planes}
A generic cylindrical almost complex structure on $W$ admits no
nonconstant finite-energy holomorphic plane with one negative Reeb end.
\end{lemma}

\begin{proof}
Each Reeb component of the cotangent-torus boundary is $N=T^2$, with
unparametrized orbit space $S=N/S^1=S^1$. Positive curvature gives a
positive shear in the contact linearization, whose Robbin--Salamon
index in the cotangent trivialization is $\mu_{\mathrm{RS}}=1/2$.
An iterate multiplies the shear coefficient by a positive integer
and has the same index. The Maslov-zero condition identifies the
induced determinant trivialization with the one fixed on $M$.

For free Morse--Bott asymptotics in dimension four, the index is
\[
 \operatorname{ind}(u)
 =-\chi(\dot\Sigma)+2c_1^\tau(u^*TW)
  +\sum_+\left(\mu_{\mathrm{RS}}+\tfrac12\dim S\right)
  -\sum_-\left(\mu_{\mathrm{RS}}-\tfrac12\dim S\right).
\]
For a plane with one negative end, the relative Chern number is zero,
so this gives
\[
 \operatorname{ind}(u)=-1-(\tfrac12-\tfrac12)=-1.
\]
A plane cannot lie entirely in the negative cylindrical end, by the
maximum principle for $\ell$. It therefore crosses a compact collar
where the exterior structure is free to vary. The simple-curve
transversality theorem, applied there, excludes simple planes for
generic $J$ even with $J_0$ fixed at infinity. A multiply covered plane factors
through a simple punctured curve. The covering extends over the
punctures; its target has genus zero and at most one puncture. Since
the negative end is nonremovable, that target is again a plane.
This excludes multiple covers as well. We use here the Fredholm,
genericity and factorization results of
\cite[Sections~8.3.2--8.3.3 and Remark~8.52]{WendlLowDimensions}.
\end{proof}
We call an exterior almost complex structure $J$ with the property of
Lemma~\ref{cor:no-exterior-planes} \emph{plane-free}. We may choose
the exterior structure to be sphere-free as well as plane-free.
Indeed, a simple holomorphic sphere has index $-2$, and a
nonconstant sphere cannot be contained in the exact cylindrical end.
In the noncompact case we first choose $J_0$ generically, preserving
the geometric bounds at infinity.

The prescribed exterior does not obstruct regularity of the strict
differential. A nonstationary cylinder with local asymptotics has a
regular point on the inner side of the cutting annulus, by the density
of regular points \cite[Theorem~4.3]{FHS95}. Variations of $J_t$
supported there suffice for the usual universal transversality
argument. Stationary cylinders are regular by nondegeneracy of the
local orbits. The generic one-parameter family $J_t$ can also be
chosen sphere-free, since the parametrized index of a simple sphere
is $-1$. We use the same $J$ as we shrink the bend.
\begin{proposition}[Confinement in a finite window]
  \label{prop:rc-confinement}
  Fix an ordinary window $[a,b)$ and a plane-free cylindrical
  almost complex structure $J$ on $W=M\setminus D$.
  
  There exist $\delta_0>0$ and $C_0>0$ with the following
  property. Let $0<\delta<\delta_0$ and let $(H,J_\delta)$
  be regular strictly $S$-shaped Floer data such that:
  \begin{enumerate}
  \item The outer bend of $H$ lies inside $D_\delta$, and
        its exterior plateau has height greater than $C_0$.
  \item $J_\delta$ is a smooth compatible almost complex
        structure on $M$, possibly time-dependent on $D_\delta$,
        which agrees with $J$ outside $D_\delta$.
  \end{enumerate}
  In the geometrically finite-type case, require in addition
  that the flat exterior plateau be sufficiently wide.
  
  Then every Floer cylinder contributing to the differential
  of $C^{\mathrm{str}}_{[a,b)}(H)$ is contained in $V(D)$.
  \end{proposition}

  To extract a forbidden plane from an exterior limit, we need to
  control its conformal type and the topology of its compact
  truncations. The first point is a contact case of McLean's
  finite-type result. For a noncompact surface, genus zero means that
  every compact subsurface is planar.
  
  \begin{lemma}[Finite type, after McLean]
  \label{lem:rc-exterior-finite-type}
  Let $(W,\omega)$ be a symplectic manifold with a negative contact end
  $(-\infty,0]\times(Y,\alpha)$, where $Y$ is compact, and let $J$ be
  an $\omega$-compatible almost complex structure which is cylindrical
  on this end. Let $v:\Sigma\to W$ be a proper nonconstant
  $J$-holomorphic map from a smooth, connected, noncompact genus-zero
  Riemann surface without boundary. Assume that its image is contained
  in the union of a compact set and the negative end.
  
  Write $\ell$ for the end coordinate and
  $\Sigma_-=v^{-1}((-\infty,0]\times Y)$. Suppose that
  \begin{equation}\label{eq:rc-end-energy}
   E^-(v):=
   \int_{\Sigma_-}v^*d\alpha
   +\sup_{\substack{\psi\in C_c^\infty(\mathbb R),\ \psi\geq0\\
                     \int_{\mathbb R}\psi=1}}
         \int_{\Sigma_-}v^*\bigl(\psi(\ell)\,d\ell\wedge\alpha\bigr)
   <\infty.
  \end{equation}
  Then $\Sigma$ is biholomorphic to a sphere minus finitely many points.
  If the Reeb flow is Morse--Bott, every puncture is negative and
  asymptotic to a nonconstant closed Reeb orbit.
  \end{lemma}
  
  \begin{proof}
  The contact form $\alpha$ determines the stable Hamiltonian structure
  $(d\alpha,\alpha)$, and \eqref{eq:rc-end-energy} is its negative-end
  energy. Properness and the compact-containment assumption imply that
  the part of the domain above any fixed end level is compact, so the
  interior energy is finite as well. The finite-type assertion follows
  from \cite[Lemma~8.8]{McLeanDiscrepancy}, whose cobordism is allowed
  to be noncompact.
  
  Near each puncture, properness forces the image to leave every compact
  subset of $W$. The compact-containment assumption therefore implies
  that $\ell\circ v$ tends to $-\infty$. Thus no puncture is removable.
  When the Reeb flow is Morse--Bott,
  \cite[Proposition~5.6]{BEHWZ} gives a Reeb asymptote at every puncture,
  and all these asymptotes are negative.
  \end{proof}
  
  We next record how exhaustive convergence behaves under a compact
  transverse cut. The point is that the cut retains the entire inverse
  image of the target region, including all its connected components.
  
  \begin{lemma}[Transverse cuts under exhaustive convergence]
  \label{lem:rc-transverse-cuts}
  Let $u_k:\Sigma_k\to W_k$ be proper pseudoholomorphic curves
  converging in the exhaustive Gromov sense to a proper nodal curve
  $u_\infty:\Sigma_\infty\to W$. Here the open target regions $W_k$
  exhaust $W$, and the domains $\Sigma_k$ are smooth surfaces without
  boundary; they need not be connected.
  
  Let $K\subset W$ be a compact domain with smooth boundary. Suppose
  that $u_\infty$ is transverse to $\partial K$ on its normalization
  and that no node maps to $\partial K$. Then, for all sufficiently
  large $k$, the map $u_k$ is transverse to $\partial K$, and
  $u_k^{-1}(K)$ is diffeomorphic to a smoothing of the compact nodal
  surface $u_\infty^{-1}(K)$. Here smoothing replaces each node by an
  annulus joining its two branches. These identifications give a
  correspondence of connected components and of boundary circles.
  Near the boundary, they can be chosen so that the maps converge
  smoothly to $u_\infty$.
  \end{lemma}
  
  \begin{proof}
  Choose a compact neighborhood $K^+$ of $K$ with $K$ in its interior.
  Properness of the limit and the exhaustion in
  \cite[Definition~2.21]{FishHofer} allow us to choose a compact
  limiting domain piece containing the full inverse image of $K^+$
  in its interior. By increasing the exhaustion level, the corresponding
  approximating pieces also contain the full inverse image of $K^+$
  for all sufficiently large $k$.
  
  On these compact pieces, Gromov convergence gives diffeomorphisms
  from a smoothing of the limiting nodal domain to the approximating
  domains. After pullback, the maps converge uniformly on the whole
  smoothed domain and smoothly away from the circles introduced by
  smoothing the nodes; see \cite[Definition~2.16]{FishHofer}.
  
  The limiting inverse image of $\partial K$ is a finite union of
  smooth circles disjoint from the nodes. Choose disjoint collars of
  these circles. A defining function for $\partial K$, composed with
  $u_\infty$, has nonzero normal derivative along each circle. Smooth
  convergence therefore makes the corresponding zero set for $u_k$
  a unique nearby normal graph. A small isotopy in the collars
  identifies these circles and their interior sides. On the complement
  of the collars, the limiting map stays a positive distance from
  $\partial K$, so uniform convergence excludes any further
  intersections, including in the node-smoothing regions. The
  identifications consequently restrict to the claimed diffeomorphisms
  of the cut domains, with smooth convergence near their boundaries.
  \end{proof}
  
  For a cylinder the plateau estimate will put both ends on the same
  side of the contact cut, in a single component. We give the exterior
  argument for a punctured sphere, since in Section~\ref{sec:regular}
  the same conclusion about all the punctures will allow us to treat
  products and continuation maps as well.
  
  \begin{lemma}\label{lem:genus-zero-exterior-confinement}
  Let $J$ be a plane-free exterior almost complex structure as above,
  and let $u_k$ be a sequence of Floer solutions on spheres with
  finitely many cylindrical ends, at least one on each sphere. Suppose that,
  for $\delta_k\downarrow0$, all asymptotic orbits lie inside
  $D_{\delta_k}$ and the restrictions over $M\setminus D_{\delta_k}$
  are $J$-holomorphic, with $\omega$-area bounded by $E$. Suppose also
  that these exterior images lie in the union of a fixed compact
  subset of $W$ and its negative end.
  
  Assume that for every fixed sufficiently small collar level $r>0$,
  all the ends lie in one component of $u_k^{-1}(D_r)$ whenever
  $k$ is sufficiently large and the cut is transverse to $u_k$.
  Then the images of $u_k$ lie in $V(D)$ for all sufficiently large $k$.
  \end{lemma}
  
  \begin{proof}
  Suppose otherwise and pass to a subsequence meeting
  $M\setminus V(D)$ at every stage. We first bound the cylindrical
  energy of the exterior curves. Orient the level circles as
  boundaries of the interior sublevels and write
  \[
   \mathcal A_k(r)=\int_{u_k^{-1}(\rho=r)}u_k^*\alpha.
  \]
  For a holomorphic curve in the cylindrical collar,
  $u_k^*\alpha=-d(\ell\circ u_k)\circ j$ and $u_k^*d\alpha\geq0$.
  Hence the level flux is nonnegative, and Stokes' theorem gives, for
  regular levels $\delta_k<s<t$ in the collar,
  \begin{equation}\label{eq:rc-flux}
   \begin{gathered}
   0\leq\mathcal A_k(s)\leq\mathcal A_k(t),\\
   \int_{u_k^{-1}([s,t])}u_k^*\omega
   =e^t\mathcal A_k(t)-e^s\mathcal A_k(s)
   \geq(e^t-e^s)\mathcal A_k(s).
   \end{gathered}
  \end{equation}
  Fix two collar levels $0<r_0<r_1$ which are regular for all the $u_k$,
  and set $B=E/(e^{r_1}-e^{r_0})$. Such a common choice exists by Sard's
  theorem, since the sequence is countable. For every fixed
  $0<\epsilon<r_0$, we have $\delta_k<\epsilon$ for all sufficiently
  large $k$. Hence $\mathcal A_k(s)\leq B$ at every regular level
  $\epsilon\leq s\leq r_0$. On any such truncated end, Stokes gives
  \[
   \int u_k^*d\alpha\leq B,
   \qquad
   \int u_k^*d\bigl(\varphi(\ell)\alpha\bigr)\leq B
   \quad(0\leq\varphi\leq1,\ \varphi'\geq0).
  \]
  Both terms in
  $d(\varphi(\ell)\alpha)=\varphi(\ell)d\alpha+
  \varphi'(\ell)d\ell\wedge\alpha$ are nonnegative on the curve.
  The second inequality therefore bounds the longitudinal energy in
  \eqref{eq:rc-end-energy}, by taking $\varphi$ to be a primitive of
  the test function $\psi$.
  
  Apply exhaustive Gromov compactness
  \cite[Theorem~1]{FishHofer} to the open exterior restrictions of
  $u_k$. More precisely, choose an exhaustion whose inner boundaries
  approach $\partial D$, and pass to a subsequence so that each
  $D_{\delta_k}$ lies inside the corresponding inner boundary.
  Restriction to these open target regions gives proper holomorphic
  maps without boundary: all asymptotics lie on the interior side,
  so the preimage of a compact exterior set stays away from the
  punctures and is compact. Every component of a restriction is
  nonconstant, hence stable. The domains have genus zero and no marked
  or nodal points. On each compact target region the compatible
  metric area is bounded by the exterior area bound. These are the
  local bounds required by the theorem. The limit $u_\infty$ is
  nonempty because the escaping curves meet a fixed compact region
  outside $V(D)$ and all preimages there are retained.
  
  The limit and its normalization are proper. Every normalization
  component has genus zero, and the incidence graph of the nodal
  curve is a forest. Indeed, a compact subsurface of positive genus,
  or a finite cycle together with paths joining its nodal points in
  the normalization, is contained in the inverse image of a compact
  target domain with transverse boundary avoiding the nodal images.
  By Lemma~\ref{lem:rc-transverse-cuts}, smoothing this cut domain
  gives a subsurface of an approximating domain. The subsurface or
  cycle would force that smoothing to have positive genus, whereas
  every approximating domain is a subsurface of a punctured sphere.
  
  Every nonconstant normalization meets the fixed compact core. If one
  were contained entirely in the negative end, properness would force
  its cylindrical coordinate to attain a maximum, contrary to the
  maximum principle. The inverse image of the core is compact, so it
  meets only finitely many components of the normalization. Thus there
  are finitely many nonconstant components.
  
  A constant normalization is compact by properness, and is therefore
  a sphere. Each connected subtree of constant spheres maps to a single
  point and is finite, again by properness. If it has $m$ vertices and
  $e$ attachments to nonconstant components, stability gives
  $3m\leq2(m-1)+e$, hence $e\geq m+2$. Since the incidence graph is a
  forest and there are finitely many nonconstant vertices, there are
  only finitely many constant components and nodes as well.
  
  The end-energy bounds pass to every normalization of $u_\infty$:
  first pass to the limit on a fixed truncated end, and then exhaust
  the negative end, using positivity.
  Lemma~\ref{lem:rc-exterior-finite-type} therefore makes every
  noncompact normalization a sphere with finitely many negative Reeb
  punctures. Compact normalizations are spheres. We have obtained a
  finite nodal curve with finitely many punctures.
  
  Choose disjoint asymptotic annuli at all punctures and a sufficiently
  small regular collar level $r$ below every nodal image and the
  compact remainder of the curve, cutting each asymptotic annulus
  once. The exterior images above this level, together with a thin band
  around the cut, lie in a common compact subset of $W$. We may
  therefore choose a compact target domain whose inner boundary is
  $\partial D_r$, which contains all the exterior images above the
  cut, and whose remaining boundary misses them.
  Lemma~\ref{lem:rc-transverse-cuts} identifies the entire exterior
  domain above the cut, for all sufficiently large $k$, with a
  smoothing of the truncated nodal limit. The level is then regular
  for $u_k$, all cutting circles correspond, and the node-smoothing
  regions lie above the cut. Smooth convergence near the cutting
  circles identifies their homotopy classes in $D_r$ with those of
  the nonconstant Reeb asymptotes. A connected component of the nodal
  limit without a puncture would smooth to a closed component of the
  original connected punctured domain, which is impossible.
  
  Fill the original punctures of the domain of $u_k$. Cutting the
  resulting sphere along the level circles gives a finite tree of
  interior and exterior regions. By hypothesis, one interior region
  contains all the original punctures. Root the tree there. Each
  nonroot leaf is then an unpunctured disk. An interior leaf would
  bound in $D_r$ a loop homotopic to a nonconstant Reeb orbit. Every
  such orbit on the cotangent-torus boundary represents a nonzero
  class in $\pi_1(D_r)=\pi_1(T^2)$. An interior leaf is therefore
  impossible, and there is an exterior leaf disk.

  The complete cut configuration is connected, with only the original
  Floer ends left open; every other end is matched to another component
  or capped. Figure~\ref{fig:exterior-building} encodes that topology for
  a cylinder: starting at the main interior component, either side
  branch passes through an exterior cylinder and an interior cylinder
  before terminating in an exterior disk. All seven circles along the
  dashed line are matched, and gluing recovers a cylinder with exactly
  the original input and output ends. The two rows indicate whether a
  piece maps inside or outside $D_r$. For an operation with $k$
  original ends one changes the ends attached to the main component;
  the caption explains why a terminal disk must occur on the exterior
  side.

\begin{figure}[htbp]
\centering
\resizebox{\linewidth}{!}{%
\begin{tikzpicture}[x=1cm,y=1cm,
  font=\small,
  ink/.style={draw=black!80,line width=.7pt},
  guide/.style={draw=black!35,densely dashed,line width=.45pt},
  leader/.style={draw=black!55,line width=.45pt},
  tube/.style={draw=black!80,line width=.7pt,
               double distance=4.8mm,line cap=butt,line join=round}]
\definecolor{inside}{RGB}{230,239,230}
\definecolor{outside}{RGB}{226,235,245}

\draw[guide] (-6.2,0) -- (5.8,0);
\node[anchor=east,align=right] at (-6.35,1.15)
  {Exterior\\$M\setminus D_r$};
\node[anchor=east,align=right] at (-6.35,-1.05)
  {Interior\\$D_r$};

\path[fill=inside,ink]
  (-1.24,0) -- (-1.24,-.48)
  .. controls (-1.52,-.72) and (-1.55,-1.10) .. (-1.20,-1.40)
  .. controls (-1.04,-1.55) and (-.96,-1.67) .. (-.96,-1.91)
  -- (-.96,-2.45) -- (-.48,-2.45) -- (-.48,-1.91)
  .. controls (-.48,-1.63) and (-.20,-1.56) .. (0,-1.56)
  .. controls (.20,-1.56) and (.48,-1.63) .. (.48,-1.91)
  -- (.48,-2.45) -- (.96,-2.45) -- (.96,-1.91)
  .. controls (.96,-1.67) and (1.04,-1.55) .. (1.20,-1.40)
  .. controls (1.55,-1.10) and (1.52,-.72) .. (1.24,-.48)
  -- (1.24,0) -- (.76,0) -- (.76,-.43)
  .. controls (.76,-.61) and (.24,-.61) .. (.24,-.43)
  -- (.24,0) -- (-.24,0) -- (-.24,-.43)
  .. controls (-.24,-.61) and (-.76,-.61) .. (-.76,-.43)
  -- (-.76,0) -- cycle;

\draw[tube,double=inside]
  (-5,0) -- (-5,-.30)
  .. controls (-5,-1.43) and (-3.4,-1.43) .. (-3.4,-.30)
  -- (-3.4,0);
\draw[tube,double=inside]
  (3.4,0) -- (3.4,-.30)
  .. controls (3.4,-1.43) and (5,-1.43) .. (5,-.30)
  -- (5,0);

\draw[tube,double=outside]
  (-3.4,0) -- (-3.4,.35)
  .. controls (-3.4,2.50) and (-1,2.50) .. (-1,.35)
  -- (-1,0);
\draw[tube,double=outside]
  (1,0) -- (1,.35)
  .. controls (1,2.50) and (3.4,2.50) .. (3.4,.35)
  -- (3.4,0);

\foreach \x in {-5,0,5}{
  \begin{scope}[xshift=\x cm]
    \path[fill=outside,ink]
      (-.24,0) -- (-.24,.35)
      .. controls (-.24,.64) and (-.62,.68) .. (-.62,1.00)
      .. controls (-.62,1.61) and (.62,1.61) .. (.62,1.00)
      .. controls (.62,.68) and (.24,.64) .. (.24,.35)
      -- (.24,0) -- cycle;
  \end{scope}
}

\foreach \x in {-5,-3.4,-1,0,1,3.4,5}{
  \fill[white] (\x,0) ellipse [x radius=.24,y radius=.075];
  \draw[ink,dash pattern=on 1.4pt off 1.1pt]
    (\x-.24,0) arc[start angle=180,end angle=0,
                    x radius=.24,y radius=.075];
  \draw[ink]
    (\x-.24,0) arc[start angle=180,end angle=360,
                    x radius=.24,y radius=.075];
}

\foreach \x in {-.72,.72}{
  \filldraw[fill=white,ink]
    (\x,-2.45) ellipse [x radius=.24,y radius=.075];
  \draw[guide] (\x-.24,-2.47) -- (\x-.24,-2.80);
  \draw[guide] (\x+.24,-2.47) -- (\x+.24,-2.80);
}
\draw[-{Stealth[length=1.5mm]},line width=.55pt]
  (-.72,-2.84) -- (-.72,-2.52);
\draw[-{Stealth[length=1.5mm]},line width=.55pt]
  (.72,-2.52) -- (.72,-2.84);
\node[below] at (-.72,-2.84) {$\gamma_-$};
\node[below] at (.72,-2.84) {$\gamma_+$};

\node[align=center] at (0,-1.04)
  {main interior\\component};
\node[align=center] at (-4.2,-1.65) {interior cylinder};
\node[align=center] at (4.2,-1.65) {interior cylinder};
\node at (-2.2,2.55) {exterior cylinder};
\node at (2.2,2.55) {exterior cylinder};
\node at (-5,1.03) {disk};
\node at (0,1.03) {disk};
\node at (5,1.03) {disk};

\node[anchor=west,align=left] at (2.05,-2.48)
  {All original ends lie on\\the main interior component.};
\draw[leader] (2.00,-2.39) -- (1.36,-2.39) -- (1.10,-1.65);
\node[anchor=east,align=right] at (-2.05,-2.48)
  {Matching circles\\along $u_k^{-1}(\partial D_r)$};
\draw[leader] (-4.90,-2.48) -- (-5.75,-2.48) -- (-5.75,-.27)
  -- (-5.25,-.04);
\end{tikzpicture}%
}
\caption{Topology of a putative escaping cylinder. The main interior
component carries both original Floer ends, and every circle on the
dashed line is matched. The dual graph is a tree. A nonroot leaf is
a disk; it cannot be interior, since its boundary represents a
nontrivial Reeb class in $\pi_1(D_r)$. Thus an exterior disk cap must
occur.}
\label{fig:exterior-building}
\end{figure}
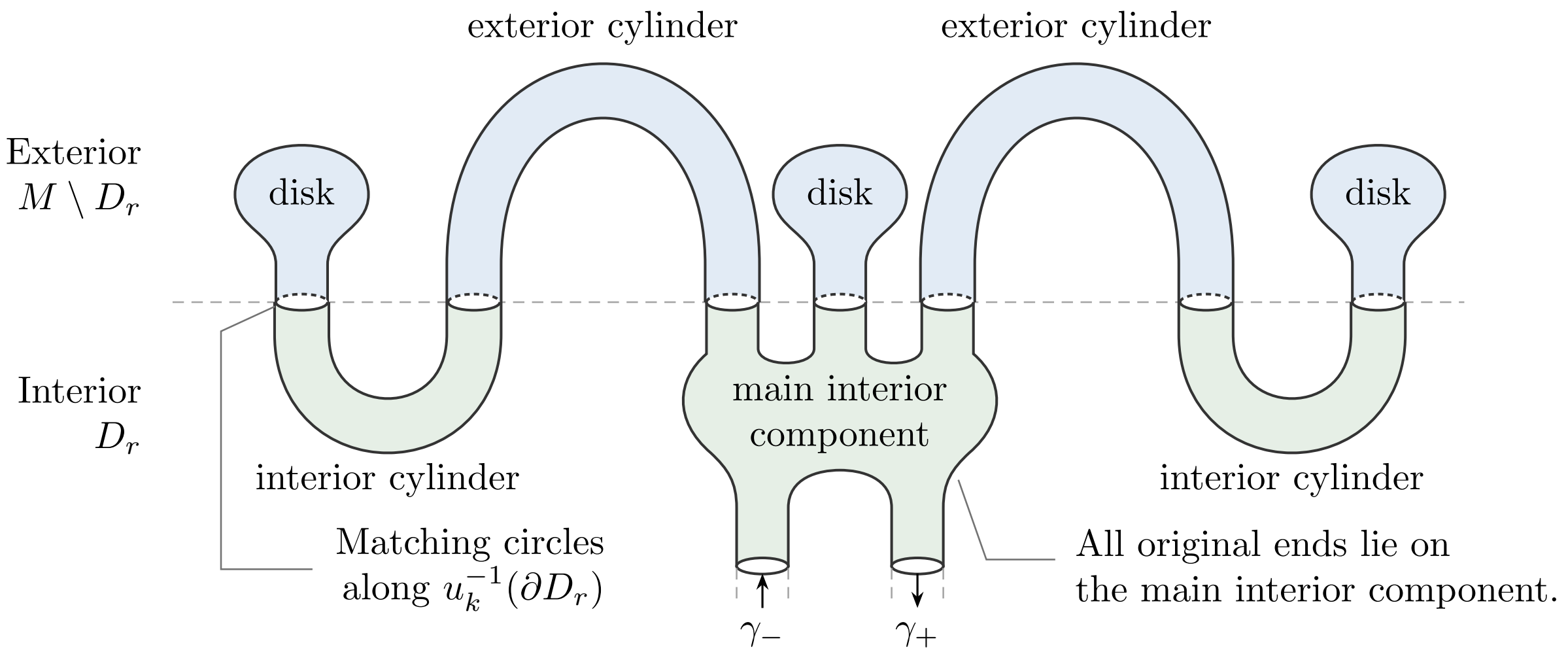

  Under the smoothing identification, this disk comes from a connected
  component of the nodal limit with precisely one negative puncture.
  The normalization carrying that puncture is nonconstant by properness,
  has genus zero, and has no other puncture. Its nodal points are
  ordinary interior points of the normalization. It is a finite-energy
  plane, contradicting the choice of $J$.
  \end{proof}

\begin{proof}[Proof of Proposition~\ref{prop:rc-confinement}]
  Put $E=b-a$. In the noncompact case choose the flat plateau as in
  Lemma~\ref{lem:finite-type-flat-exterior}, so that every cylinder
  of energy at most $E$ lies in a fixed compact set before the
  exhaustion tail.
  
  Suppose the proposition fails. Then there are regular strictly
  $S$-shaped Floer data $(H_k,J_k)$ with outer bends contained in
  $D_{\delta_k}$, where $\delta_k\downarrow0$, and plateau heights
  $C_k:=C(H_k)>b+k$, for which there are contributing Floer
  cylinders $u_k$ meeting $M\setminus V(D)$. The structures $J_k$
  agree with the fixed exterior structure $J$ outside
  $D_{\delta_k}$.
  
  By \eqref{eq:rc-window-energy},
  \[
   E_{\mathrm{geo}}(u_k)=E_{\mathrm{top}}(u_k)<E.
  \]
  Since $H_k$ is spatially constant on the exterior plateau,
  the restrictions of $u_k$ outside $D_{\delta_k}$ are
  $J$-holomorphic. In the noncompact case the choice of plateau
  ensures that the curves do not reach the exhaustion tail.
  Their exterior images are therefore contained in the union
  of a fixed compact subset of $W$ and its negative cylindrical
  end, and their exterior areas are bounded by $E$.
  
  Fix a sufficiently small collar level $r>0$. For sufficiently
  large $k$ we have $\delta_k<r$, so near this level $H_k=C_k$
  and $J_k=J$. Whenever the level is regular for $u_k$,
  Lemma~\ref{lem:rc-cut-component}, applied with $\sigma_r=0$
  and using $C_k>b$, puts both asymptotic ends in the same
  connected component of $u_k^{-1}(D_r)$.
  
  All the hypotheses of
  Lemma~\ref{lem:genus-zero-exterior-confinement} are thus
  satisfied. The lemma implies that the images of $u_k$ lie
  in $V(D)$ for sufficiently large $k$, a contradiction.
  \end{proof}
\begin{proof}[Proof of Theorem~\ref{thm:dimension-four-collapse}]
Suppose $D_{\mathrm{def}}\ne0$ and choose $a,b,\lambda$ as in
Proposition~\ref{prop:rc-strict-detection}. Choose a plane-free and sphere-free
exterior almost complex structure and the confinement bounds of
Proposition~\ref{prop:rc-confinement} for the larger window
$[a,b+\lambda)$. By Proposition~\ref{prop:rc-strict-detection} we
can find a strictly $S$-shaped Hamiltonian satisfying these bounds
for which multiplication by $T^\lambda$ into the shifted window is
not injective on cohomology.

But all cylinders contributing in the larger window stay in $V(D)$,
where $\omega=d\theta$, so by Stokes $E_\theta(u)=0$. The differential
therefore preserves the actual relative exponent. In particular,
multiplication by $T^\lambda$ from $[a,b)$ to
$[a+\lambda,b+\lambda)$ identifies the source, as a complex over $\kk$, with the direct
summand of the target consisting of relative exponents at least $\lambda$.
Indeed, the projection onto these exponents commutes with the
differential, and division by $T^\lambda$ shifts the ordinary window
back to $[a,b)$. The map is thus injective on cohomology, a
contradiction. Hence $D_{\mathrm{def}}=0$, and the relative-filtered
retract \eqref{eq:rc-retract} proves collapse at $E_1$.
\end{proof}

For reconstruction we also need to compare restriction, product
and BV with their intrinsic counterparts. We turn to these operations
in Section~\ref{sec:regular}.

\section{Reconstruction over the regular locus}
\label{sec:regular}

Let $M$ have real dimension four and let its regular fibers be
Maslov-zero Lagrangian two-tori. Over a convex polytope in a regular
affine chart, we will identify supported Floer cohomology with
analytic polyvector fields on a domain in the dual torus. To do so
we extend the confinement argument of Section~\ref{sec:collapse} to
continuation maps, the product and the BV operator, first in finite windows and then on
completed cohomology. Compatibility with restriction will allow us
to glue these torus domains over $B_{reg}$. We will show that
the logarithmic volume forms of the oriented charts glue as well.

\subsection{Regular torus domains and tropical markings}
\label{subsec:regular-leading-terms}

A \emph{regular convex polytope} is a compact convex polytope with
rational facet directions in an integral-affine chart of
$B_{\rm reg}$ containing a neighborhood of the polytope. Supporting
constants may be arbitrary real numbers, and faces and points are
allowed. Write
\[
 K_Q=\pi^{-1}(Q),\qquad \cF^*(Q)=SH_M^*(K_Q;\Lambda).
\]
In an action--angle chart with torus fiber $L$, put
$M_{\mathbb Z}=H_1(L;\mathbb Z)$ and
$N_{\mathbb Z}=M_{\mathbb Z}^{\vee}$. The dual analytic torus has
characters $z^m$, $m\in M_{\mathbb Z}$, and valuation map
\[
 \langle m,\Log(x)\rangle=-\log|z^m(x)|.
\]

Let $A(Q)$ be the algebra of convergent Laurent series on
$\Log^{-1}(Q)$, with norm
\[
 \left\|\sum_m a_mz^m\right\|_Q
 =\sup_m |a_m|e^{-s_Q(m)},
 \qquad s_Q(m)=\min_{b\in Q}\langle m,b\rangle.
\]
The normalized character $e_{Q,m}=T^{-s_Q(m)}z^m$ has norm one.
Write $\vartheta_{Q,m}$ for its image in
$\widetilde A(Q)=A(Q)^\circ/A(Q)^{\circ\circ}$, where the two
lattices consist of elements of norm at most one and strictly less
than one, respectively. If $F_Q(m)$ is the minimizing face of
$\langle m,\cdot\rangle$, then
\begin{equation}\label{eq:regular-symbol-product}
 \begin{gathered}
 \vartheta_{Q,m}\vartheta_{Q,m'}=
 \begin{cases}
 \vartheta_{Q,m+m'},&F_Q(m)\cap F_Q(m')\ne\varnothing,\\
 0,&\text{otherwise},
 \end{cases}\\[2pt]
 r_{PQ}(e_{P,m})=T^{s_Q(m)-s_P(m)}e_{Q,m}
 \qquad(Q\subset P).
 \end{gathered}
\end{equation}
The restriction exponent is nonnegative; the reduced symbol
survives precisely when $s_Q(m)=s_P(m)$.

The analytic polyvectors are
\[
 \cA^*(Q)=\bigwedge_{A(Q)}^*\Der_{\Lambda,\mathrm{cont}}A(Q),
\]
with the coefficient norm in the logarithmic frame
$E_i=z_i\partial_{z_i}$. Put $N_\kk=N_{\mathbb Z}\otimes\kk$ and
identify $H^*(L;\kk)$ with $\bigwedge^*N_\kk$, so that the
cohomology classes dual to the chosen homology basis correspond to
the $E_i$. The reduction of the unit ball of polyvectors is
\[
 \cA^*(Q)^\circ/\cA^*(Q)^{\circ\circ}
 =\widetilde A(Q)\otimes_\kk\bigwedge^*N_\kk
 =\bigoplus_{m\in M_{\mathbb Z}}
       \kk\vartheta_{Q,m}\otimes_\kk\bigwedge^*N_\kk.
\]
These normalized symbols, with their lattice labels, coefficients
and restriction weights, form the \emph{analytic tropical marking}.

The Floer marking can be read on sufficiently small regular convex
sets. On such a set $Q$ the cotangent-torus calculation in
\cite[Theorem~5 and Section~6]{Paper1} gives the unweighted symbols
\[
 \kk\varepsilon_{Q,m}\otimes_\kk H^*(L;\kk),
 \qquad m\in M_{\mathbb Z}.
\]
Here $m$ is the winding class and $m=0$ denotes the critical
block. For a smooth contact-type approximation, the short-window
identification in Theorem~\ref{thm:background-indicator-input} is
made after fixing a finite Reeb-period cutoff. On a sufficiently
small regular convex set the uniform version
\cite[Theorem~5]{Paper1} includes all labels. In either calculation
the contribution of a label present in the window is
\[
 \bigl(\kk\varepsilon_{Q,m}\otimes_\kk H^*(L;\kk)\bigr)
       \otimes_\kk\Lambda_I^{\log}.
\]
Here $I$ is indexed by log norm, as in
Section~\ref{subsec:background-filtrations}.

We give $\varepsilon_{Q,m}$ the action normalization corresponding
to $e_{Q,m}$. For $Q\subset P$, restriction transports the class on the
fiber and acts on the truncated Novikov coefficient by $T^{s_Q(m)-s_P(m)}$. The unweighted local blocks,
with these labels, coefficient identifications and weights, form
the \emph{Floer tropical marking}. Its distinguished degree-zero
element in the block labelled by $m$ is
$\varepsilon_{Q,m}\otimes1$.

A Floer--analytic comparison \emph{intertwines the tropical
markings} if its normalized local leading-term identifications send
\[
 \varepsilon_{Q,m}\otimes\alpha
 \longmapsto\vartheta_{Q,m}\otimes\alpha,
 \qquad \alpha\in H^*(L;\kk)=\bigwedge^*N_\kk,
\]
compatibly with restriction weights and integral-affine transport.
For a larger regular polytope, the assertion about leading terms
can be tested at its vertices. Indeed, a normalized character survives
at precisely the vertices of its minimizing face. The Laurent norm
on the polytope is the maximum of the vertex norms, so reduction
at the vertices detects its entire leading term.
The source terms here are unweighted local classes; their
ordinary-window realizations include the truncated coefficient
factor displayed above.

The standard logarithmic volume form is
\[
 \Omega_0=d\log z_1\wedge d\log z_2.
\]
For a nowhere-vanishing volume form $\Omega$, our divergence
convention is
\begin{equation}\label{eq:regular-divergence}
 \diver_\Omega|_{\PV^q}=(-1)^q I_{q-1}^{-1}dI_q
 \quad(q\geq1),\qquad I_q(\xi)=\iota_\xi\Omega,
\end{equation}
and $\diver_\Omega$ is zero on functions. We use
$\iota_{X_1\wedge\cdots\wedge X_q}
 =\iota_{X_q}\cdots\iota_{X_1}$.
Set $\Delta_0=\diver_{\Omega_0}$. We fix the overall sign of Floer
BV by orienting its circle parameter so that, on an unweighted
cotangent-torus block,
\[
 \Delta_F(\varepsilon_{Q,m}\otimes\alpha)
 =-\varepsilon_{Q,m}\otimes\iota_m\alpha.
\]
Here $\iota_m$ is contraction under the ordinary duality between
$H^1(L;\kk)$ and $H_1(L;\mathbb Z)$.
The product has degree zero and BV has degree $-1$. With the bracket
convention
\begin{equation}\label{eq:bracket}
 \{x,y\}=(-1)^{|x|}\bigl(
 \Delta(x*y)-\Delta(x)*y-(-1)^{|x|}x*\Delta(y)\bigr),
\end{equation}
the analytic bracket is the Schouten bracket. In particular,
\[
 \Delta_0(fE_i)=-E_i(f),\qquad \{E_i,f\}=E_i(f).
\]

\subsection{The regular reconstruction theorems}

\begin{theorem}[Local marked regular reconstruction]
\label{thm:local-regular-reconstruction}\label{thm:fixed-j-locality}
Fix a regular convex polytope $P_0$ and an action--angle chart over
a neighborhood of $P_0$. There is a family of isometric
isomorphisms of graded $\Lambda$-BV algebras
\[
 \Phi_Q:
 \bigl(\cF^*(Q),\Delta_F\bigr)
 \xrightarrow{\ \sim\ }
 \bigl(\cA^*(Q),\diver_{\Omega_0}\bigr),
 \qquad Q\subset P_0,
\]
where $Q$ ranges over regular convex polytopes. These isomorphisms
intertwine the Floer and analytic tropical markings and are natural
for every inclusion $Q_1\subset Q_2\subset P_0$.
\end{theorem}

The local comparisons may depend on the choices in their construction.
The maps used for gluing come from Floer restriction. For $Q\subset P$,
the map of affinoids $Y_Q\to Y_P$ identifies, in local torus
coordinates, with the inclusion of the domain over $Q$ in the one
over $P$.

\begin{theorem}[Regular reconstruction and volume form]
\label{thm:global-reconstruction}\label{thm:bv-reconstruction}
The affinoids $Y_Q=\operatorname{Sp}\cF^0(Q)$ glue by Floer
restriction to a smooth rigid analytic surface $Y_{\rm reg}$ with
a continuous proper surjection
\[
 p:Y_{\rm reg}^{\rm Berk}\longrightarrow B_{\rm reg},
 \qquad p^{-1}(Q)=Y_Q^{\rm Berk}
\]
for every regular convex polytope $Q$. The torus identifications
supplied by Theorem~\ref{thm:local-regular-reconstruction}
identify $p|_{Y_Q^{\rm Berk}}$ with $\Log$. Finite regular
polyhedral coverings pull back to admissible affinoid coverings.

The chosen orientation of $B_{\rm reg}$ determines a
nowhere-vanishing analytic two-form $\Omega_{\rm reg}$, given by
$d\log z_1\wedge d\log z_2$ in every oriented torus chart supplied
by that theorem. There are natural isomorphisms
\[
 \cF^q(Q)\xrightarrow{\ \sim\ }
 \Gamma\bigl(Y_Q,\textstyle\bigwedge^q T_{Y_{\rm reg}}\bigr)
\]
in all degrees. These comparisons commute with restriction,
intertwine the tropical markings, and together are isomorphisms of
BV algebras with analytic BV operator
$\diver_{\Omega_{\rm reg}}$. In degree zero they identify
$\cF^0(Q)$ with $\Gamma(Y_Q,\cO_{Y_{\rm reg}})$.

The class $s_Q\in\cF^2(Q)$ corresponding to the bivector dual to
$\Omega_{\rm reg}|_{Y_Q}$ generates $\cF^2(Q)$ freely over
$\cF^0(Q)$. These classes are BV-closed and compatible with
restriction.

The space, projection and volume form are independent of the
auxiliary choices and the regular cover. Reversing the chosen
orientation changes the sign of $\Omega_{\rm reg}$ and leaves
its divergence unchanged.
\end{theorem}

\subsection{Locality from neck stretching}

Every regular convex polytope $P_0$ is contained in the interior of a
smooth strictly convex domain $P$ in the same chart, whose boundary has
positive curvature. The torus fibration over a neighborhood of $P$
admits action--angle coordinates. Indeed, the base is contractible, so
the torus bundle is trivial; the closed base two-form in a
trivialization is exact and is removed by translating the section. We
may thus write $D=\pi^{-1}(P)=L\times P$, with its action--angle
primitive, and choose a slightly larger exact neighborhood $V(D)$.
Centering the primitive at an interior point of $P$ makes $D$ a
cotangent-torus Liouville domain.

The exterior-plane argument of Lemma~\ref{cor:no-exterior-planes}
applies to this domain: it uses positive curvature, the Maslov-zero
grading and dimension four, without a smallness condition on $P$.
Choose a cylindrical, plane-free and sphere-free exterior almost complex structure $J$
as in Section~\ref{sec:collapse}; in the noncompact case it agrees with
the fixed geometrically bounded structure outside a compact set. Use
this same $D$ and $J$ for every $Q\subset P_0$.

For the comparison with $\widehat D$, we will count the curves for
continuation maps, the product and BV operator inside $V(D)$, where the
ambient and intrinsic Floer data agree. We use continuation to identify
the two cohomologies in each ordinary window on the relative lattice.
After passing to completed cohomology, we will return to a single
window to test the product and the BV operator.

\paragraph{Floer data for the comparison.}
We use the strictly $S$-shaped Hamiltonians and finite-window complexes
of Section~\ref{sec:collapse}. For operations, choose split Hamiltonian
data on a genus-zero punctured surface $\Sigma$:
\begin{equation}\label{eq:regular-split-data}
 \mathcal K=F\beta,\qquad d\beta=0,\qquad
 d_\Sigma F_x\wedge\beta\geq0\quad(x\in M).
\end{equation}
Here $F:\Sigma\times M\to\mathbb R$, and near a puncture
$\beta=w_i\,dt$ and $F\beta=H_i\,dt$. With output weight one, take
the $m$ input weights to sum to one and satisfy $w_i\geq2^{1-m}$
as in \cite[Equation~(2.4)]{AbouzaidGromanVarolgunes24}. For the
pair of pants take both input weights equal to $1/2$.
The energy identity of
\cite[Equations~(3.16)--(3.17)]{AbouzaidGromanVarolgunes24} gives,
on every subdomain $S\subset\Sigma$,
\begin{equation}\label{eq:regular-operation-energy}
 E_{\mathrm{top}}(u|_S)
 =E_{\mathrm{geo}}(u|_S)
   +\int_S d_\Sigma F_{u(z)}\wedge\beta
 \geq E_{\mathrm{geo}}(u|_S)\geq0.
\end{equation}

For the differential use the stationary Floer datum. For operations
requiring a later output, choose its Hamiltonian to satisfy
the domination condition of
\cite[Definition~2.18]{AbouzaidGromanVarolgunes24} and
\[
 \min_t H_0(x)>\max_i\frac{\max_t H_i(x)}{w_i}
 \quad\text{for all }x
\]
for the weights in the family. The inputs are strictly negative on
the support, where the output can be taken closer to zero. On the
remaining compact region this choice follows from cofinality. For
the BV family we impose the same inequality for all time rotations
of the input.

Lemma~2.19 of \cite{AbouzaidGromanVarolgunes24} constructs the split
datum \eqref{eq:regular-split-data}. We need it to be spatially
constant on a common exterior plateau. This holds in the proof of
that lemma upon choosing its intermediate Hamiltonian constant there.
For compact families use the convex fibers in the proof of Lemma~2.21
of the same paper and the collar gluing of its Section~2.2, as in
Proposition~2.20. For each fixed $\beta$, spatial constancy on the
plateau is preserved by convex interpolation, and the end data agree
there under gluing. Thus the construction works relative to the
prescribed faces with this exterior condition. When all outer bends
lie inside $D_\delta$, the smooth Floer structures can be chosen
equal to $J$ on $M\setminus D_\delta$.
One can extend the compatible metrics from the ends and parameter
faces, keeping the prescribed metric on the exterior collar, and
take the associated $\omega$-compatible almost complex structures.

On the common plateau the Hamiltonian one-form equals the scalar
one-form $\kappa=F(\,\cdot\,,x_0)\beta$ on $\Sigma$, for any fixed
point $x_0$ of the plateau. In particular $d\kappa\geq0$, and
$\kappa=C_i\,dt$ at the $i$-th end of plateau height $C_i$.
In the finite-type case we put the exhaustion tails beyond the
common plateau. For inputs with tail coefficient $\varepsilon$, an
output coefficient $2^{m-1}\varepsilon$ and a sufficiently high
plateau give the required domination at infinity. All coefficients
in a finite diagram can be chosen admissibly small.
Lemma~\ref{lem:rc-regular-companion} compares the resulting supported
windows using a sufficiently wide common plateau.

For regularity we reserve a small spatially constant
increase of $F$ along an incoming cylindrical end,
using the strict domination by the output Hamiltonian.
This does not change the equation, but permits small
domain-dependent Hamiltonian perturbations on the inner
side of the cutting annulus while preserving
$d_\Sigma F_x\wedge\beta\geq0$.

For a fixed energy bound, compactness and the exclusion
of breaking at exterior orbits allow us to choose the
perturbation collar long enough that every relevant
solution meets this interior region there. Variations
of the spatial differential of $F$ then give universal
transversality, including for horizontal solutions.
The choices are made relative to the regular parameter
faces as in \cite[Section~3.4]{AbouzaidGromanVarolgunes24}.
A sphere not excluded by exactness or by the fixed
sphere-free exterior meets a region where the almost
complex structure can vary. The simple-sphere attachment
locus has codimension two, so the zero- and one-dimensional
moduli spaces have no sphere bubbles.

\begin{proposition}[Confinement for operations in a finite window]
\label{prop:regular-operation-confinement}
Fix a plane-free and sphere-free cylindrical almost complex structure
$J$ on $W=M\setminus D$. Consider a finite collection of genus-zero
operations with at least one input and a single output, and their compact
comparison families. For each operation fix an ordinary lower bound
$a_i$ at each input and an ordinary upper bound $b$ at the output.

There are $\delta_0>0$ and $C_0>0$ with the following property.
For $0<\delta<\delta_0$, take regular strict data as above with all
outer bends inside $D_\delta$ and all input plateau heights greater
than $C_0$. Use a smooth family of compatible almost complex
structures $J_{\delta,z}$ on $M$, indexed by the domain point $z$
and, for comparison families, by the parameter. For every such
choice of indices require
\[
 J_{\delta,z}|_{M\setminus D_\delta}=J|_{M\setminus D_\delta}.
\]
In the finite-type case take the common flat plateau sufficiently
wide for the prescribed bounds.

Then every contributing solution with local asymptotic orbits in
$D_\delta$, inputs of nonnegative relative valuation, and the
prescribed ordinary bounds has image in $V(D)$. This holds for the
differential, monotone continuation, the product, the BV operator,
and the comparison families for these maps.
\end{proposition}

\begin{proof}
Write an input as $T^{q_i}\widetilde a_{\gamma_i}$, with $q_i\geq0$,
and put
\[
 A_i=A_{H_i,\theta}(\gamma_i),\qquad
 \nu_i=q_i-A_i.
\]
For a contribution whose output has ordinary depth $\nu_0<b$, the
energy identity is
\begin{equation}\label{eq:regular-total-window-energy}
 E_{\mathrm{top}}(u)=\nu_0-\sum_i\nu_i
       < b-\sum_i a_i.
\end{equation}
Equation~\eqref{eq:regular-operation-energy} therefore bounds the
geometric energy uniformly. If the right-hand side is nonpositive
there are no such contributions, so assume it is positive.

We must show that all the asymptotic ends belong to one component
of $u^{-1}(D_r)$. Fix a small regular collar level $r>\delta$ on the
flat plateau and consider the component containing the output. If an
input lies elsewhere, some other component $S$ contains a nonempty
set $I$ of inputs and no output. Orient its cut circles as boundaries
of $S$. Near the cut the equation is $J$-holomorphic, so the
contact-type inequality gives $\int_{\partial_{\rm cut}S}u^*\theta
\geq0$, as in Lemma~\ref{lem:rc-cut-component} with zero slope.
Stokes' theorem and $d\kappa\geq0$ give
\[
 \int_{\partial_{\mathrm{cut}}S}\kappa
       =\sum_{i\in I}C_i+\int_Sd\kappa
       \geq\sum_{i\in I}C_i.
\]
Using the action convention of Section~\ref{sec:collapse}, we obtain
\[
 \begin{split}
 E_{\mathrm{top}}(u|_S)
 &=\sum_{i\in I}A_i+
       \int_{\partial_{\mathrm{cut}}S}(u^*\theta+\kappa)\\
 &\geq\sum_{i\in I}(A_i+C_i).
 \end{split}
\]
The complementary subdomain has nonnegative topological energy by
\eqref{eq:regular-operation-energy}. Comparing the last inequality
with the total energy and substituting $\nu_i=q_i-A_i$ for $i\in I$
gives
\begin{equation}\label{eq:regular-input-subset-estimate}
 \begin{split}
 \sum_{i\in I}C_i
 &\leq\nu_0-\sum_{i\in I}q_i
                    -\sum_{j\notin I}\nu_j\\
 &<b-\sum_{j\notin I}a_j.
 \end{split}
\end{equation}
Choose the input plateaux so that
\begin{equation}\label{eq:regular-input-plateaux}
 \sum_{i\in I}C_i>b-\sum_{j\notin I}a_j
 \qquad(\varnothing\ne I\subset\{1,\ldots,m\}).
\end{equation}
There are only finitely many subsets and operations, so one $C_0$
suffices. The two inequalities rule out a non-empty $S$. Thus every input lies
in the component containing the output, independently of the
interior Hamiltonian slopes.

If no $\delta_0$ worked, there would be escaping solutions $u_k$
for data with $\delta_k\downarrow0$. Outside $D_{\delta_k}$ the
smooth Floer structures agree with the same $J$, and the spatial
Hamiltonian vector field vanishes on the plateau. In the finite-type
case the monotonicity argument of Lemma~\ref{lem:finite-type-flat-exterior}
on the common plateau keeps curves with this energy bound in a fixed
compact region before the tail. Thus the exterior portions are
$J$-holomorphic and their $\omega$-areas are bounded by
\eqref{eq:regular-total-window-energy}.

For every fixed sufficiently small $r>0$ we eventually have
$\delta_k<r$. Whenever the cut is transverse, the calculation above
puts all ends of $u_k$ in one component of $u_k^{-1}(D_r)$.
Lemma~\ref{lem:genus-zero-exterior-confinement} now puts the whole
image in $V(D)$ for large $k$, a contradiction. The same reasoning
applies with the parameter varying along the sequence, including
when the domains degenerate. This gives the asserted uniformity.
\end{proof}

\paragraph{Computing operations with strict data.}
We record how to compute a product or a BV image for specified
inputs at a finite Hamiltonian stage. For the given ordinary bounds,
Proposition~\ref{prop:regular-operation-confinement} confines the
strict solutions with local ends to the exact neighborhood $V(D)$.
There Stokes' theorem gives
\[
 E_\theta(u)=0,\qquad q_{\mathrm{out}}=\sum_iq_i.
\]
Thus these solutions preserve the relative exponent in the count.
To use the count we must also exclude breaking at an omitted
exterior orbit and compare with the operation for regular data.
For both purposes we will cut on the protected annulus, on the
inner side of the omitted orbits.

To identify this count with the usual operation, use the
continuation isomorphisms from each strict end to its regular
companion in Lemma~\ref{lem:rc-regular-companion}. We compare
continuing the inputs before the operation with continuing the
output afterwards. For this homotopy we also need monotonicity
after subtracting the value of the Hamiltonian on the cutting
annulus, as in \eqref{eq:rc-normalized-monotonicity}.

For these unary and binary comparisons, choose the end
Hamiltonians with a common normalized exterior profile $f$ after dividing by their weights and subtracting their
cut values. Thus if the normalized slope on the protected annulus
is $\sigma$, the actual end Hamiltonian has slope $w_i\sigma$.
Choose $\sigma$ small enough for every end package. For the product
these are the input and doubled output packages of
\cite[Sections~3.6 and~4.3]{GromanBoundaryDepth}, in our normalization
of the weights. Take a small nonnegative companion perturbation $h$
which vanishes near the annulus. At an end of weight $w_i$ the
perturbation of the end Hamiltonian is $w_i h$. On a mixed surface
choose the normalized exterior profile to be
\[
 F_z-F_z|_{\mathrm{cut}}=f+\rho(z)h,
 \qquad d\rho\wedge\beta\geq0,
\]
where $\rho$ goes from $0$ at strict ends to $1$ at regular ends.
Increasing cutoffs on the continuation collars give these choices
on the two broken configurations. The collar construction and
convex interpolation for fixed $\beta$ extend them over the
homotopy. Since $h\geq0$, the normalized curvature on the exterior
is $h\,d\rho\wedge\beta\geq0$.

To see what this gives at a regular level $r$ in the protected
annulus, write $\theta(X_{F_z})=-\sigma_r$ there. The common radial
profile makes $\sigma_r$ independent of $z$. With the cut oriented
as the boundary of the interior, the contact inequality is
$u^*\theta+\sigma_r\beta\geq0$ on its oriented tangent. For a
component carrying inputs $I$ and no output, Stokes gives
$\int_{\rm cut}\beta=\sum_{i\in I}w_i$. Thus the same argument as
in \eqref{eq:regular-input-subset-estimate} excludes such a component
as soon as
\[
 \sum_{i\in I}\bigl(H_i|_r-\sigma_r w_i\bigr)
       > b-\sum_{j\notin I}a_j .
\]
We choose the cutting values this high for all relevant subsets,
as in \eqref{eq:rc-high-cut}. An exterior output below the cutoff
is thereby excluded. Breaking at an omitted exterior orbit would
give a neck with a whole circle beyond this cut. In the unbroken
surface that circle separates a nonempty set of original inputs
from the output, since each Floer vertex has an input and the
domain has genus zero. The same inequality excludes this breaking
for the strict data ($h=0$) and throughout the comparison homotopy.

For local asymptotics the exterior has no original ends, so the
total integral of $\beta$ over its cut boundary is zero. Applying
Stokes on the interior and the energy identity on the exterior
therefore gives
\[
 E_\theta(u)\ \geq\
 E_{\mathrm{geo}}(u|_{\mathrm{ext}})
 +\int_{\mathrm{ext}}d_\Sigma(F-F|_r)_{u(z)}\wedge\beta
 \ \geq\ 0.
\]
Here $\mathrm{ext}=u^{-1}(M\setminus D_r)$, and the final inequality
uses the normalized curvature just computed. The comparison homotopy
therefore preserves the relative lattice. In the noncompact case
take $h$ compactly supported and keep $J_0$ at infinity. On the tail
the base Hamiltonian has a fixed normalized exhaustion profile, so
its spatial Hamiltonian vector field is bounded and its curvature is
the curvature of the scalar cut value. On each domain chart the
graph metric is therefore uniformly comparable to the product metric
by \cite[Lemma~A.1]{GromanBoundaryDepth}.
These bounds persist on gluing collars, where the data are cylindrical.
They give the common taming data at infinity for
\cite[Definition~C.12 and Proposition~C.18]{AbouzaidGromanVarolgunes24},
and hence the required compactness for the mixed homotopy.
Its boundary identifies the
strict product and BV counts on these window classes with the
operations computed using regular companions.

We now use the differentials and continuation maps to compare the
directed systems defining supported cohomology. We will need this
comparison to commute with changes of window and restrictions of
the support.

\begin{proposition}[Comparison in an ordinary window]\label{prop:regular-window-locality}
For every $Q\subset P_0$ as above and every ordinary window $[a,b)$,
there are isomorphisms
\[
 \phi^+_{J,Q;a,b}:
 SH^+_{M,\theta,[a,b)}(K_Q)
 \xrightarrow{\ \sim\ }
 SH^+_{\widehat D,\theta,[a,b)}(K_Q).
\]
They commute with upper truncations, inclusions of lower bounds,
Novikov scalar shifts, and actual restriction maps. On the local
leading blocks they are the intrinsic identifications of the
normalized orbit orientation lines.
\end{proposition}

\begin{proof}
For the given window, choose strict acceleration data with the
heights, collar width $\delta$ and plateau widths of
Proposition~\ref{prop:regular-operation-confinement}. Match the
ambient and intrinsic data on $V(D)$, including the smooth Floer
structures $J_{\delta,z}$ there. The two strict complexes have the
same local generators, orientation lines and action normalizations.
Their exterior constants are above the upper cutoff on the
nonnegative relative lattice. The proposition confines the ambient
cylinders for the differentials, structural continuations and their
comparison homotopies to $V(D)$.

Intrinsically, first choose a flat plateau wide enough that the
bounded-energy curves cannot reach the exhaustion tail. An exterior
portion beyond a holomorphic contact cut has no original puncture
and satisfies
\[
 \int_{\mathrm{exterior}}u^*d\theta
       =-\int_{\mathrm{cut}}u^*\theta\leq0.
\]
Its holomorphic area vanishes, so it cannot cross the cut.
This applies to continuation and to the strict operations above,
since the spatial Hamiltonian vector field vanishes beyond the cut.
In particular the unary counts identify the differentials and
structural continuation maps.

Lemma~\ref{lem:rc-regular-companion} and
\eqref{eq:rc-completed-window} give
\[
 SH^+_{M,\theta,[a,b)}(K_Q)
 \cong\varinjlim_H HF^{+,\mathrm{str}}_{\theta,[a,b)}(H),
\]
and the analogous intrinsic formula. Taking direct limits gives
$\phi^+_{J,Q;a,b}$. This also applies to faces and points. The
Hamiltonians approximate the indicator of $K_Q$, while the cutting
boundary is the fixed smooth boundary of $D$.

For finitely many windows, use common data with plateaux above all
their upper endpoints. The generator identifications commute with
the endpoint maps and with Novikov multiplication into shifted
windows. For $Q_1\subset Q_2$, ordered acceleration data compute
the actual restriction, and confinement identifies these continuation
maps as well. To compare two choices of cutting annulus for a
fixed window, choose both outer bends strictly inside its confinement
bound. For a map between two windows choose intermediate data with
a collar narrow enough and a plateau high enough for both windows,
and compare to each of the original choices within its own window.
In the noncompact case use a common admissible tail coefficient.
Now use the bridge construction of
\cite[Lemma~3.10(2)]{GromanBoundaryDepth}. Place the additional
protected shell and the final taper beyond the original shells but
inside the confinement bound $D_{\delta_0}$. Between the shells keep the radial derivative positive and below
the least Reeb period. The two cuts on each simultaneous ray
therefore give exactly the same local generators. Compare each
original ray to its simultaneous ray using its original shell.
Beyond the final taper use
a common flat plateau and retain $J$. Within each package the common
lower profiles and interpolations in the proof of
\cite[Proposition~3.9]{GromanBoundaryDepth} can be taken constant
on that plateau. Thus the strict continuation maps and their
homotopies along the bridge satisfy
Proposition~\ref{prop:regular-operation-confinement} in the window.
They give the same comparison on either side, with $J$ fixed.
Finally, the generator identifications
preserve the normalized orientation lines, giving the asserted
identification of the local leading blocks.
\end{proof}

\subsection{Passage to the completed theory}

We have compared the subcomplexes of nonnegative relative valuation
in every finite ordinary window. To recover the full ordinary window we
allow division by Novikov monomials. Density then suffices, since
each element has only finitely many coordinates visible in a finite
window.

\begin{lemma}\label{lem:regular-localization-windows}
Let $C$ be the ordinary completed Floer telescope and let
$C^+\subset C$ be its nonnegative relative lattice, with the
ordinary filtration induced from $C$. For every $a<b$, there is a
natural identification of complexes
\begin{equation}\label{eq:regular-saturate-window}
 G^aC/G^bC
 \cong
 \varinjlim_{\lambda\geq0}
 T^{-\lambda}C^+_{[a+\lambda,b+\lambda)}.
\end{equation}
Consequently,
\[
 H^*(G^aC/G^bC)
 \cong
 \varinjlim_{\lambda\geq0}
 T^{-\lambda}H^*(C^+_{[a+\lambda,b+\lambda)}).
\]
These identifications commute with endpoint maps and coefficient
multiplication.
\end{lemma}

\begin{proof}
Each term on the right of~\eqref{eq:regular-saturate-window}
embeds in $G^aC/G^bC$. For $\mu\geq\lambda$, the transition map
is multiplication by $T^{\mu-\lambda}$ before rescaling by
$T^{-\mu}$.

Let $x\in G^aC$. Finite-coordinate chains are dense in the ordinary
completion, so there is such a chain $y$ with $\nu(x-y)\geq b$.
Then $\nu(y)\geq a$. Since $y$ has only finitely many orbit
coordinates, their relative valuations have a common lower bound.
Thus $T^\lambda y\in C^+$ for some $\lambda\geq0$, proving
surjectivity. The cohomological assertion follows from exactness of
filtered direct limits. All maps used in the construction respect
the indicated compatibilities.
\end{proof}

Applying the lemma to the ambient and intrinsic complexes identifies
their full ordinary windows. To recover completed cohomology we must
also account for the $\varprojlim^1$ term in the Milnor exact sequence.
Here we use that intrinsic cohomology is Hausdorff.

\begin{lemma}\label{lem:regular-window-completion}
Let $C$ and $D$ be complete, separated filtered cochain complexes
over $\Lambda$, with filtration-preserving differentials. Write
$G^aC=\{x:\nu(x)\geq a\}$, and similarly for $D$, and suppose
these filtrations are exhaustive. Assume that there are
isomorphisms
\[
 \phi_{a,b}^k:
 H^k(G^aC/G^bC)\longrightarrow H^k(G^aD/G^bD),
 \qquad a<b,
\]
compatible with both endpoint maps, coefficient multiplication,
and scalar shifts of the windows. If $H^*(D)$ is Hausdorff for its
quotient seminorm, these isomorphisms induce a canonical isometric
$\Lambda$-linear isomorphism
\[
 \Phi:H^*(C)\longrightarrow H^*(D).
\]
The construction is natural for bounded maps whose window maps
commute with the $\phi_{a,b}^k$.
\end{lemma}

\begin{proof}
Fix $a$. The filtration balls are closed, and completeness gives
\[
 G^aC=\varprojlim_{b\to\infty}G^aC/G^bC.
\]
The quotient maps are degreewise surjective. Taking any countable
cofinal sequence of upper endpoints, the Milnor exact sequence
therefore gives
\begin{equation}\label{eq:regular-window-milnor}
 0\longrightarrow
 \varprojlim_b{}^1 H^{k-1}(G^aC/G^bC)
 \longrightarrow H^k(G^aC)
 \longrightarrow
 \varprojlim_b H^k(G^aC/G^bC)
 \longrightarrow0;
\end{equation}
see~\cite[Tag~0CQC]{Stacks}. We shall take this sequence to the direct limit over the
lower endpoints $a$, ordered by decreasing $a$. For full cohomology
we take all $a\in\mathbb R$. For the open lattice used below we take
only $a>-t$. Denote the chosen range by $I$ and write
$C_I=\bigcup_{a\in I}G^aC$, and similarly for $D$. Exactness of
filtered direct limits gives
\begin{equation}\label{eq:regular-window-milnor-full}
\begin{split}
 0\longrightarrow
 \varinjlim_{a\in I}\varprojlim_b{}^1 H^{k-1}(G^aC/G^bC)
 \longrightarrow H^k(C_I)\\
 \longrightarrow
 \varinjlim_{a\in I}\varprojlim_b H^k(G^aC/G^bC)
 \longrightarrow0.
\end{split}
\end{equation}

For $D$, a class coming from the left term is represented at some
$a\in I$ by a cycle $z\in G^aD$ which is a boundary modulo $G^bD$
for arbitrarily large $b$. Thus its image in $H^k(D)$ has quotient
norm zero and vanishes by Hausdorffness. Whenever
$H^k(D_I)\to H^k(D)$ is injective, the class already vanishes in
$H^k(D_I)$. The left term of \eqref{eq:regular-window-milnor-full}
for $D$ is then zero. The window isomorphisms respect both endpoints,
so they identify this term for $C$ with the one for $D$. Comparing
the right terms now gives canonical identifications
\begin{equation}\label{eq:regular-window-recovery}
 H^k(C_I)\cong\varinjlim_{a\in I}\varprojlim_b H^k(G^aC/G^bC)
 \cong H^k(D_I).
\end{equation}
Taking $I=\mathbb R$ gives the asserted map $\Phi$.

The surjection in~\eqref{eq:regular-window-milnor} shows that
$\Phi$ identifies
\[
 \operatorname{im}\bigl(H^k(G^aC)\to H^k(C)\bigr)
 \quad\hbox{and}\quad
 \operatorname{im}\bigl(H^k(G^aD)\to H^k(D)\bigr)
\]
for every $a$. The quotient valuation of a cohomology class is the
supremum of these possible lower bounds, so $\Phi$ is isometric.
Compatibility with coefficient multiplication and scalar shifts
gives $\Lambda$-linearity.

Naturality for bounded maps follows by taking the same limits.
\end{proof}

Section~\ref{sec:detection} will also require the cohomology of the
open subcomplex of chains of norm less than one. On the intrinsic
side, zero boundary depth gives the injectivity needed in
\eqref{eq:regular-window-recovery}.

\begin{lemma}[Ordinary open lattices]\label{lem:ordinary-open-lattice}
Under the hypotheses of Lemma~\ref{lem:regular-window-completion},
suppose in addition that $D$ has boundary depth zero. For every
$t\in\mathbb R$, put
\[
 C_t=\{x\in C:\|x\|<e^t\},\qquad
 D_t=\{x\in D:\|x\|<e^t\}.
\]
The window comparisons induce isomorphisms
$H^*(C_t)\cong H^*(D_t)$ compatible with $\Phi$. In particular,
$H^*(C_t)\to H^*(C)$ is injective, and $H^*(C_t)$ is torsion-free
over $\Lambda_{\geq0}$.
\end{lemma}

\begin{proof}
First, $H^*(D_t)\to H^*(D)$ is injective. If a nonzero cycle
$z\in D_t$ is a boundary in $D$, zero boundary depth gives, for
every $\epsilon>0$, a primitive $y$ with
$\nu(y)\geq\nu(z)-\epsilon$. Choose
$\epsilon<\nu(z)+t$. Then $y\in D_t$, as required.

Since $C_t=\bigcup_{a>-t}G^aC$ and $D_t=\bigcup_{a>-t}G^aD$,
\eqref{eq:regular-window-recovery} gives $H^*(C_t)\cong H^*(D_t)$.
Inclusion of the lower endpoints $a>-t$ into the full system gives
compatibility with $\Phi$. Thus $H^*(C_t)$ embeds in the
$\Lambda$-vector space $H^*(C)$ and is torsion-free over
$\Lambda_{\geq0}$.
\end{proof}

We will also need the ordinary open lattice for a regular vertex
cap whose corners have been rounded. Let $U$ be such a convex
domain. In the same action--angle chart choose decreasing rational
polytopes $U_n$ with $U\subset\operatorname{int}U_n$ and
$\bigcap_n U_n=U$, all inside one surrounding cotangent domain
$D$. A Hamiltonian strictly negative on $K_U$ is negative on
$K_{U_n}$ for all sufficiently large $n$. Thus, in each finite
ordinary window, the indicator systems for the $K_{U_n}$ together
compute that for $K_U$. Indeed, truncation gives the algebraic
telescope, whose cohomology is the filtered colimit. This holds in both $M$
and $\widehat D$. Since the regular window comparisons commute
with restriction and with the endpoint maps, they give the same
comparisons for $K_U$.

The intrinsic supported complex of $K_U$ has boundary depth zero
by \cite[Corollary~5.7]{GromanBoundaryDepth}, which applies to
every compact convex base. In the exact completion the change of orbit basis
$a_\gamma\mapsto T^{-A_{H,\theta}(\gamma)}a_\gamma$ identifies
the intrinsic action norm after Novikov extension with the ordinary
norm of the weighted telescope
\cite[Section~2.5 and Lemma~5.5]{GromanBoundaryDepth}.
Its boundaries are closed by
\cite[Lemma~5.4]{GromanBoundaryDepth}. Lemma~\ref{lem:ordinary-open-lattice}
now shows that the ordinary open lattice of the rounded cap
injects into its ambient cohomology and is torsion-free.

\begin{proof}[Proof of Theorem~\ref{thm:local-regular-reconstruction}]
Apply Lemma~\ref{lem:regular-localization-windows} to the relative
lattices and the compatible identifications of
Proposition~\ref{prop:regular-window-locality}. This gives
isomorphisms of all full ordinary window cohomologies, compatible
with both endpoints and scalar shifts.

The action--angle identification extends along the Liouville flow
to an exact symplectic identification of $\widehat D$ with $T^*L$.
The intrinsic calculation
\cite[Proposition~5.2 and Theorem~5.6]{GromanVarolgunes2022}
therefore gives isometric BV isomorphisms
\[
 \Psi_Q:SH_{\widehat D}^*(K_Q;\Lambda)
       \xrightarrow{\ \sim\ }\cA^*(Q),
\]
compatible with restriction, with target operator
$\diver_{\Omega_0}$. After translating the cotangent coordinates,
the reference's winding label is $-m$, its character is $z^m$,
and its degree-$q$ coefficient comparison carries a factor
$(-1)^q$. Multiplying that comparison by $(-1)^q$ and reversing
the Floer BV sign gives the conventions fixed in
Section~\ref{subsec:regular-leading-terms}; the two changes
compensate because BV lowers degree by one. Thus the ordinary
cohomology coefficient is sent to the corresponding logarithmic
polyvector. Changes of the cotangent origin are absorbed in the
action normalization of the characters.

In particular, the intrinsic quotient seminorm is Hausdorff.
Lemma~\ref{lem:regular-window-completion} now gives isometric
isomorphisms
\[
 \Phi_{J,Q}:SH_M^*(K_Q;\Lambda)
       \xrightarrow{\ \sim\ }SH_{\widehat D}^*(K_Q;\Lambda).
\]
The linear naturality in that lemma gives compatibility with
restriction. We next verify multiplicativity for a specified pair
of homogeneous classes $x,y$. Write $\Phi=\Phi_{J,Q}$.

Choose ambient cocycles for $x,y$ of ordinary depths at least
$a,c$, respectively. Their product has depth at least $a+c$ by
\eqref{eq:regular-operation-energy}. By the surjection in
\eqref{eq:regular-window-milnor}, choose intrinsic cocycles for
$\Phi x,\Phi y,\Phi(x*y)$ at these same lower bounds, realizing
the images of the three compatible systems of window classes.
The difference between the third cocycle and the product of the
first two represents
\[
 w=\Phi(x*y)-\Phi(x)*\Phi(y).
\]
If $w\ne0$, its quotient valuation $\nu_H(w)$ is finite by
Hausdorffness. Take $B>\nu_H(w)$. The difference is already
nonzero in the intrinsic window $[a+c,B)$. Otherwise a primitive
in that quotient would lift to a cochain whose coboundary could
be subtracted to give a representative of $w$ of depth at least
$B$, contrary to the choice of $B$.

For either ordinary telescope $C$, the product in this output
window is determined by the operation
\begin{equation}\label{eq:regular-window-product}
 H^*(C_{[a,B-c)})\otimes_{\Lambda_0}H^*(C_{[c,B-a)})
       \longrightarrow H^*(C_{[a+c,B)}).
\end{equation}
Indeed, changing the first input in depth at least $B-c$, or the
second in depth at least $B-a$, changes the product only in depth
at least $B$. The operation and its continuation homotopies have
this filtration bound by \eqref{eq:regular-operation-energy}.

Let $\bar x,\bar y$ be the two input window classes.
By Lemma~\ref{lem:regular-localization-windows}, after multiplying
them by $T^\lambda,T^\mu$ for suitable $\lambda,\mu\geq0$, we can
represent them in the nonnegative relative lattice. The input
windows shift by $\lambda$ and $\mu$, and the output window shifts
by $\lambda+\mu$. By \eqref{eq:rc-completed-window} and
Lemma~\ref{lem:rc-regular-companion}, the two shifted classes have
representatives at finite strict Hamiltonian stages. Continue both
to a common input stage and choose a later output Hamiltonian
dominating the two weighted inputs.

We may use the data of Proposition~\ref{prop:regular-operation-confinement}
for this finite diagram, by the continuation compatibility of
Proposition~\ref{prop:regular-window-locality}. Ambient confinement
and the intrinsic contact cut identify the oriented pair-of-pants
counts on these representatives. The comparison with regular
companions above identifies the strict counts with the actual
products. The input and output identifications in this calculation
are precisely the window maps used to construct $\Phi$. Passing
to the full ordinary output window and dividing by $T^{\lambda+\mu}$
therefore gives equality of the two product classes in $[a+c,B)$.
This contradicts the choice of $B$ and proves multiplicativity.

For BV, a nonzero difference $\Phi(\Delta_F x)-\Delta_F(\Phi x)$
would likewise survive in a window $[a,B)$. The framed-cylinder
operation preserves ordinary depth, so only the input in $[a,B)$
is needed. After a single Novikov shift choose a strict input stage
representing that class and an output dominating all its time
rotations. The same confinement and companion comparisons identify
the framed-cylinder counts in the shifted window, giving the
contradiction. Thus $\Phi$ intertwines the BV operators.

It follows that $\Phi_Q=\Psi_Q\circ\Phi_{J,Q}$ is an isometric BV
isomorphism, and in particular preserves the Gerstenhaber bracket.

Proposition~\ref{prop:regular-window-locality} identifies the
ambient and intrinsic normalized local orbit blocks with their
cohomology coefficients and restriction weights. The intrinsic
calculation sends these blocks to the normalized monomial symbols
$\vartheta_{Q,m}\otimes\alpha$. Consequently $\Phi_Q$
intertwines the Floer and analytic tropical markings, as required.
\end{proof}

\subsection{Assembly and canonical gluing}
\label{subsec:canonical-affinoid}

\begin{proof}[Proof of the gluing assertions in
Theorem~\ref{thm:global-reconstruction}]
For regular convex polytopes $Q\subset P$, let
$\iota_{QP}:Y_Q\to Y_P$ be induced by actual Floer restriction.
The local theorem identifies this map with
$\Log^{-1}(Q)\hookrightarrow\Log^{-1}(P)$: an inequality
$\langle m,b\rangle\geq c$ defining $Q$ becomes the Weierstrass
inequality $|T^{-c}z^m|\leq1$. Thus $\iota_{QP}$ is an affinoid
domain embedding, and
\[
 \iota_{RP}=\iota_{QP}\circ\iota_{RQ}\qquad(R\subset Q\subset P).
\]
All these identifications hold simultaneously in one surrounding
cotangent domain for $P$.

Let $p_P:Y_P^{\rm Berk}\to P$ be the local valuation map. The torus model gives
\begin{equation}\label{eq:regular-gluing-domains}
 p_P^{-1}(Q)=\iota_{QP}^{\rm Berk}(Y_Q^{\rm Berk})
\end{equation}
and, for $Q,R\subset P$,
\begin{equation}\label{eq:regular-gluing-intersections}
 \iota_{QP}(Y_Q)\cap\iota_{RP}(Y_R)
 =\iota_{Q\cap R,P}(Y_{Q\cap R}).
\end{equation}
Finite polyhedral coverings of $P$ likewise pull back to admissible
affinoid coverings.
Since convex subpolytopes separate points of $P$,
\eqref{eq:regular-gluing-domains} determines $p_P$ from the actual
restriction maps. It is therefore independent of the chosen torus
coordinates and satisfies $p_P\circ\iota_{QP}^{\rm Berk}=p_Q$.

Two polytopes $P,Q$ need not lie in one surrounding chart. Cover
$P\cap Q$ by finitely many small polytopes $R_a$ convex in both
charts. Such a covering is obtained by intersecting $P$ and $Q$ with
small rational boxes in common coordinate neighborhoods: on each
connected overlap the coordinate change is integral affine.
Compactness gives a finite covering, and the construction also gives
common refinements. Set
\[
 U_{PQ}=\bigcup_a\iota_{R_aP}(Y_{R_a}),
 \qquad U_{PQ}^{\rm Berk}=p_P^{-1}(P\cap Q).
\]
This is an admissible open subset with the displayed finite affinoid
covering, independent of the refinement. The identifications between
the images of each $Y_{R_a}$ in $Y_P$ and $Y_Q$ agree on intersections by
\eqref{eq:regular-gluing-intersections} and functoriality of Floer
restriction. They therefore define an analytic isomorphism
\[
 \psi_{PQ}:U_{PQ}\xrightarrow{\sim}U_{QP},
 \qquad \psi_{PQ}\circ\iota_{RP}=\iota_{RQ}
\]
for every small convex $R\subset P\cap Q$. In particular,
$p_Q\circ\psi_{PQ}^{\rm Berk}=p_P$.

For three polytopes the overlap domains satisfy
\[
 \psi_{PQ}(U_{PQ}\cap U_{PR})=U_{QP}\cap U_{QR},
\]
because their Berkovich realizations are the corresponding inverse
images of $P\cap Q\cap R$. On this domain,
$\psi_{QR}\circ\psi_{PQ}=\psi_{PR}$: after a finite common
polyhedral refinement, both maps are induced by the same actual
restrictions to each small polytope. The rigid gluing theorem
\cite[\S 5.3, Proposition~5]{Bosch14} now gives a space $Y_{\rm reg}$
with these charts and transitions. It is smooth, since the charts
are torus affinoids.

The valuation maps agree on overlaps and give a continuous projection
$p:Y_{\rm reg}^{\rm Berk}\to B_{\rm reg}$. A point of $Y_P^{\rm Berk}$
over $Q$ lies in $U_{PQ}^{\rm Berk}$, where the transition identifies
it with a point of $Y_Q^{\rm Berk}$. Thus $p^{-1}(Q)=Y_Q^{\rm Berk}$.
In torus coordinates the rigid point $(T^{b_1},T^{b_2})$ lies over
$(b_1,b_2)$, giving surjectivity.

These affinoid charts also give properness. Points with different
images under $p$ can be separated in the base. Two points over the
same $b$ lie in the open set $p^{-1}(\operatorname{int}Q)$ for some
polytope $Q$ containing $b$ in its interior. This set lies in the
Hausdorff affinoid $Y_Q^{\rm Berk}$, so the realization is Hausdorff.
For a compact $K\subset B_{\rm reg}$, choose finitely many polytopes
$Q_i$ whose interiors cover $K$. Then $p^{-1}(K)$ is closed in the
compact union $\bigcup_iY_{Q_i}^{\rm Berk}$ and hence is compact.

The construction is independent of the regular cover. For any
locally finite cover by small polytopes whose interiors cover the
base, each compact regular polytope has a finite refinement subordinate
to that cover. Its affinoids form an admissible covering of $Y_Q$,
so gluing just the chosen charts yields the same space. Coordinate
independence follows because all embeddings and transitions were
defined by actual restriction.

By construction,
\[
 \Gamma(Y_Q,\cO_{Y_{\rm reg}})=\cF^0(Q),
\]
naturally under restriction. The higher-degree comparisons also glue:
an analytic polyvector is determined by its iterated brackets with
functions, so two local Gerstenhaber identifications agree wherever
their degree-zero maps agree. We obtain the marked BV comparisons
on the domains $Y_Q$, compatible with every inclusion of regular
polytopes.
\end{proof}

In particular, the regular convex strips used later have their
standard torus affinoids. Their injective affine development and
compact regular neighborhood place them in a surrounding action--angle
domain, regardless of their length.

\subsection{The regular volume form}

Use the chosen orientation of $B_{\rm reg}$ to order the two
characters in each marked chart. We show that the logarithmic
volume form in these coordinates is independent of the chart.

\begin{proof}[Proof of the volume-form assertions in
Theorem~\ref{thm:global-reconstruction}]
On a regular polytope small enough for two local comparisons from
Theorem~\ref{thm:local-regular-reconstruction}, write $u_1,u_2$ and
$v_1,v_2$ for their coordinates and set
\[
 \Omega_u=d\log u_1\wedge d\log u_2,\qquad
 \Omega_v=d\log v_1\wedge d\log v_2=g\Omega_u.
\]
Both comparisons identify Floer BV with the corresponding logarithmic
divergence. Their comparison on polyvectors is induced by the analytic
change of coordinates, as in the gluing proof of
Theorem~\ref{thm:global-reconstruction}.
Thus $\diver_{g\Omega_u}X=\diver_{\Omega_u}X$ for each vector field
$X$, which gives $X(g)=0$. Applying this to
$u_j\partial_{u_j}$ in the Laurent expansion shows that
$g\in\Lambda^\times$, since $\kk$ has characteristic zero.
This is the volume comparison of
\cite[Proposition~3.6]{GromanVarolgunes2022}.

To determine the scalar, restrict to the torus domain over one
point of the polytope. The tropical marking identifies the reductions
of the normalized characters, so on this domain
\[
 v_i=c_i u_1^{m_{i1}}u_2^{m_{i2}}(1+h_i),\qquad
 c_i\in\Lambda^\times,\quad \|h_i\|<1,
\]
where $M=(m_{ij})\in\mathrm{GL}_2(\bZ)$ is the change of basis on
the character lattice. Since the valuation is trivial on $\kk$,
the logarithms $L_i=\log(1+h_i)$ converge. With
$\alpha_i=\sum_j m_{ij}\,d\log u_j$ we have
$d\log v_i=\alpha_i+dL_i$, and hence
\[
 (g-\det M)\Omega_u
 =d\bigl(L_1\alpha_2-L_2\alpha_1+L_1dL_2\bigr).
\]
The coefficient of $u_1^0u_2^0$ in an exact two-form in the frame
$\Omega_u$ is zero. Indeed, its coefficient function is of the form
$u_1\partial_{u_1}b-u_2\partial_{u_2}a$ for Laurent series $a,b$.
It follows that $g=\det M$. For oriented coordinates $\det M=1$,
so the local forms agree and give $\Omega_{\rm reg}$.

The local BV formula now gives the asserted global divergence.
In oriented coordinates the dual bivector is
$z_1\partial_{z_1}\wedge z_2\partial_{z_2}$, so it generates the
top polyvectors and the marking identifies its zero Laurent term
with the normalized orientation class of the fiber. Its contraction
with $\Omega_{\rm reg}$ is $1$, hence its divergence vanishes.
Naturality of the regular comparison gives compatibility of the
classes $s_Q$ under restriction.
\end{proof}

\begin{proof}[Proof of
Theorem~\ref{thm:intro-regular-reconstruction}]
The local comparisons and their compatibility with restriction and
the tropical markings are given by
Theorem~\ref{thm:local-regular-reconstruction}. In degree zero they
identify $\cF^0(Q)$ with the affinoid algebra of a domain in the dual
analytic torus. Theorem~\ref{thm:global-reconstruction} gives the
gluing of these affinoids, the projection to $B_{\rm reg}$, the
volume form and the global compatibility of the BV comparisons.
\end{proof}

For a polygon containing a singular value we can still use the
regular comparison near its vertices. In Section~\ref{sec:detection}
we ask whether restriction to these neighborhoods detects a class
over the polygon and its norm.

\section{Detection by vertex caps}
\label{sec:detection}\label{SecSYZ}\label{SecIntAffGeom}\label{Sec4DimSYZ}

Let $P$ be a polygon whose boundary lies in the regular locus of
the torus fibration $\pi:M\to B$. By Section~\ref{sec:regular}, over
a small convex polygon at each vertex we can identify Floer cohomology
with analytic polyvector fields on the corresponding torus domain.
When $P$ contains a singular value, we will use restriction to these
vertex caps to detect classes in $\cF^j(P)$ and compute their norms.
The local convexity and Delzant conditions below are those of
\cite[\S7.1.1]{GromanIndicatorV1}.

\begin{df}[Admissible polygon]
\label{def:boundary}
Let $P\subset B$ be a compact topological disk containing at most
one singular value in its interior. We say that $P$ is an
\emph{admissible polygon} if the following conditions hold.
\begin{enumerate}
\item $\partial P\subset B_{reg}$.

\item The boundary of $P$ consists of finitely many integral affine
segments of rational direction and is locally convex. More
precisely, orienting $\partial P$ counterclockwise in positively
oriented local affine charts, at each vertex the tangent vector
to the outgoing edge points to the left of the line through the
incoming edge. Moreover, the polygon is Delzant: the primitive
integral covectors annihilating the two edges at each vertex
form a basis of the integral cotangent lattice.

\item The integral affine structure is conical near $\partial P$.
That is, there is an affine atlas on a regular annular neighborhood
of $\partial P$ whose transition maps preserve the Euler vector
field
\[
 \mathcal E=\sum_i p_i\frac{\partial}{\partial p_i},
\]
and this vector field points outward along every edge of $P$.
\end{enumerate}
An admissible polygon is called \emph{regular} if it contains no
singular value. In this case we also fix a conical center in its
convex affine interior.
\end{df}

For reconstruction across an isolated singular value $s$ we use
three conditions on a polygonal neighborhood of $s$.

\begin{definition}[Boundary geometry $(G_\partial)$]
\label{def:singular-boundary-geometry}
A polygon $P$ satisfies $(G_\partial)$ if it is admissible and has
at least three labeled edges with distinct adjacent labeled pairs.
At a singular value $s$ we require that $P$ contain $s$ in its interior
and no other singular value.
\end{definition}

\begin{definition}[Detection in degree $j$]
\label{def:topological-detection}
For an admissible polygon $P$, put $K_P=\pi^{-1}(P)$ and let $L$
be a regular fiber in its boundary collar. We say that $(T)_j$ holds
if the ordinary restriction
\[
 H^j(K_P;\kk)\longrightarrow H^j(L;\kk)
\]
is injective. The condition can also be formulated with $\kk$
replaced by the coefficient ring of the Floer calculation.
\end{definition}

The condition $(T)_j$ is independent of the choice of a polygonal
neighborhood containing the given singular value and no others.
Indeed, given two polygonal disks containing the same singular value
and no others, choose a smaller disk about that value contained in
both interiors. The annular regions between the respective
boundaries lie in the regular torus bundle. Their retractions onto
the inner boundary lift to retractions of the corresponding fibered
neighborhoods, and transport identifies the restrictions to nearby
regular fibers.

\begin{definition}[Full radial conicity $(G_{\rm rad})$]\label{def:radial}
The conical structure of the boundary collar extends over the
punctured neighborhood $P\setminus\{s\}$, and its inward flow gives
embedded homothetic polygons $P_r$, $0<r\leq1$, with $P_1=P$ and
$P_{r'}\subset\operatorname{int}P_r$ for $r'<r$. Their boundaries
are regular, with the same labeled edge directions and Delzant pairs.
They form a neighborhood basis of $s$ and $\bigcap_{r>0}P_r=\{s\}$.
\end{definition}

The Liouville lifts of the Euler field on regular collars are
described in Lemma~\ref{lmContactBoundary}. We use the radial
condition only at positive radii, in Section~\ref{sec:radial}.

At a vertex $v$ of an admissible polygon, take the primitive inward
affine coordinates $x,y$ vanishing on its two incident edges. By the
Delzant condition these are integral affine coordinates. For
$\epsilon>0$ sufficiently small, the triangle
$Q_v=\{x,y\ge0,\ x+y\le\epsilon\}$ in this chart is a regular
Delzant polygon. We call it a \emph{vertex cap}.

Each $K_{Q_v}$ retracts onto a fiber, and transport in the collar
identifies the kernels of restriction from $H^j(K_P;\kk)$ to the
cohomology of the different fibers. Thus one can equally test $(T)_j$
by restriction to the vertex caps. For a regular convex polygon
$K_P$ retracts onto a fiber, so $(T)_j$ holds in every degree.

We give $\cF^j(P)$ the quotient seminorm induced by the ordinary
Novikov norm on $SC_M^*(K_P;\Lambda)$. Denote by
\[
 \cF^j_0(P)
 =H^j\!\bigl(\{x\in SC_M^*(K_P;\Lambda):\ \|x\|<1\}\bigr)
\]
the ordinary open lattice in degree $j$.

\begin{theorem}[Strict detection]\label{thm:strict-detection}
Let $P$ be an admissible polygon and choose a
sufficiently small cap $Q_v$ at each vertex. In every degree $j$ in
which $(T)$ holds, Floer restriction
\[
 R_P^j:\cF^j(P)\longrightarrow\bigoplus_{v\in V(P)}\cF^j(Q_v)
\]
is an isometric embedding with closed image for the maximum norm
on the right. The ordinary open lattice $\cF^j_0(P)$ is
torsion-free over $\Lambda_{\geq0}$.
\end{theorem}

We first round the corners of $P$ and its vertex caps. At a fixed
period cutoff, the local Floer calculation gives injectivity in
sufficiently short Novikov windows. We then enlarge the windows,
remove the cutoff, and pass from the roundings to $P$.

\subsection{Rounding the polygon and its caps}\label{sec:polygon-windows}

Only the boundary of a rounding needs to lie in the regular locus.
We describe its convexity in affine charts along the boundary.

\begin{df}[Integral affine convex boundary]\label{dfIAConvex}
Let $D\subset B$ be a compact domain with $C^1$ boundary.
We say that $\partial D$ is \emph{integral affine convex} if
$\partial D\subset B_{reg}$ and for every $p\in\partial D$
there exists an integral affine chart
$\varphi:V\to\bR^2$ with $p\in V\subset B_{reg}$ such that
$\varphi(D\cap V)$ is contained in a closed affine half-plane
whose boundary line is tangent to
$\varphi(\partial D\cap V)$ at $\varphi(p)$.
\end{df}

\begin{rem}
The definition is independent of the chosen chart because
$Aff(\bZ^2)$ sends lines to lines and half-planes to half-planes,
preserving tangency.
\end{rem}

\begin{df}\label{dfAdmissible}
A compact domain $D\subset B$ is \emph{admissible} if:
\begin{enumerate}
\item $\partial D\subset B_{reg}$,
\item $\partial D$ is a smooth embedded circle and is convex
with respect to the integral affine structure in the sense of
Definition~\ref{dfIAConvex},
\item The integral affine structure is conical near $\partial D$:
there exists an atlas where the Euler vector field
$\sum p_i\frac{\partial}{\partial p_i}$ is preserved by transition
maps and points outward from $D$.
\end{enumerate}
\end{df}

The conical structure gives a flow transverse to the boundary.
We use its lift to the torus fibration to obtain the contact
geometry required by the smooth-domain Floer theory.

\begin{lemma}[Contact collar]\label{lmContactBoundary}
The conical boundary collar lifts to a Liouville collar of $K_P$.
\end{lemma}
\begin{proof}
A torus bundle over an annulus admits a section, since its base has
the homotopy type of a circle and its fiber is connected. The pullback
of $\omega$ by a smooth section is exact, as $H^2$ of the annulus
vanishes. Translating the section by a one-form on the base whose
differential cancels this pullback gives a Lagrangian section.
Use it as the zero section in action--angle coordinates. The conical
base changes are linear, and the angle coordinates then transform by
the inverse transpose. Thus the local forms $\theta=\sum p_j\,dq_j$
agree. Their Liouville field lifts the Euler field and is transverse
to each boundary face. For any smooth rounding $D$ in the collar
transverse to the Euler field, $\theta|_{\partial K_D}$ is a contact form.
Changing the Lagrangian section identifies the resulting collars by
a fiberwise symplectic translation.
\end{proof}

We choose the approximation so that each edge still has a flat
portion on the rounded boundary.

\begin{proposition}[Conical smoothing]
  \label{prpCornerSmoothing}
  There is a nested outward family $P_t$, $0\le t<\varepsilon$,
  with $P_0=P$, whose boundary for $t>0$ is smooth and locally
  affine-convex. Outside $O(t)$ vertex neighborhoods its facets
  are the conical dilates of those of $P$; inside them the
  adjacent facets are joined by smooth convex rounded arcs.
  \end{proposition}
  
  \begin{proof}
  Let $\mathcal E$ be the Euler vector field near $\partial P$.
  For each edge $E_i$ choose an affine defining function $\ell_i$,
  positive on the exterior side of $E_i$, normalized so that
  \[
  \mathcal E\ell_i=\ell_i+1.
  \]
  Such a normalization is possible since $\mathcal E$ points
  outward along $E_i$.
  
  Choose an even, nonnegative smooth function $\rho$ of integral
  one, supported on $[-1/2,1/2]$ and positive on its interior.
  For $t>0$ put
  \[
  m_t(a,b)=\int_{\mathbb R}
  \max\{a+tu,b-tu\}\rho(u)\,du,
  \qquad
  m_0(a,b)=\max\{a,b\}.
  \]
  The function $m_t$ is smooth, convex and nondecreasing in each
  variable. It satisfies
  \[
  \max\{a,b\}\le m_t(a,b)\le\max\{a,b\}+t/2,
  \]
  and agrees with $\max\{a,b\}$ whenever $|a-b|\ge t$.
  Moreover,
  \[
  |m_t(a,b)-m_s(a,b)|\le |t-s|/2.
  \]
  
  Near the interior of an edge $E_i$ define $P_t$ by
  $\ell_i\le t$, and near a vertex $E_i\cap E_{i+1}$ define it by
  \[
  m_t(\ell_i,\ell_{i+1})\le t.
  \]
  Choose the edge and vertex neighborhoods so that on their
  overlaps the defining function of the other incident edge
  is strictly smaller. For sufficiently small $t$ the two
  definitions therefore agree on the overlaps and define
  $P_t$ in a fixed collar of $\partial P$. Extend $P_t$ by
  including the interior of $P$ beyond this collar.
  
  The estimates above give $P\subset P_t$, and for $s<t$,
  \[
  m_t(a,b)-t
  \le m_s(a,b)-s-\tfrac12(t-s).
  \]
  Thus $P_s\subset P_t$. The same estimates show that
  the domains decrease to $P$ as $t\downarrow0$.
  
  The boundary is locally convex since its defining functions
  are convex in affine coordinates. Furthermore, near each
  vertex, $\mathcal E\ell_i$ and $\mathcal E\ell_{i+1}$ are
  positive. Since the partial derivatives of $m_t$ are
  nonnegative and sum to one,
  \[
  \mathcal E m_t(\ell_i,\ell_{i+1})>0.
  \]
  The same holds for the defining function along each edge.
  Consequently the boundary is smooth and transverse to the
  Euler field. It is a graph over $\partial P$ along the Euler
  flow, and hence bounds a disk.
  
  Away from the vertices the boundary is given by $\ell_i=t$,
  which is the image of $\ell_i=0$ under the Euler flow for
  time $\log(1+t)$.
  
  Finally, in affine coordinates centered at a vertex, the
  identity
  \[
  m_t(ta,tb)=t\,m_1(a,b)
  \]
  shows that the rounding is obtained from a fixed profile
  by dilation by $t$. The rounded arc lies in an $O(t)$
  neighborhood of the vertex and has positive curvature
  between its flat ends: indeed,
  \[
  m_t(a,b)=b+f_t(a-b),
  \qquad
  f_t''(r)=\frac{1}{2t}\rho\left(-\frac{r}{2t}\right)>0
  \quad\text{for } |r|<t,
  \]
  and the two incident defining functions have independent
  differentials.
  \end{proof}

\begin{figure}[htbp]
\centering
\begin{tikzpicture}[x=1.16cm,y=1.16cm,font=\small,
    sharp/.style={black!75,thin,densely dashed},
    front/.style={black!85,line width=.75pt},
    common/.style={teal!65!black,line width=1.8pt,line cap=round}]
  \def\corner{(-.28,.18) .. controls (-.28,-.10) and (-.10,-.28)
                                      .. (.18,-.28)}
  \def\outer{\corner -- (3.02,-.28)
    .. controls (3.30,-.28) and (3.48,-.10) .. (3.48,.18) -- (3.48,3.02)
    .. controls (3.48,3.30) and (3.30,3.48) .. (3.02,3.48) -- (.18,3.48)
    .. controls (-.10,3.48) and (-.28,3.30) .. (-.28,3.02) -- cycle}
  \def\capfront{\corner -- (2.22,-.28)
    .. controls (2.42,-.28) and (2.46,-.17) .. (2.33,-.04)
    -- (-.04,2.33)
    .. controls (-.17,2.46) and (-.28,2.42) .. (-.28,2.22) -- cycle}
  \begin{scope}
    \fill[black!8] \outer;
    \fill[white] (0,0) rectangle (3.2,3.2);
    \draw[front] \outer;
    \draw[sharp] (0,0) rectangle (3.2,3.2);
    \node at (1.60,1.60) {$P$};
    \node at (1.60,-.74) {(a) Outward rounding};
  \end{scope}
  \begin{scope}[xshift=5.65cm]
    \fill[teal!13] \capfront;
    \draw[front] \outer;
    \draw[front] \capfront;
    \draw[common] (-.28,2.22) -- \corner -- (2.22,-.28);
    \draw[sharp] (0,0) rectangle (3.2,3.2);
    \draw[sharp] (0,2) -- (2,0);
    \fill (0,0) circle (1.2pt);
    \node[above right,inner sep=2pt] at (0,0) {$v$};
    \node at (.72,.72) {$Q_v$};
    \node at (2.42,2.34) {$P$};
    \node at (1.60,-.74) {(b) Compatible rounding of a cap};
  \end{scope}
\end{tikzpicture}
\caption{Dashed lines mark the original edges. In (a) the shaded collar
is the region added to $P$ to obtain $P_t$. In (b) the shaded domain
is $Q_{v,t}$. It contains $Q_v$ and lies in $P_t$. The heavy line is
common to both rounded boundaries, including the entire arc at $v$.}
\label{fig:outward-compatible-rounding}
\end{figure}
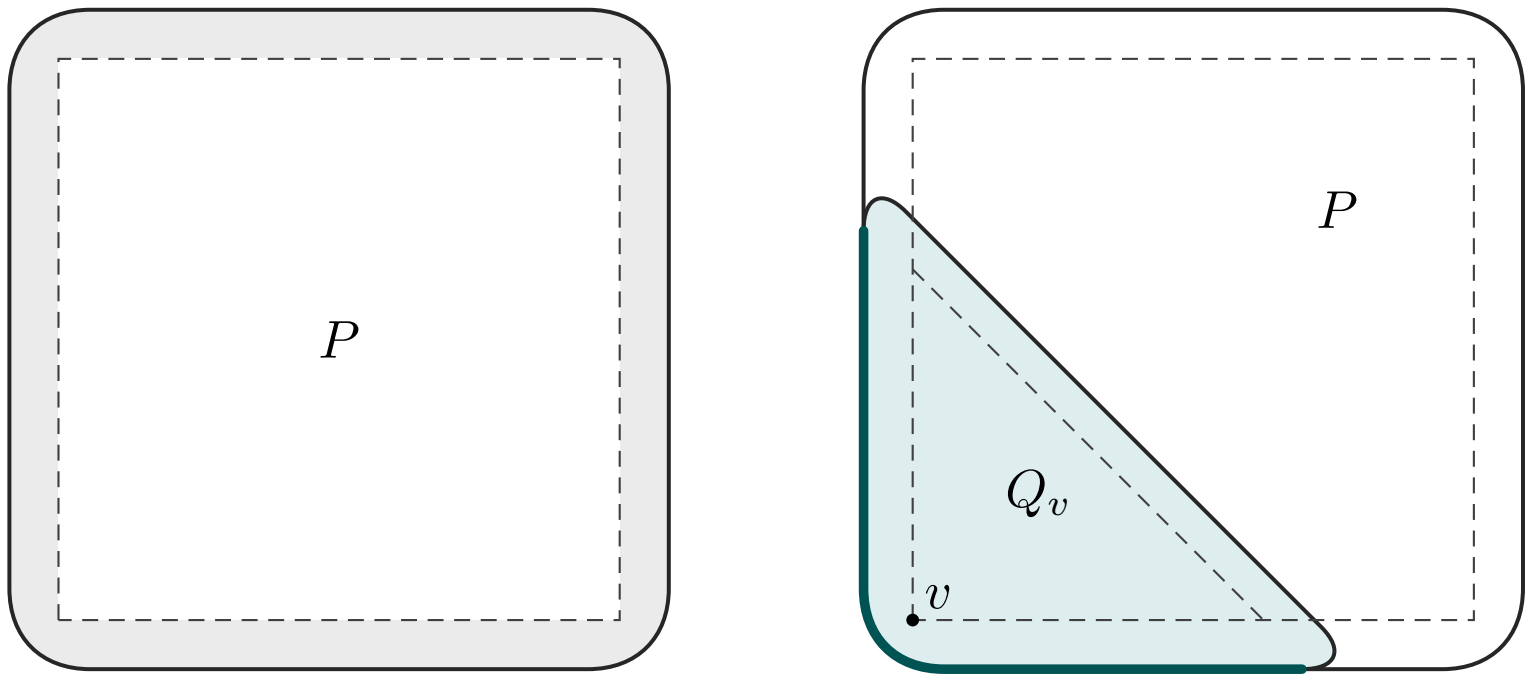

To round the caps we need the same construction for an inclusion
$Q\subset P$. Where $Q$ and $P$ share a vertex, we want their
roundings to agree on the whole rounded arc. Along an edge of $P$
the common boundary may consist of a proper subsegment.
Figure~\ref{fig:outward-compatible-rounding} shows the rounding of $P$
and of a cap at a shared vertex.

We say that an inclusion $Q\subset P$ of admissible polygons is
\emph{admissible}, or \emph{facewise admissible}, if the common
boundary $\partial Q\cap\partial P$ is a union of nondegenerate
closed intervals whose endpoints lie in the interiors of edges
of $P$, with the empty and whole-boundary cases also allowed.
Furthermore, its intersection with each edge of $P$ must be
either empty or a single nondegenerate segment.

\begin{figure}[htbp]
 \centering
 \begin{tikzpicture}[scale=0.95]
   \begin{scope}
     \fill[gray!35] (0,0)--(1.8,0)--(0,1.8)--cycle;
     \draw[thick] (0,0) rectangle (3.2,3.2);
     \draw[thick] (1.8,0)--(0,1.8);
     \draw[very thick] (1.8,0)--(0,0)--(0,1.8);
     \node at (0.55,0.55) {$Q$};
     \node at (2.35,2.35) {$P$};
     \node[below] at (1.6,-0.3) {Admissible};
   \end{scope}
   \begin{scope}[xshift=5cm]
     \fill[gray!35]
       (0,0)--(1.7,0)--(2.7,1)--(1,1)--cycle;
     \draw[thick] (0,0) rectangle (3.2,3.2);
     \draw[thick] (1.7,0)--(2.7,1)--(1,1)--(0,0);
     \draw[very thick] (0,0)--(1.7,0);
     \node at (1.35,0.5) {$Q$};
     \node at (2.35,2.35) {$P$};
     \node[below] at (1.6,-0.3) {Inadmissible};
   \end{scope}
 \end{tikzpicture}
 \caption{On the left both edge germs at the shared vertex are
 retained. On the right the common boundary terminates at a
 vertex of the larger polygon.}
 \label{FigAdm}
\end{figure}
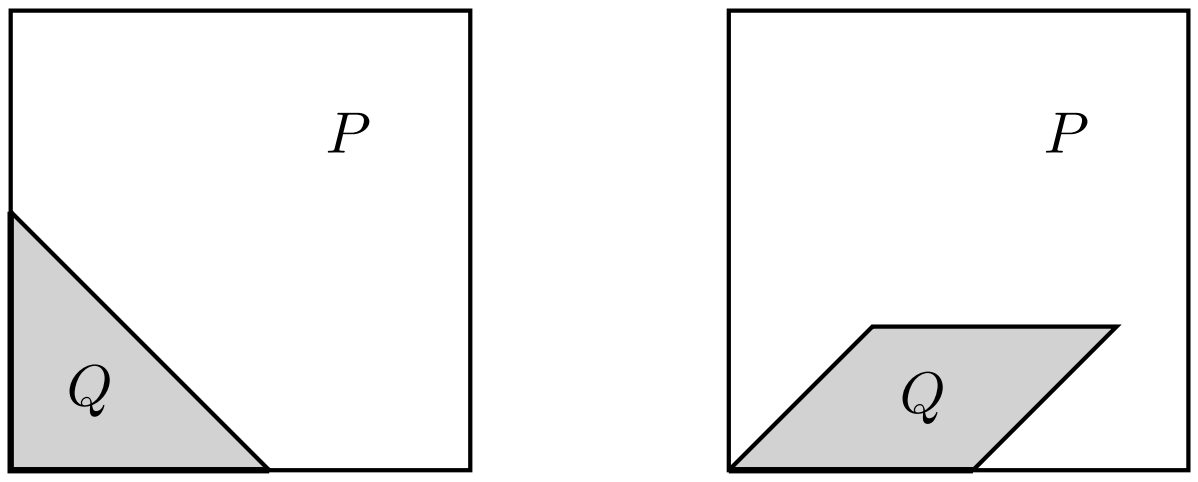

For example, the inclusion of a vertex cap $Q_v$ is admissible. Its
common boundary with $P$ consists of the two edge segments issuing
from $v$, as on the left of Figure~\ref{FigAdm}.
More generally, if a vertex of $P$ belongs to $\partial Q$, both incident
edge germs must be shared. The requirement on the endpoints
concerns the connected components of the common boundary,
rather than each individual shared edge. Each such component
is a union of successive edges of $Q$.

\begin{proposition}[Compatible smoothings]
  \label{prpCompatibleCornerSmoothing}
  A facewise admissible inclusion $P_1\subset P_2$ admits
  nested outward smoothings $P_{1,t}\subset P_{2,t}$ with
  smooth locally convex regular boundaries, $O(t)$ facet
  displacements, and fixed local rounding profiles at scale
  $t$. On a neighborhood of every compact subarc of the
  relative interior of their common boundary, the smoothings
  agree for sufficiently small $t$. In particular, the
  entire rounded arc at every shared vertex is common to
  both boundaries. The smoothing of $P_2$ may be chosen
  as in Proposition~\ref{prpCornerSmoothing}.
  \end{proposition}
  
  \begin{proof}
  Choose the smoothing of $P_2$ as in
  Proposition~\ref{prpCornerSmoothing}.
  For each edge of $P_1$ contained in an edge of $P_2$,
  use the defining function of the corresponding edge of
  $P_2$. For the remaining edges of $P_1$, choose affine
  defining functions positive on the exterior side.
  
  Apply the construction of
  Proposition~\ref{prpCornerSmoothing} to these defining
  functions, using the same regularized maximum $m_t$.
  The proof of outwardness, nestedness, smoothness and
  convexity is unchanged. In particular, normalization of
  the defining functions with respect to the Euler field
  of $P_1$ is not needed: outward transversality follows
  from positivity of its derivative on each incident
  defining function in a sufficiently small neighborhood
  of the boundary.
  
  It remains to check containment and agreement near the
  common boundary. At a shared vertex of $P_2$, admissibility
  implies that both incident edges are shared. The two
  smoothings are therefore given by the same inequality
  \[
  m_t(\ell_i,\ell_{i+1})\le t.
  \]
  They agree on a fixed neighborhood of the vertex, including
  the entire rounded arc for small $t$.
  
  Along the interior of a shared edge segment, away from
  its endpoints, both smoothings are given by $\ell_i\le t$.
  
  The only remaining points to consider are the endpoints
  of the common boundary intervals. By admissibility these
  lie in the interiors of edges of $P_2$. Near such a point
  the smoothing of $P_2$ is given by $\ell_i\le t$, while
  that of $P_1$ is given by
  \[
  m_t(\ell_i,h)\le t,
  \]
  where $h$ is the defining function of the other incident
  edge of $P_1$. Since
  \[
  m_t(\ell_i,h)\ge\ell_i,
  \]
  the latter domain is contained in the former.
  
  Away from a neighborhood of the common boundary, the
  boundary of $P_1$ lies a positive distance inside $P_2$.
  Thus, by compactness, for sufficiently small $t$ we have
  \[
  P_1\subset P_{1,t}\subset P_{2,t}.
  \]
  This proves the claimed compatibility.
  \end{proof}

\subsection{Restriction to a rounded vertex cap}
\label{subsec:rounded-cap-restriction}

Choose vertex caps in disjoint regular neighborhoods. Fix a single
outward family for $P$ and use it to round each cap by
Proposition~\ref{prpCompatibleCornerSmoothing}. At a fixed parameter
write $D$ and $D_v\subset D$ for the resulting domains. Which Reeb
families of $\partial K_D$ survive restriction to $K_{D_v}$ in a
short window?

In action--angle coordinates $(p,q)$, the characteristic direction
over a point of a smooth boundary is determined by the conormal
to the boundary at that point. Along a straight edge this direction
is constant. When a corner is rounded convexly, the conormal
direction sweeps through the cone generated by the conormals of
the two adjacent edges. A characteristic closes when its conormal
spans a rational ray. The primitive integral generator labels the
simple orbits, and its positive multiples label their covers.

Denote the edges of $P$ in cyclic order by
$e_1,\ldots,e_N$, and let $\xi_i$ be the primitive outward integral
conormal to $e_i$. To each edge associate the monoid
\[
 C_{e_i}=\bN\xi_i,
\]
and to the vertex $e_i\cap e_{i+1}$ associate the monoid
\[
 C_{i,i+1}=\bN\xi_i+\bN\xi_{i+1}.
\]
By the Delzant condition, the latter consists of all the integral
points of the corresponding cone.

We define the \emph{integral tropicalization} of $P$ to be
\begin{equation}\label{eqPtropUnion}
 P_{trop}(\bZ)=\bigcup_i C_{i,i+1},
\end{equation}
where the copy of $C_{e_i}$ in $C_{i-1,i}$ is identified with its
copy in $C_{i,i+1}$, and all the origins are identified.
Addition is defined for elements belonging to a common cone.
Each cone uses the integral lattice at its own vertex, and the
affine structure along $e_i$ identifies the two copies of
$C_{e_i}$.

\begin{lemma}\label{lmMxUwIso}
  The Reeb components of $\partial K_D$ are in bijection
  with the nonzero elements of $P_{trop}(\bZ)$.
  Denote by $\gamma_m$ the component corresponding to $m$
  and choose a point $b_m$ in its projection to the base.
  
  For every finite period cutoff $a$ there is a $\delta>0$
  such that for every ordinary Novikov window $I$ of width
  less than $\delta$ there is an isomorphism of graded modules
  \[
   SH^{\leq a}_{M,I}(K_D)
   \simeq
   \left(
   H^*(K_D;\kk)\oplus
   \bigoplus_{\substack{
   m\in P_{trop}(\bZ)\setminus\{0\}\\
   \mathcal T_\theta(\gamma_m)\leq a}}
   H^*(L_{b_m};\kk)
   \right)\otimes_\kk\Lambda_I^{\log}.
  \]
  Under this isomorphism the zero label corresponds to the
  critical block and each nonzero label $m$ corresponds to
  the Reeb block of $\gamma_m$.
  \end{lemma}
  
  \begin{proof}
  In action--angle coordinates the characteristic direction
  over a point of $\partial D$ is determined by the
  conormal to the boundary. Along a flat edge this
  direction is constant, while along a rounded corner
  it sweeps through the interior of the cone generated
  by the two incident conormals. The characteristic
  closes precisely when the conormal spans a rational
  ray. The primitive integral generator labels the
  simple orbits and its positive multiples label their
  covers. This gives the claimed bijection.
  
  There are finitely many components below any finite
  period cutoff. Indeed, the pairing of the Euler field
  with the outward unit conormal is bounded below by a
  positive constant along the boundary. The period of
  an orbit is therefore bounded below by a positive
  multiple of the length of its integral conormal.
  By Theorem~\ref{thm:background-indicator-input}, it
  remains to identify the unweighted local Floer
  cohomology of each component.
  
  For a Reeb component over a rounded corner, evaluation
  at time zero identifies the family of parametrized
  orbits with the fiber $L_{b_m}$.
  Choose action--angle coordinates and write a $k$-fold
  orbit with primitive outward conormal $n$ as $q(t)=q_0-kn t$.
  Indeed, our convention gives $X_H=-\nabla_p H\,\partial_q$.
  In the moving frame the linearized normal model is
  the cotangent-torus model with positive momentum
  Hessian. Its normal Morse index is zero and the
  moving frame has zero Maslov rotation with respect
  to the Maslov-zero grading. Thus the Floer grading
  agrees with the ordinary Morse grading on $L_{b_m}$;
  see also \cite[Section~6]{Paper1}.
  
  To compute the orientation local system, use the spin
  framing induced by the chosen action--angle coordinates.
  As $q_0$ varies, the Fredholm capping operators are
  conjugated by translations preserving this framing.
  Along a loop in $L_{b_m}$ the conjugation returns to
  the identity, so the determinant line has trivial
  monodromy. The space of positive normal Hessians is
  contractible, so the normal determinant factor can
  be oriented consistently. The Morse--Bott comparison
  therefore gives an isomorphism
  \[
   SH^*_{\mathrm{uw}}(\gamma_m;\bZ)
   \simeq H^*(L_{b_m};\bZ).
  \]
  The same argument applies to every multiple cover,
  since the normal shear is rescaled by a positive
  factor and we do not quotient by rotation of the
  orbits.
  
  For a component over a flat edge, the family of
  parametrized orbits is the torus bundle over the
  flat segment. Slightly bowing the edge gives an
  isolating isotopy to a component supported over a
  single point of the same regular affine chart.
  Local Floer transport and the preceding calculation
  give the same identification with fiber cohomology.
  
  The desired formula now follows from
  Theorem~\ref{thm:background-indicator-input}.
  \end{proof}

We now consider restriction for compatible smoothings of a facewise
admissible inclusion. The Reeb components over a shared rounded corner are strictly
common. Along a partially common edge, however, the two Reeb
components need not coincide. We show that they are weakly common
in the sense of Section~\ref{subsec:background-local-indicator}.
Here and below $\Lambda_{[r,s)}$ denotes the ordinary log-norm window
of Section~\ref{subsec:background-filtrations}.

\begin{lemma}[Restriction along a partially common edge]
  \label{lmEdgeRestrictionIsotopy}
  Let $Q\subset P$ be a facewise admissible inclusion with
  $Q$ a regular convex polygon, and let $Q_t\subset P_t$
  be compatible smoothings. Let $e$ be an edge of $P$
  whose intersection with $Q$ is a nondegenerate segment.
  
  For each positive multiple of the primitive label of $e$,
  in a sufficiently short Novikov window the restriction
  between the corresponding Reeb blocks is given by local
  Floer transport, with zero action shift in the common
  normalization.
  \end{lemma}

\begin{proof}
Fix $t$ and write $D=P_t$ and $D'=Q_t$.
The edge and its endpoints lie in a single affine chart. Indeed,
the developing map on a thin simply connected neighborhood is
injective along the compact edge and hence on a smaller neighborhood.
In this chart write the common
supporting line as $v=h$, with $n=dv$ the primitive conormal.
Choose $c$ in the interior of $D'$ just below the relative
interior of its flat face. We may use the primitive
$\theta_c=\sum_j(p_j-c_j)dq_j$ on the target, since $D'$ is
convex and contains $c$ in its interior. The affine coordinates
extend along the source edge, so this primitive also defines a
radial coordinate near the source family. For a different target
primitive the identity inclusion of $K_{D'}$ identifies the two
presentations: every component is stationary and strictly common
\cite[Theorem~3]{Paper1}.

Cut both domains by a line parallel to and slightly below the
common supporting line, in the developed edge neighborhood.
This gives nested convex caps $B\subset A$ containing $c$ in
their interiors. On rays from $c$ meeting the source flat face
or its endpoint germs, the two outer fronts are the respective
boundaries of $D$ and $D'$. Indeed these rays point towards the supporting
line, so neither front meets the lower cut. Shrink to an open
sector of such rays containing both families. On this sector
the gauges $g_A,g_B$ about $c$ are smooth and $g_A\le g_B$.
The unit levels of
\[
 g_\tau=(1-\tau)g_A+\tau g_B,\qquad 0\le\tau\le1,
\]
are smooth radial graphs: their logarithmic heights are
$-\log g_\tau$ on rays of unit length. They decrease from the
source front to the target front.

The supporting level for $n$ is $v=h$ throughout. At $\tau=0$
its exposed face is the source face. For $\tau>0$ it is the
target face: both gauges are one on that face, whereas on the
rest of the source face $g_B>1$. For the label $mn$, the orbits
have fiber class $-mn$ and their $\theta_c$-period is
$m(h-v(c))$ at every $\tau$. This family is compact, including
at $\tau=0$. It is isolated in its fiber class: on either side
of the exposed face the conormal turns away from $n$, and
convexity excludes another face with this conormal. An isolating
neighborhood can be chosen in the sector, with its two ends
beyond the source face. This gives the radial-front continuation
of \cite[Definitions~4.5 and~5.13]{Paper1}.

The components are therefore weakly common (or strictly common
if their boundary germs agree). Theorem~\ref{thm:background-indicator-input}
now gives the restriction. The action difference in
the common local primitive is zero. To identify the unweighted
transport, bow the fronts slightly with their supporting point
in the interior of the common segment. For example, in affine
coordinates $(u,v)$ replace $v$ by $v+\varepsilon(u-u_0)^2$
in the equations of the fronts, with $u_0$ in that segment.
Near the exposed faces the upper graphs remain concave and nested.
For $\varepsilon>0$ they have the same germ at the supporting torus
over $u_0$. The families stay isolated as $\varepsilon$ tends to
zero. Transport in this two-parameter family intertwines the flat
and bowed comparisons by the two endpoint bowing isomorphisms of
Lemma~\ref{lmMxUwIso} \cite[Proposition~4.7]{Paper1}. In the bowed
family the torus and its framing are stationary, so the induced map
on the cohomology of the supporting fiber is the identity.
\end{proof}

\begin{corollary}[Restriction to rounded vertex caps]
  \label{tmInfShPtrop}
  Fix compatible roundings $D_v\subset D$ of the vertex
  caps and the polygon. For every finite period cutoff $a$
  there is a $\delta>0$ such that for every ordinary
  Novikov window $I$ of width less than $\delta$ the following
  holds.
  
  Under the decomposition of Lemma~\ref{lmMxUwIso},
  restriction
  \[
   SH^{\leq a}_{M,I}(K_D)
   \longrightarrow SH^*_{M,I}(K_{D_v})
  \]
  is given on the critical block by ordinary cohomological
  restriction. On the nonzero blocks whose labels belong
  to the vertex cone at $v$, it is given by local Floer
  transport to the corresponding blocks of the cap.
  It vanishes on all remaining blocks.
  
  In particular, in every degree $j$ in which $(T)$ holds,
  restriction induces an injection
  \[
   SH^{j,\leq a}_{M,I}(K_D)
   \hookrightarrow
   \bigoplus_{v\in V(P)} SH^j_{M,I}(K_{D_v}).
  \]
  \end{corollary}
  
  \begin{proof}
  Fix a vertex $v$. The Reeb components corresponding
  to labels in the interior of the vertex cone at $v$
  are strictly common. Those corresponding to nonzero
  labels on its two rays are common by
  Lemma~\ref{lmEdgeRestrictionIsotopy}.
  All remaining Reeb components are non-common.
  Choose a target period cutoff containing the components
  corresponding to the common source components of period
  at most $a$. The restriction formula now follows from
  Theorem~\ref{thm:background-indicator-input}.
  Since there are finitely many caps, we can choose the
  window width $\delta$ to work for all of them simultaneously.
  
  Every nonzero label belongs to at least one vertex cone.
  Moreover, distinct labels in a given vertex cone correspond
  to distinct blocks of the corresponding cap. Thus restriction
  to all the caps together is injective on the sum of the
  Reeb blocks.
  
  Finally, $K_D$ retracts onto $K_P$ through the regular
  collar, and each $K_{D_v}$ retracts onto a regular fiber.
  By $(T)$, restriction to the caps is also injective on
  the critical block in degree $j$.
  \end{proof}
  
\subsection{Proof of Theorem~\ref{thm:strict-detection}}
\label{subsec:cap-window-enlargement}
\label{subsec:cap-polygon-limit}

By Corollary~\ref{tmInfShPtrop}, restriction to the vertex
caps of a fixed rounding is injective up to any finite
period cutoff in sufficiently short Novikov windows.
The difficulty is that the allowed window width may tend
to zero as the period cutoff increases. Using
torsion-freeness of the ordinary open lattices of the caps,
we first enlarge the window at a fixed period cutoff and
then remove the cutoff. This allows us to pass to the polygon
itself in any finite window. We then deduce injectivity on
the ordinary open lattice and, finally, show that restriction
is an isometric embedding with closed image.

\begin{proof}
Choose compatible nested outward roundings
$Q_{v,\rho}\subset P_\rho$ as in
Section~\ref{sec:polygon-windows}, using the same rounding
$P_\rho$ for all the caps. Write
$K_\rho=\pi^{-1}(P_\rho)$ and
$K_{v,\rho}=\pi^{-1}(Q_{v,\rho})$.
Throughout the proof fix a degree $j$ in which $(T)$ holds.

\medskip
\noindent\emph{Step 1. Enlarging the window at a fixed period cutoff.}
Fix $\rho$ and write
\[
 C=SC_M^*(K_\rho;\Lambda),\qquad
 D=\bigoplus_{v\in V(P)}SC_M^*(K_{v,\rho};\Lambda).
\]
For either complex $B$ write
\[
 B_t=\{x\in B\mid \|x\|<e^t\},\qquad
 B_{[r,s)}=B_s/B_r.
\]

By the regular comparison of Section~\ref{sec:regular},
$H^*(D_0)$ is torsion-free over $\Lambda_{\geq0}$.
Multiplication by Novikov monomials identifies the inclusion
$D_r\hookrightarrow D_s$, for $r<s$, with multiplication
by $T^{s-r}$ on $D_0$. Thus inclusion is injective on
cohomology. By the long exact sequence of the window,
\[
 H^*(D_{[r,s)})
 \simeq H^*(D_s)/H^*(D_r).
\]
In particular, for $r<s<u$ the natural map
\[
 H^*(D_{[r,s)})\longrightarrow H^*(D_{[r,u)})
\]
is injective.

Fix a contact primitive near the boundary of $K_\rho$.
By \cite[Theorem~2 and its proof]{Paper1}, for some
$\hbar>0$ the complex $C_{[0,\hbar)}$ carries an
exhaustive action filtration by subcomplexes.
In the telescope presentation used there, this filtration
is obtained by taking the inverse image of the action
filtration on the quotient by the concave generators.
Since the filtration selects generators and passage to a
smaller Novikov window truncates their coefficients,
we obtain the corresponding filtration on every subwindow.
By Novikov translation we obtain such a filtration
$C_I^{\leq a}$ on every window $I$ of width less than
$\hbar$.

For each finite cutoff $a$ outside the period spectrum,
\cite[Theorem~2]{Paper1} identifies the cohomology
of $C_I^{\leq a}$, by inclusion, with the period-filtered
summand $SH^{*,\leq a}_{M,I}(K_\rho)$ in sufficiently
short windows. By Corollary~\ref{tmInfShPtrop}, there is
therefore a $\delta(a)>0$ such that
\[
 H^j(C_I^{\leq a})\longrightarrow H^j(D_I)
\]
is injective whenever $|I|<\delta(a)$.

We show that for fixed $a$ this map is in fact injective
for every window of width less than $\hbar$.
Suppose that injectivity holds on the windows $[r,s)$
and $[s,u)$, where $u-r<\hbar$.
The induced filtrations give a short exact sequence
\[
 0\longrightarrow C_{[r,s)}^{\leq a}
 \longrightarrow C_{[r,u)}^{\leq a}
 \longrightarrow C_{[s,u)}^{\leq a}
 \longrightarrow0.
\]
Let $x\in H^j(C_{[r,u)}^{\leq a})$ map to zero in
$H^j(D_{[r,u)})$. By injectivity on the window $[s,u)$,
the image of $x$ in $H^j(C_{[s,u)}^{\leq a})$ vanishes.
By the long exact sequence, $x$ comes from a class
$y\in H^j(C_{[r,s)}^{\leq a})$.
The image of $y$ in $H^j(D_{[r,s)})$ vanishes after
enlargement to $[r,u)$. By the injectivity established
above, it already vanishes in $H^j(D_{[r,s)})$.
Injectivity on the first window therefore gives $y=0$,
and hence $x=0$.

Subdividing any window of width less than $\hbar$ into
windows of width less than $\delta(a)$ and applying
this argument inductively proves the assertion.

\medskip
\noindent\emph{Step 2. Removing the period cutoff.}
By exhaustiveness of the filtration, for every window
$I$ of width less than $\hbar$,
\[
 C_I=\varinjlim_a C_I^{\leq a}.
\]
Since filtered colimits of modules are exact,
\[
 H^j(C_I)=\varinjlim_a H^j(C_I^{\leq a}).
\]
By Step 1 each of the groups on the right maps
injectively to $H^j(D_I)$. Hence so does $H^j(C_I)$.

Applying the subdivision argument of Step 1 to the full
window complexes, without the period cutoff, we now
obtain injectivity in arbitrary finite windows. Thus
\begin{equation}
 SH^j_{M,I}(K_\rho)
 \lhook\joinrel\longrightarrow
 \bigoplus_{v\in V(P)} SH^j_{M,I}(K_{v,\rho})
 \label{eq:tors-smooth-injection}
\end{equation}
for every finite window $I$.

\medskip
\noindent\emph{Step 3. Passing to the polygon.}
Fix a finite window $I$. The strictly admissible
Hamiltonians for the roundings form a cofinal family
for $K_P$, since a Hamiltonian strictly negative on
$K_P$ remains negative on $K_\rho$ for sufficiently
small $\rho$. The same holds for the caps.

In a finite window the completed telescope agrees
with the algebraic telescope, since each element has
only finitely many coordinates surviving the lower
log-norm cutoff \cite[Lemma~2.3]{Paper1}. Consequently,
\[
 SH^j_{M,I}(K_P)
 =\varinjlim_{\rho\downarrow0}SH^j_{M,I}(K_\rho),
\]
and similarly for the caps. By functoriality,
the maps \eqref{eq:tors-smooth-injection} commute
with the maps in these directed systems. Taking the
direct limit, we obtain
\begin{equation}
 SH^j_{M,I}(K_P)
 \lhook\joinrel\longrightarrow
 \bigoplus_{v\in V(P)}SH^j_{M,I}(K_{Q_v})
 \label{eq:tors-polygon-injection}
\end{equation}
for every finite window $I$.

\medskip
\noindent\emph{Step 4. Passing to the ordinary open lattice.}
Now write
\[
 C=SC_M^*(K_P;\Lambda),\qquad
 D=\bigoplus_{v\in V(P)}SC_M^*(K_{Q_v};\Lambda),
\]
and denote the restriction map by $\phi:C\to D$.

Let $z_0\in C_0^j$ be a cocycle whose image is exact
in $D_0$. We inductively construct chains
$y_k\in C_{-k}^{j-1}$ such that
\[
 z_{k+1}=z_k-dy_k\in C_{-k-1}^j.
\]
Suppose that $z_k\in C_{-k}^j$ has been constructed.
Its image is exact in $D_0$. Since the ordinary open
lattices of the caps are torsion-free by
Lemma~\ref{lem:ordinary-open-lattice}, the argument at
the beginning of Step 1 shows that
$H^j(D_{-k})\to H^j(D_0)$ is injective.
Thus $\phi(z_k)$ is already exact in $D_{-k}$.
By \eqref{eq:tors-polygon-injection}, the reduction of
$z_k$ in the window $[-k-1,-k)$ is exact, giving the
required chain $y_k$.

Since $\|y_k\|<e^{-k}$, completeness implies that the
series
\[
 y=\sum_{k=0}^\infty y_k
\]
converges in $C_0$. As $z_k\to0$, we have $dy=z_0$.
This proves injectivity of
\[
 \cF^j_0(P)\longrightarrow
 \bigoplus_{v\in V(P)}\cF^j_0(Q_v).
\]
The target is torsion-free, so the source is
torsion-free as well.

\medskip
\noindent\emph{Step 5. Norm detection and closedness.}
By Novikov translation, the injectivity just proved
holds at every filtration level. Since every cocycle
and every primitive in the completed complex belongs
to some filtration level, restriction is also
injective on full cohomology.

We show that restriction detects the ordinary
filtration on cohomology. Let $r<s$ and suppose that
$x\in H^j(C_s)$ has image in
$\operatorname{im}(H^j(D_r)\to H^j(D_s))$.
Then its image in $H^j(D_{[r,s)})$ vanishes.
By \eqref{eq:tors-polygon-injection}, its reduction in
$H^j(C_{[r,s)})$ vanishes as well. The long exact
sequence of the window therefore shows that $x$
comes from $H^j(C_r)$.

Recall that the ordinary filtration on $\cF^j(P)$
is given by
\[
 F_t\cF^j(P)=
 \operatorname{im}(H^j(C_t)\to H^j(C)).
\]
By torsion-freeness of the target lattice, each
$H^j(D_t)$ injects into $H^j(D)$. Thus, applying the
preceding argument with $s$ large enough to represent
a given class, we obtain
\[
 x\in F_t\cF^j(P)
 \quad\Longleftrightarrow\quad
 R_P^j(x)\in\bigoplus_{v\in V(P)}F_t\cF^j(Q_v).
\]
Consequently,
\[
 \|x\|_P=\max_{v\in V(P)}
 \|x|_{Q_v}\|_{Q_v}.
\]

In fact, the caps can be shrunk arbitrarily towards
their vertices. By the regular reconstruction,
the analytic polyvector norms vary continuously
with the base point. Since there are finitely many
vertices, letting the caps shrink gives
\[
 \|x\|_P=\max_{v\in V(P)}\|x|_v\|_v.
 \tag{VN}
\]
Here $\|x|_v\|_v$ denotes the limit of
$\|x|_{Q_v}\|_{Q_v}$ as the vertex cap $Q_v$ shrinks
to $v$.
In degree zero these are the multiplicative Gauss
norms, and in higher degrees they are the regular
polyvector norms in the normalized torus frame.

It remains to establish completeness. The kernel of
the Floer differential is closed. Suppose that a
sequence of boundaries converges to a cocycle $z$
in $C$.
By Theorem~\ref{thm:local-regular-reconstruction},
the cohomology of each regular cap is Hausdorff in
the quotient topology. Thus its boundaries are closed,
and the restriction of $z$ to each cap is a boundary.
By injectivity of restriction on
cohomology, $z$ is itself a boundary. Thus the image
of the differential is closed in the source as well,
and the quotient norm on $\cF^j(P)$ is complete.
An isometric image of a complete normed space is
closed, completing the proof of the theorem.
\end{proof}

\begin{corollary}[Concentration]\label{cor:concentration}
Let $j\notin\{0,1,2\}$. If ordinary restriction to a boundary
fiber is injective in degree $j$, then $\cF^j(P)=0$.
\end{corollary}
\begin{proof}
The regular cap groups vanish in degree $j$ by
Theorem~\ref{thm:local-regular-reconstruction}. The joint restriction
to these groups is injective by Theorem~\ref{thm:strict-detection}.
\end{proof}

\subsection{Examples satisfying the hypotheses}
\label{subsec:reconstruction-examples}

\paragraph{Almost toric fibrations.}
Consider a singular fiber with $n\geq1$ focus--focus
points in an almost toric fibration. Its base has the
standard affine slit model, with the two banks of the
slit identified by the shear
$(x,y)\mapsto(x+ny,y)$
\cite[Section~3]{LeungSymington}.
We take the slit along the positive $x$-axis, which is fixed
by this shear.

Dilations about the singular value descend to the
affine base. Consider a sufficiently small dilate of
the diamond with vertices
\[
 (\pm a,0),\qquad (0,\pm1),
\]
where $a$ is a positive integer with $2a>n$.
At the vertex on the slit the determinant of the
incoming and outgoing tangent directions, computed
after identifying the banks, is either $2a-n$ or
$2a+n$, according to the slit convention.
It is positive. Thus the polygon is strictly locally
convex, including at the vertex on the slit.

We can make this polygon Delzant by truncating its
corners. Indeed, for two consecutive primitive outward
conormals $u,v$ with $|\det(u,v)|>1$, choose a primitive
lattice vector $w$ in the interior of their cone such
that
\[
 0<|\det(u,w)|,\ |\det(w,v)|<|\det(u,v)|.
\]
Such a vector is obtained from a nonzero lattice point
in the fundamental parallelogram of $u,v$.
A sufficiently shallow truncation introduces an edge
with outward conormal $w$, replacing the original
corner by two corners of smaller determinants.
Induction gives a Delzant polygon.

These truncations may be chosen so that the Euler
field remains strictly outward-pointing on every edge.
Indeed, the new outward conormal is a positive linear
combination of the two original conormals, and the
supporting constant can be chosen arbitrarily close
to its value at the original corner.
The resulting polygon has the required cyclic sequence
of labeled edges. Its inward dilates give a neighborhood
basis of the singular value. This verifies
$(G_\partial)$ and $(G_{\rm rad})$.

For the topological conditions, a sufficiently small
fibered neighborhood is obtained, up to homotopy,
from a regular torus by attaching $n$ Lefschetz
thimbles along parallel primitive vanishing cycles.
In the homology sequence of this pair the boundary
homomorphism
\[
 H_2(K_P,L;\kk)\simeq\kk^n
       \longrightarrow H_1(L;\kk)
\]
has rank one. In particular,
$H_1(L;\kk)\to H_1(K_P;\kk)$ is surjective and
\[
 b_2(K_P;\kk)=1+n-1=n.
\]
The space $K_P$ is connected, so $(T)_0$ and $(T)_1$
hold. The map $H^2(K_P;\kk)\to H^2(L;\kk)$ is
surjective and is injective precisely when $n=1$.

\paragraph{Elliptic surfaces.}
Let a hyperk\"ahler surface have an elliptic fibration
in complex structure $I$. Write
$\omega_I,\omega_J,\omega_K$ for its K\"ahler forms.
The regular elliptic fibers are Lagrangian with respect
to $\omega_J$, since the holomorphic symplectic form
$\omega_J+i\omega_K$ restricts to zero on them.
In complex structure $J$, the form
$\omega_K+i\omega_I$ has constant phase on the oriented
fibers. Thus the fibration-induced grading trivialization
extends over the whole surface.

Consider a nonmultiple Kodaira fiber of a relatively
minimal elliptic surface with smooth total space.
Write its divisor as
\[
 F=\sum_{\nu=1}^c m_\nu C_\nu.
\]
A sufficiently small fibered neighborhood retracts onto
the reduced fiber. One way to see this is to use the
gradient flow of $|t\circ\pi|^2$, where $t$ is a
holomorphic coordinate centered at the singular value.
The {\L}ojasiewicz gradient inequality near the compact
zero set gives finite-length convergence to that set and a
deformation retraction for a sufficiently small
sublevel set; compare also the local-neighborhood discussion
in \cite[Section~3.3]{BuzzardLu}.

The inclusion of a nearby regular fiber induces a
surjection on fundamental groups
\cite[Section~3.3, Proposition~3.13]{BuzzardLu}.
In particular, $(T)_0$ and $(T)_1$ hold.

The component classes $[C_\nu]$ form a basis of
$H_2(K_P;\kk)$. Indeed, the reduced fiber is a union
of irreducible compact curves meeting in finitely
many points, and each irreducible component contributes
one fundamental class in degree two.
Specialization of the regular fiber gives
\[
 [L]=\sum_{\nu=1}^c m_\nu[C_\nu].
\]
Consequently, in the dual basis the degree-two
restriction is
\[
 \kk^c\longrightarrow\kk,\qquad
 (a_1,\ldots,a_c)\longmapsto
 \sum_{\nu=1}^c m_\nu a_\nu.
\]
Since $\kk$ has characteristic zero, this map is
surjective. Among the nonmultiple singular Kodaira
fibers, $(T)_2$ holds precisely for the irreducible
types $I_1$ and $II$
\cite[Section~4]{SchuettShioda}.

We next verify the geometric conditions in several
examples with finite affine monodromy. For the local period
asymptotics at such fibers on elliptic K3 surfaces, see
\cite[\S4.2]{LinTropical}.

Suppose that a local developing map of the punctured
base, in a complex coordinate $u$ centered at the
singular value, is given by
\[
 u\longmapsto u^\beta
\]
in the real period plane, where the period lattice
is preserved by the monodromy.
The affine base is obtained from a sector of angle
$2\pi\beta$ by identifying its banks by the
corresponding linear monodromy. In particular, the
Euler field descends.

For the exponents
\[
 \beta=\frac56,\quad\frac34,\quad\frac23,\quad\frac12,
\]
one can choose respectively a monodromy-invariant
Delzant hexagon, square, hexagon or centrally
symmetric Delzant hexagon in the period plane,
with the origin in its interior.

To see this in the cases of order six and three
monodromy, choose a lattice basis in which the
primitive outward conormals, in cyclic order, are
\[
 e_1,\ e_1+e_2,\ e_2,\ -e_1,\ -e_1-e_2,\ -e_2.
\]
Consecutive conormals form integral bases.
The monodromy permutes the six rays, and the
supporting constants may all be taken equal.
For order four monodromy, use the square with
conormals $\pm e_1,\pm e_2$. For monodromy $-1$,
the centrally symmetric hexagon above works.

Choose the slit along a ray through a vertex of
the chosen polygon. Taking the five, three, four
or three consecutive sides lying in the respective
developed sectors gives a polygon whose boundary
closes up after identifying the banks.
Local convexity and the Delzant condition at the
identified vertex follow from monodromy invariance.
The Euler field is outward-pointing, and its
inward flow gives the required neighborhood basis.
Thus $(G_\partial)$ and $(G_{\rm rad})$ hold in
all four cases.

An isotrivial elliptic K3 surface is given by the
Weierstrass equation
\[
 y^2=x^3+b(t),
\]
where $b$ is a section of $\cO_{\bP^1}(12)$ with
simple zeros. Its twelve singular fibers are
cuspidal, of type $II$
\cite[Lemma~7]{SawonIsotrivial}.

Near a zero of $b$, the changes of variables
$x=b(t)^{1/3}X$ and $y=b(t)^{1/2}Y$ identify the
smooth elliptic fibers with a fixed elliptic curve.
Their periods of $dx/y$ are therefore constant
multiples of $b(t)^{-1/6}$.
For the hyperk\"ahler-rotated symplectic form, the
action coordinates are the real parts of the
periods of a primitive of the holomorphic
symplectic form. They are thus obtained from
\[
 \int_0^t b(v)^{-1/6}\,dv.
\]
As the zero of $b$ is simple, this integral is
$t^{5/6}$ times a holomorphic unit.
After a local change of base coordinate it becomes
$u^{5/6}$. The preceding sector construction applies.
Since the cuspidal fiber is irreducible, all the
topological conditions hold as well.

For another compact example, take the Kummer K3
obtained by resolving $(E\times E')/\{\pm1\}$,
with its projection to $E'/\{\pm1\}$
\cite[Section~2]{SawonIsotrivial}.
In a coordinate $w$ on $E'$ centered at a branch
point, the holomorphic symplectic form upstairs is
a constant multiple of $dw\wedge dz$, where $dz$
is the invariant one-form on $E$. The action
coordinates are consequently real parts of
complex linear functions of $w$.
Downstairs, with base coordinate $u=w^2$, the
developing map is of the form $u^{1/2}$.

The centrally symmetric hexagon construction
therefore gives the required polygons.
Each singular fiber is of type $I_0^*$, with five
irreducible components. Thus $(T)_0$ and $(T)_1$
hold, whereas the kernel of ordinary degree-two
restriction has dimension four.

\paragraph{Noncompact elliptic and Hitchin fibrations.}
Consider the relatively minimal elliptic surfaces
obtained from the Weierstrass equations
\[
 y^2=x^3+t,\qquad
 y^2=x^3+tx,\qquad
 y^2=x^3+t^2
\]
by compactifying the base to $\bP^1$ and taking
the minimal regular models. Their fibers at zero
are of types $II$, $III$ and $IV$, respectively,
and their fibers at infinity are of types
$II^*$, $III^*$ and $IV^*$.
Remove the fiber at infinity.

These open surfaces carry complete hyperk\"ahler
metrics with ALG asymptotics
\cite{HeinInstantons}. Cheung and Lin study the family Floer mirrors
of these three surfaces in \cite{CheungLin}.

Rescaling the Weierstrass coordinates identifies
the smooth fibers in each example with a fixed
elliptic curve. The periods of $dx/y$ are constant
multiples of
\[
 t^{-1/6},\qquad t^{-1/4},\qquad t^{-1/3},
\]
respectively. Integrating gives developing maps
of the form
\[
 t^{5/6},\qquad t^{3/4},\qquad t^{2/3}.
\]
These are the period exponents computed in
\cite[Lemmas~5.6, 6.2 and~7.2]{CheungLin}.
The sector constructions above verify the geometric
conditions. In all three examples $(T)_0$ and
$(T)_1$ hold, and $(T)_2$ holds for the type $II$
fiber.

For the four-punctured sphere, consider the strongly parabolic
$\mathrm{SL}(2,\bC)$ Hitchin system with generic weights studied
in \cite[Section~8]{FMSWParabolic}.
Its Hitchin base has coordinate $q$ and its only
singular fiber is the nilpotent fiber.
The regular spectral curves become isomorphic
to a fixed elliptic curve after the base change
$q=w^2$. Under these identifications the spectral
differential is multiplied by $w$.
The developing map is therefore of the form
$q^{1/2}$.

In Hausel's model \cite[\S4.2]{HauselToy}, the moduli space
is the resolution of $(\bC_w\times E)/\{\pm1\}$ and the
central fiber is of type $I_0^*$.
The same sector construction as for the Kummer
example verifies $(G_\partial)$ and $(G_{\rm rad})$.
The conditions $(T)_0$ and $(T)_1$ hold, but
$(T)_2$ fails.

We finally verify the geometric finite-type condition
for these noncompact examples.
The ALG asymptotics of the metrics
\cite{HeinInstantons,FMSWParabolic}
give an end modeled on a flat torus fibration
over a Euclidean cone, with asymptotic control
of the metric and its derivatives.
In particular, the metrics have bounded curvature
and a positive lower bound for the injectivity radius.

Let $r$ be the radial coordinate of the model end.
It extends, after modification on a compact set,
to a proper bounded-below smooth function $f$.
The critical set of $f$ is compact. Moreover,
the metric and its first derivatives are uniformly
controlled relative to the model, in which
$|dr|=1$ and the Hessian of $r$ is bounded.
It follows that $|\nabla f|$ and
$|\nabla^2 f|$ are uniformly bounded.

For the rotated hyperk\"ahler structure the
compatible complex structure is parallel.
Consequently, the Hamiltonian vector field $X_f$
and its covariant derivative are uniformly bounded.
This is the geometric finite-type condition of
Section~\ref{sec:floer-background}.

In Section~\ref{sec:affinoid} we use torsion-freeness to lift the
primitive edge classes from a finite window to $\cF^0(P)$. We can
then compute the norms of their products on the vertex caps.

\section{The Floer algebra of an admissible polygon}\label{sec:affinoid}

Let $P$ be an admissible polygon with edges $e_1,\ldots,e_N$ in cyclic
order. The regular mirror $Y_{\rm reg}$ of
Section~\ref{sec:regular} is defined over $\partial P$. To extend
the construction over $P$, we first determine the algebra
$A_P=\cF^0(P)$. For an edge or vertex $F$, let
$Y_F=\operatorname{Sp}\cF^0(F)$ be its affinoid domain in the
regular mirror, and put $U_\partial=\bigcup_iY_{e_i}$.
Their Berkovich realizations are $p^{-1}(F)$ and
$p^{-1}(\partial P)$, respectively. Restriction from $P$ to its edges gives
functions on the $Y_{e_i}$ that agree at the common vertices. So we obtain
a map
\[
 r_\partial:A_P\longrightarrow B_\partial=\Gamma(U_\partial,\cO)
\]
of algebras. We equip $B_\partial$ with the maximum of the Gauss
norms on the edge domains.

Let $\xi_i$ be the primitive outward conormal to $e_i$. In a regular
affine chart along $e_i$ take a small regular domain $Q$ having $e_i$ as a
supporting face, and use a tropically marked Laurent presentation of
$\cF^*(Q)$ from Theorem~\ref{thm:local-regular-reconstruction}.
The geometric edge label is $\xi_i$. Its Floer orbits have winding
$-\xi_i$ by Lemma~\ref{lmMxUwIso}, so the marking sends the
degree-zero unit of this block to the reduction of the normalized
character of exponent $-\xi_i$.
A \emph{normalized edge lift} for $e_i$ is a norm-one class in $A_P$
with this leading term along $e_i$.

Recall from \eqref{eqPtropUnion} that $P_{trop}(\bZ)$, which labels the blocks of local Floer cohomology of $\pi^{-1}(P)$, is the union of
the vertex cones $C_{i,i+1}=\bN\xi_i+\bN\xi_{i+1}$, with the two
copies of each edge ray identified. Given edge lifts $z_i$, put
\[
 b_0=1,\qquad b_\gamma=z_i^a z_{i+1}^b
 \quad\text{for }\gamma=a\xi_i+b\xi_{i+1}\in C_{i,i+1}.
\]
On a ray common to two cones both formulas give $b_{k\xi_i}=z_i^k$.
For a countable set $I$, write $c_0(I,\Lambda)$ for the Banach space whose elements are 
families $(a_i)_{i\in I}$ such that for each $\epsilon>0$ only
finitely many $i$ satisfy $|a_i|\ge\epsilon$. We give it the
supremum norm. Put $E=c_0(P_{trop}(\bZ),\Lambda)$.

\begingroup
\renewcommand{\thetheorem}{\thesection}
\begin{theorem}[The Floer algebra of an admissible polygon]\label{tmAffinoidPolygon}
Assume that the admissible polygon $P$ has at least three labeled
edges and distinct adjacent labeled pairs, as in $(G_\partial)$,
and that $(T)$ holds in degrees 0 and 1. Use the tropical marking of
Section~\ref{sec:regular} on its regular collar.
\begin{enumerate}[label=\textup{(\roman*)},ref=(\roman*)]
\item\label{part:polygon-basis}
There are normalized edge lifts $z_1,\ldots,z_N$ for which
\[
 S:E\xrightarrow{\ \sim\ }A_P,\qquad
 (c_\gamma)\longmapsto\sum_\gamma c_\gamma b_\gamma
\]
is an isometric isomorphism of Banach spaces. Thus each class has a
unique such expansion, of norm $\max_\gamma|c_\gamma|$.
\item\label{part:polygon-affinoid}
The homomorphism $\Lambda\langle Z_1,\ldots,Z_N\rangle\to A_P$,
$Z_i\mapsto z_i$, is surjective. Its quotient norm is the spectral
norm and equals $\|\cdot\|_P$. In particular, $A_P$ is a strict
affinoid algebra, and $\dim A_P=2$.
\item\label{part:polygon-boundary-functions}
Floer restriction $r_\partial:A_P\to B_\partial$ is an isometric
isomorphism of $\Lambda$-algebras.
\item\label{part:polygon-extension}
Let $U$ be a connected smooth region of the regular mirror containing
$U_\partial$ as an admissible domain, and write $j:U_\partial\hookrightarrow U$
for the inclusion.
If $U$ receives a Floer restriction $r_U:A_P\to\Gamma(U,\cO)$
satisfying $j^*r_U=r_\partial$, then both arrows
\[
 A_P\xrightarrow{\,r_U\,}\Gamma(U,\cO)\xrightarrow{\,j^*\,}B_\partial
\]
are isomorphisms of $\Lambda$-algebras.
\item\label{part:polygon-normal}
$A_P$ is an integrally closed domain. Thus $X_P=\operatorname{Sp}A_P$
is a normal integral affinoid surface.
\end{enumerate}
\end{theorem}
\endgroup

For the nodal models the boundary extension in
\ref{part:polygon-boundary-functions} was proved in
\cite[Propositions~3.19 and~9.1]{GromanVarolgunes2022}.
Here we prove the theorem by lifting the primitive edge classes and
reading Laurent coefficients at the vertices.
\subsection{Lifting the edge classes}

Choose a small vertex cap at each vertex. Cover the rest of the
boundary by finitely many small polygons contained in $P\cap B_{reg}$,
each meeting just one edge in an interval. These polygons and the
vertex caps cover a relative neighborhood of $\partial P$. We call
them the \emph{testing polygons}. For each testing polygon $Q$ fix a marked presentation
$\Phi_Q:\cF^0(Q)\xrightarrow{\sim}A(Q)$. In conical coordinates the
supporting line of $e_i$ is $\langle\xi_i,b\rangle=h_i$, with $h_i>0$
by outward transversality. When $Q$ meets $e_i$, put
\[
 x_i(b)=h_i-\langle\xi_i,b\rangle,\qquad
 m_{i,Q}=\Phi_Q^{-1}(e_{Q,-\xi_i}),\qquad
 e_{Q,-\xi_i}=T^{h_i}z^{-\xi_i}.
\]
Thus $|m_{i,Q}|=e^{-x_i}$ in this presentation. The representative
$m_{i,Q}$ uses the chosen chart. Its leading symbol with coefficient
one is fixed by the tropical marking. When $Q$ is understood we write $m_i$.

\begin{lemma}[Normalized edge lifts]\label{lm:finite-tropical-lift}
There are $z_i\in A_P$ of norm one and a common $q<1$ such that on
every testing polygon $Q$,
\[
 \|z_i|_Q-m_{i,Q}\|_Q\le q\quad(Q\cap e_i\ne\varnothing),\qquad
 \|z_i|_Q\|_Q\le q\quad(Q\cap e_i=\varnothing).
\]
\end{lemma}

\begin{proof}
      Fix an edge $e_i$. We first construct the desired class in a short
      Novikov window. To this end, enlarge $P$ to a smooth admissible
      domain $D_i$ whose boundary meets $P$ precisely along $e_i$, with
      both endpoints of $e_i$ in the interior of a flat face of $D_i$.
      This is done as in Proposition~\ref{prpCornerSmoothing}, keeping
      the supporting line of $e_i$ fixed, displacing the remaining
      supporting lines outward, and rounding at a scale much smaller
      than the displacements.
      
      For each testing polygon $Q$, choose a smooth convex domain
      $Q\subset R\subset D_i$ as in
      Figure~\ref{fig:edge-lift-construction}. If $Q$ meets $e_i$,
      we take the common boundary of $R$ and $D_i$ to be a segment
      in the interior of the distinguished flat face of $D_i$.
      Otherwise we take $R$ in the interior of $D_i$. The testing
      polygons can be chosen small enough that $R$ and $Q$ lie in a
      single regular neighborhood at the scale of
      \cite[Theorem~5]{Paper1}.
      
      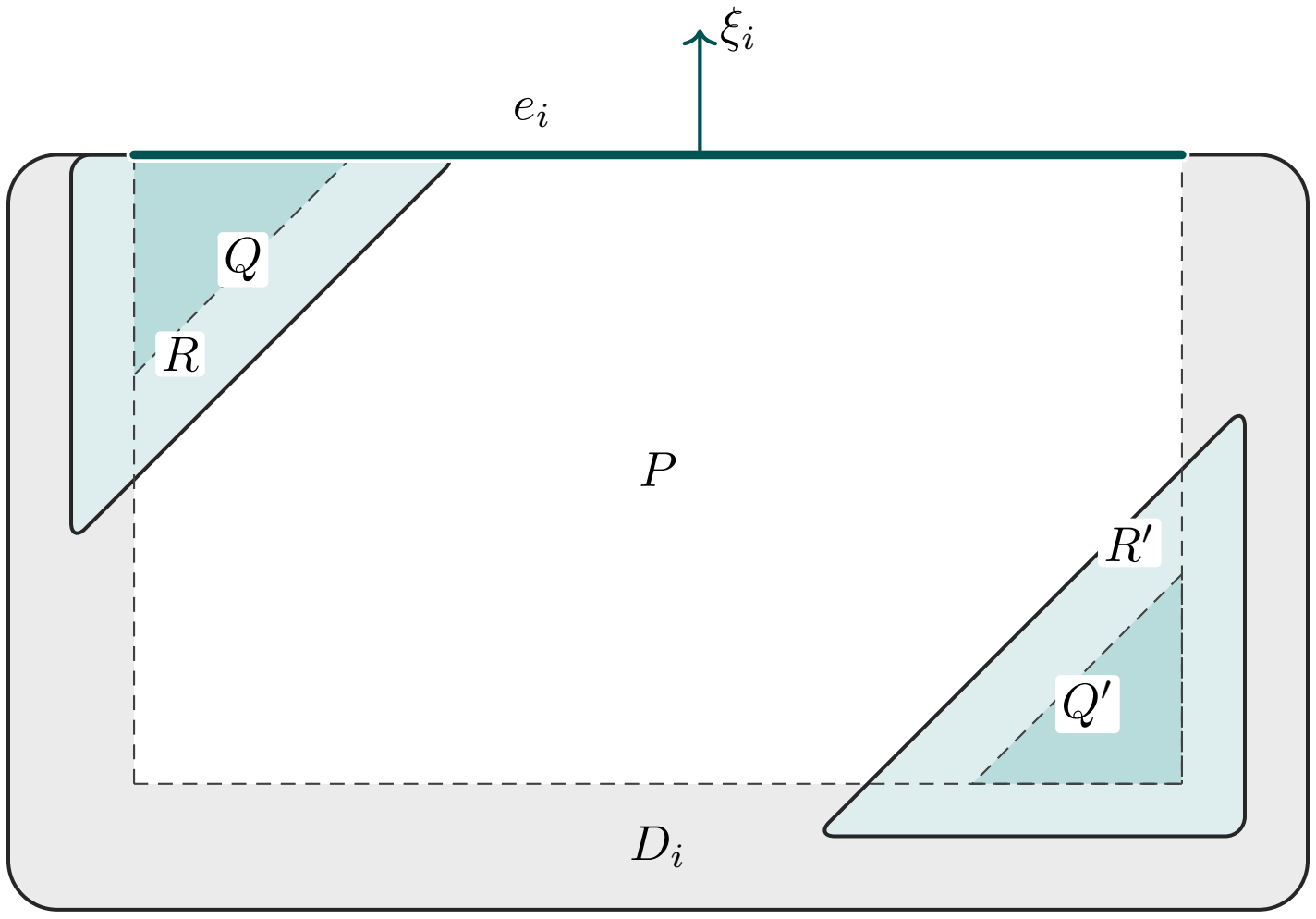
\begin{figure}[htbp]
      \centering
      \begin{tikzpicture}[
          x=1.5cm,y=1.5cm,font=\small,
          front/.style={black!85,line width=.75pt},
          polygon/.style={black!75,thin,densely dashed},
          common/.style={teal!65!black,line width=1.8pt,
                         line cap=round},
          lab/.style={fill=white,inner sep=1.2pt,rounded corners=1pt}
        ]
      
        \fill[black!8,rounded corners=10pt]
          (-.6,-.6) rectangle (5.6,3);
        \fill[white] (0,0) rectangle (5,3);
      
        \fill[teal!13,rounded corners=4pt]
          (-.3,3) -- (1.55,3) -- (-.3,1.15) -- cycle;
        \draw[front,rounded corners=4pt]
          (-.3,3) -- (1.55,3) -- (-.3,1.15) -- cycle;
      
        \fill[teal!13,rounded corners=4pt]
          (3.25,-.25) -- (5.3,-.25) -- (5.3,1.8) -- cycle;
        \draw[front,rounded corners=4pt]
          (3.25,-.25) -- (5.3,-.25) -- (5.3,1.8) -- cycle;
      
        \fill[teal!28]
          (0,3) -- (1.05,3) -- (0,1.95) -- cycle;
        \draw[polygon]
          (0,3) -- (1.05,3) -- (0,1.95) -- cycle;
      
        \fill[teal!28]
          (5,0) -- (4,0) -- (5,1) -- cycle;
        \draw[polygon]
          (5,0) -- (4,0) -- (5,1) -- cycle;
      
        \draw[front,rounded corners=10pt]
          (-.6,-.6) rectangle (5.6,3);
        \draw[polygon] (0,0) rectangle (5,3);
      
        \draw[white,line width=3.2pt,line cap=round] (0,3) -- (5,3);
        \draw[common] (0,3) -- (5,3);
        \node[above,lab] at (1.9,3.1) {$e_i$};
        \draw[->,teal!65!black,thick] (2.7,3) -- (2.7,3.6)
          node[black,right] {$\xi_i$};
      
        \node at (2.5,1.5) {$P$};
        \node at (2.5,-.3) {$D_i$};
      
        \node[lab] at (.52,2.5) {$Q$};
        \node[lab] at (.22,2.05) {$R$};
      
        \node[lab] at (4.55,.38) {$Q'$};
        \node[lab] at (4.75,1.15) {$R'$};
      
      \end{tikzpicture}
      \caption{The enlargement $D_i$ of $P$ keeps only the edge
      $e_i$ fixed. Two choices of testing polygon are shown. On the
      left, $Q\subset R$ meets $e_i$ and $R$ shares a flat segment
      with $D_i$. On the right, $Q'\subset R'$ is disjoint from
      $e_i$, with $R'$ contained in the interior of $D_i$.}
      \label{fig:edge-lift-construction}
      \end{figure}
      
      Consider the primitive Reeb block associated with the flat face
      of $D_i$ containing $e_i$. If $Q$ meets $e_i$, the argument of
      Lemma~\ref{lmEdgeRestrictionIsotopy} applies to the corresponding
      blocks of $K_{D_i}$ and $K_R$: the two faces lie on the same
      supporting line in a regular affine chart. Thus they are weakly
      common, with zero action shift, and local Floer transport identifies
      their degree-zero units. If $Q$ is disjoint from $e_i$, the source
      block is non-common.
      
      Choose a finite source period cutoff containing the primitive block
      and target cutoffs containing the corresponding common blocks.
      Take a log-norm window $I=[-\epsilon,\delta)$ sufficiently
      short for all these comparisons and for restriction from $R$ to $Q$.
      By Theorem~\ref{thm:background-indicator-input}, the degree-zero unit
      of the primitive block defines a class in $SH^0_{M,I}(K_{D_i})$ whose
      restriction to $K_R$ is the corresponding unit if $Q$ meets $e_i$,
      and zero otherwise.
      
      In the former case, restriction further to $Q$ gives the
      degree-zero unit in the unweighted block of winding $-\xi_i$,
      tensored with the class of $1$ in $\Lambda_I^{\log}$.
      Indeed, $R$ and $Q$ have the same
      supporting line, so their restriction weight in
      \cite[Theorem~5]{Paper1} is one.
      
      Let $u_i$ be the restriction of the constructed window class to
      $K_P$. By functoriality, its restriction to each $K_Q$ is the class
      just computed.
      
      It remains to lift $u_i$ to Floer cohomology. Write $C$ for the
      completed Floer complex of $K_P$. In the exact sequence
      \[
       H^0(C_\delta)
       \longrightarrow H^0(C_\delta/C_{-\epsilon})
       \longrightarrow H^1(C_{-\epsilon})
       \longrightarrow H^1(C_\delta),
      \]
      the last arrow is injective. Indeed, by Novikov translation it
      corresponds to multiplication by $T^{\epsilon+\delta}$ on
      $H^1(C_0)$, which is torsion-free by
      Theorem~\ref{thm:strict-detection}. Hence $u_i$ lifts to
      $H^0(C_\delta)$. Denote its image in $A_P$ by $z_i$.
      
      For a testing polygon $Q$ disjoint from $e_i$, exactness of the
      corresponding window sequence on $Q$ gives
      \[
       \|z_i|_Q\|_Q<e^{-\epsilon}.
      \]
      If $Q$ meets $e_i$, the coefficient $1$ comes from every
      subwindow $[-\epsilon,\eta)$ with $0<\eta<\delta$.
      On $Q$ the sequence of open lattices is injective on
      cohomology in every degree by Lemma~\ref{lem:ordinary-open-lattice}.
      Thus each of these subwindow classes lifts to a class of norm
      less than $e^\eta$. Its difference from $z_i|_Q$ has norm
      less than $e^{-\epsilon}$. Letting $\eta$ decrease to zero
      gives $\|z_i|_Q\|_Q\leq1$. The marking identifies its reduction
      with that of $m_{i,Q}$, whence
      \[
       \|z_i|_Q-m_{i,Q}\|_Q<1.
      \]
      Here the passage from valuation windows to the open log-norm window
      is given by the exact filtered colimit
      \[
       C_\delta/C_{-\epsilon}
       =
       \varinjlim_{\eta\downarrow0}
       G^{-\delta+\eta}C/G^{\epsilon+\eta}C.
      \]
      
      In particular, the restriction of $z_i$ to either incident vertex
      cap has norm one, while its restriction to every other vertex cap
      has norm less than one. By Theorem~\ref{thm:strict-detection},
      $\|z_i\|_P=1$. Since there are only finitely many edges and testing
      polygons, all the strict inequalities above are bounded by a
      common $q<1$.
      \end{proof}

\subsection{Laurent coefficients on the boundary}

Choose the lifts $z_i$ of Lemma~\ref{lm:finite-tropical-lift} and recall
that $b_{a\xi_i+b\xi_{i+1}}=z_i^az_{i+1}^b$ for $a,b\geq0$.
We wish to define a linear map $R:B_\partial\to E$ which reads the
coefficient of each of these monomials at a vertex where it is supported.
The indexing set for the coefficients has the decomposition
\[
 P_{trop}(\bZ)=\{0\}\sqcup
 \bigsqcup_i\bZ_{>0}\xi_i\sqcup
 \bigsqcup_i C_{i,i+1}^{\circ}(\bZ).
\]
For a label $\gamma$ in the interior of $C_{i,i+1}$ we take
$v_\gamma=e_i\cap e_{i+1}$. For $\gamma=k\xi_i$, $k>0$, we take
one of the two endpoints of $e_i$, choosing the same endpoint for all
$k$. For $\gamma=0$ we fix an arbitrary vertex $v_0$.

First note that the restrictions of $z_i,z_{i+1}$ to the algebra
$\Gamma(Y_v,\cO)$ at $v=e_i\cap e_{i+1}$ give Laurent coordinates.
Indeed, by the Delzant condition and
Theorem~\ref{thm:local-regular-reconstruction}, the normalized
characters $m_i,m_{i+1}$ of the two incident edges in the cap
presentation restrict to Laurent coordinates of norm one on $Y_v$.
In particular their integral powers and products form an orthonormal
basis. On $Y_v$, Lemma~\ref{lm:finite-tropical-lift} gives
$\|z_j-m_j\|_v\leq q<1$ for $j=i,i+1$, where we suppress restriction
notation. Thus $z_j/m_j$ is a unit differing from $1$ by an element
of norm at most $q$. The same is true of each of its integral powers,
including negative powers. It follows that
\[
 \|z_i^az_{i+1}^b-m_i^am_{i+1}^b\|_v\leq q
 \qquad(a,b\in\bZ).
\]
The linear operator taking the Laurent basis in the $m_j$ to the
corresponding monomials in the $z_j$ differs from the identity in
operator norm by at most $q$. Its inverse is given by the Neumann
series. Hence the $z_i^az_{i+1}^b$, $a,b\in\bZ$, also form an
orthonormal basis of $\Gamma(Y_v,\cO)$.

Given $f\in B_\partial$, expand its restriction at each vertex
$v=e_i\cap e_{i+1}$ as
\[
 f|_{Y_v}=\sum_{a,b\in\bZ}c^{(i)}_{a,b}z_i^az_{i+1}^b.
\]
For $\gamma\in P_{trop}(\bZ)$ define $R_\gamma(f)$ to be the coefficient of $b_\gamma|_{Y_{v_\gamma}}$
in this expansion at $v_\gamma$. Put $R(f)=(R_\gamma(f))_\gamma$. The coefficients of the
Laurent expansions tend to zero, and there are finitely many vertices,
so $R:B_\partial\to E$ is well defined. Note that $R$ retains the coefficients with both exponents positive at each vertex, the positive powers of $z_i$
at the chosen endpoint of each edge $e_i$, and the constant coefficient
at $v_0$. 

\begin{lemma}[Boundary coefficients]\label{lem:boundary-coefficients}
Both $R:B_\partial\to E$ and $R\circ r_\partial:A_P\to E$ are isometric
linear embeddings.
\end{lemma}
\begin{proof}
Each closed edge is contained in a regular affine chart: the developing
map on a thin simply connected neighborhood is injective on the straight
edge and hence on a smaller neighborhood. We can therefore apply
Theorem~\ref{thm:local-regular-reconstruction} to the whole edge.
The $Y_{e_i}$ form an admissible affinoid cover of $U_\partial$ with
intersections at the common vertex domains. In particular,
$B_\partial$ is the closed equalizer of the edge-to-vertex restriction
maps. Each edge domain and vertex overlap is connected, so
$U_\partial$ is connected; it is smooth by regular reconstruction.
Along an edge the log Gauss norm is the supremum of the affine log
norms of the Laurent terms, and consequently its maximum is attained
at an endpoint. Thus
\[
 \|f\|_{U_\partial}=\max_{v\in V(P)}\|f|_{Y_v}\|_v. \tag{BC}
\]
Since we use orthonormal Laurent bases at the vertices, we already have
$\|R(f)\|\leq\|f\|_{U_\partial}$. For the other inequality, scale a
nonzero $f$ to have norm one and take a vertex $v$ where its norm is
one. Some coefficient in its Laurent expansion there has norm one.
We claim that the corresponding monomial belongs to the cone at $v$,
and, if its label lies on a ray, that its residue coefficient can be
read at either endpoint of the edge.

To see this, consider the reduction of the restriction of $f$ to
an incident edge. By \eqref{eq:regular-symbol-product}, a normalized
Laurent character on an edge restricts to a vertex with a nonzero
leading symbol only if its norm is attained at that vertex. In the
latter case its residue coefficient is unchanged. In the bases at
the vertices constructed above the residue of each monomial agrees,
by Lemma~\ref{lm:finite-tropical-lift}, with that of its normalized
character. The tropical marking identifies these residues with those
in the edge presentations. Consequently, a term with nonzero residue
at $v$ cannot have increasing log norm in either of the two directions
along the incident edges. The negative of its exponent therefore
belongs to the cone at $v$. By the Delzant condition it is of the
form $a\xi_i+b\xi_{i+1}$ with $a,b\geq0$ integral.

If $a,b>0$, its coefficient is one of those recorded by $R$ at $v$.
If the label is $k\xi_i$, $k>0$, its character has constant norm
along $e_i$. Its residue coefficient in the edge presentation is
therefore carried unchanged to the other endpoint. In particular
$|R_{k\xi_i}(f)|=1$ whichever endpoint was chosen. The same argument
for the constant character carries its residue coefficient around
$\partial P$, so if the label is zero then $|R_0(f)|=1$. In all cases
$\|R(f)\|=1$ as required.

Finally, by (BC) and the vertex norm formula \textup{(VN)} of
Section~\ref{sec:detection}, we have
$\|r_\partial x\|_{U_\partial}=\|x\|_P$ for $x\in A_P$.
\end{proof}

We shall use the following fact here and again for restriction to strips.
Write $\mathbf e_i$ for the $i$th coordinate vector of a $c_0$ space.

\begin{lemma}[Near-identity basis criterion]\label{cyEpsGoodGenCrit}
Let $J:H\to c_0(I,\Lambda)$ be an isometric embedding of a Banach
space and $u_i\in H$ elements of norm at most one. If
$\sup_i\|Ju_i-\mathbf e_i\|\le q<1$, then
$(c_i)\mapsto\sum_i c_i u_i$ is an isometric isomorphism onto $H$.
\end{lemma}
\begin{proof}
The series defines a contraction $L$. Since $\|JL-1\|\le q$,
the Neumann series makes $JL$ an isometric bijection. Hence $J$
is onto, and $L$ is an isometric bijection as well.
\end{proof}

\begin{proof}[Proof of \ref{part:polygon-basis}--\ref{part:polygon-boundary-functions}
of Theorem~\ref{tmAffinoidPolygon}]
We compare the maps
\[
 E\xrightarrow{\,S\,} A_P\xrightarrow{\,r_\partial\,}B_\partial
       \xrightarrow{\,R\,}E.
\]
At a vertex where $b_\gamma$ is supported, its restriction is exactly
a Laurent basis element in the coordinates used to define $R$.
Its only nonzero coefficient there is the one with label $\gamma$,
and this coefficient is recorded at the chosen vertex $v_\gamma$.
At any other vertex at least one edge occurring in the monomial is
not incident to the vertex. By Lemma~\ref{lm:finite-tropical-lift},
the corresponding factor has norm at most $q<1$, while all the other
factors have norm at most one. Hence
\[
 \|Rr_\partial(b_\gamma)-\mathbf e_\gamma\|\le q
 \qquad(\gamma\in P_{trop}(\bZ)).
\]
Here a common ray is recorded just once and the constant monomial is
recorded at $v_0$. Lemma~\ref{cyEpsGoodGenCrit} applied to $J=Rr_\partial$
proves \ref{part:polygon-basis} and gives $Rr_\partial(A_P)=E$.
Since $R$ is injective, $r_\partial(A_P)=B_\partial$.
Both $R$ and $Rr_\partial$ preserve norms and $r_\partial$ respects
products, proving \ref{part:polygon-boundary-functions}.

The norm-one lifts define a contraction from the Tate algebra.
Represent each $b_\gamma$ by its monomial $Z^{\alpha(\gamma)}$.
The monomials are distinct, with a common ray listed once, so
$\sum c_\gamma Z^{\alpha(\gamma)}$ lifts $\sum c_\gamma b_\gamma$
with the same norm. This gives surjectivity and identifies the
quotient norm with $\|\cdot\|_P$. The Gauss norm at every vertex is
multiplicative. Formula \textup{(VN)} gives
$\|f^n\|_P=\|f\|_P^n$, hence the spectral norm agrees as well.

For a set $I$ of edges with no common vertex, at every vertex one
factor of the following product has norm at most $q$. Thus
\begin{equation}\label{eq:nonface-bound}
 \left\|\prod_{i\in I}z_i\right\|_P\le q<1 .
\end{equation}
In the canonical reduction these products vanish. Monomials whose
edges do have a common vertex are precisely the $b_\gamma$.
Their reductions form a $\kk$-basis by \ref{part:polygon-basis}.
Hence the reduction is the face algebra of the labeled cycle, with
one coordinate plane for each maximal face $\{e_i,e_{i+1}\}$.
The dimension theorem for strict affinoid reduction
\cite[Reminder~1.3, p.~87]{LutkebohmertExtension} gives
$\dim A_P=\dim\widetilde A_P=2$. This proves \ref{part:polygon-affinoid}.
\end{proof}
\subsection{Extension and normality}\label{sec:hartogs}

\begin{lemma}\label{lem:normal-functions}
For a connected normal rigid analytic space $U$, the ring
$\Gamma(U,\cO)$ is an integrally closed domain.
\end{lemma}
\begin{proof}
By the identity theorem \cite[Lemma~2.1.4]{ConradIrreducible}, a
nonzero global function is nonzero on every nonempty connected
affinoid subdomain. The global ring is therefore a domain. If a
fraction $a/b$ is integral over it, the same monic equation on any
connected normal affinoid subdomain $V$ gives $a|_V/b|_V\in\cO(V)$.
These functions agree on overlaps, since multiplication by $b$ is
injective on each component there. They glue to a function on $U$.
\end{proof}

\begin{proof}[Proof of \ref{part:polygon-extension} and \ref{part:polygon-normal}
of Theorem~\ref{tmAffinoidPolygon}]
The identity theorem on $U$ makes $j^*$ injective. The composite
$j^*r_U=r_\partial$ is an isomorphism by
\ref{part:polygon-boundary-functions}, so both maps are isomorphisms.
Applying Lemma~\ref{lem:normal-functions} to the connected smooth
space $U_\partial$ shows that $A_P\simeq B_\partial$ is an integrally
closed domain.
\end{proof}

We have identified the functions on the regular boundary with $A_P$.
To identify that boundary as an analytic domain of $X_P$ we compute
Floer restriction to the edge strips in Section~\ref{sec:localization}.

\section{Shallow localization}\label{sec:localization}

Let $P$ satisfy the hypotheses of Theorem~\ref{tmAffinoidPolygon}.
We have identified $A_P=\cF^0(P)$ with the functions on the regular
mirror over $\partial P$. To glue that mirror to $X_P=\operatorname{Sp}A_P$
we must identify the domains on which their functions agree. We do this
by cutting a thin strip from each edge of $P$. The two sides of the cut
will be distinguished by the absolute value of the same edge lift.

Choose the normalized edge lifts $z_i$ of
Lemma~\ref{lm:finite-tropical-lift}. In the regular collar of $e_i$,
let $x_i$ be the inward affine coordinate of Section~\ref{sec:affinoid},
so the edge is $x_i=0$. For small $d_i\geq0$ move this edge to $x_i=d_i$.
Write $Q$ for the resulting polygon and $S_i=S_i(d_i)$ for the closed
strip $0\leq x_i\leq d_i$ in $P$ when $d_i>0$. At a vertex shared by
$e_i$ and $e_j$ the new vertex has coordinates $(d_i,d_j)$.
For a regular polytope $H$, write $Y_H=\operatorname{Sp}\cF^0(H)$
for the affinoid domain in the regular mirror with Berkovich
realization $p_{\rm reg}^{-1}(H)$. Write $\operatorname{res}_{P,H}:A_P\to A_H$ for
Floer restriction whenever $H\subset P$.

\begin{theorem}[Shallow localization]\label{thm:shallow-localization}
For all sufficiently small cut vectors $d=(d_1,\ldots,d_N)$ the
following hold. The Tate-algebra quotients carry their separated quotient norms.
\begin{enumerate}[label=\textup{(\roman*)},ref=(\roman*)]
\item\label{part:shallow-inward}
Restriction to $Q$ induces the isometric isomorphism
\[
 \psi: W_d:=A_P\langle X_1,\ldots,X_N\rangle/(z_i-T^{d_i}X_i)_i
     \xrightarrow{\ \sim\ }A_Q,\qquad
 X_i\longmapsto w_i:=T^{-d_i}\operatorname{res}_{P,Q}(z_i).
\]
Thus $X_Q$ is the domain $D_d=\{|z_i|\leq e^{-d_i}\ \text{for all }i\}$
in $X_P$.
\item\label{part:shallow-strip}
For $d_i>0$, restriction to $S_i$ induces the isometric isomorphism
\[
 L_{i,d_i}:=A_P\langle X\rangle/(z_iX-T^{d_i})
       \xrightarrow{\ \sim\ } A_{S_i},\qquad
 X\longmapsto T^{d_i}/\operatorname{res}_{P,S_i}(z_i).
\]
It identifies $Y_{S_i}$ with the Laurent domain
$V_i=\{|z_i|\geq e^{-d_i}\}$ in $X_P$.
\item\label{part:shallow-corners}
The domains $D_d$ and the positive-width $V_i$ are affinoid and
cover $X_P$. For adjacent $e_i,e_j$ of positive strip widths,
restriction identifies $V_i\cap V_j$ with $Y_{C_{ij}}$, where
$C_{ij}=S_i\cap S_j=[0,d_i]\times[0,d_j]$ in the vertex coordinates.
The corresponding isomorphism of affinoid algebras is isometric for
the quotient norm. A collection of strip domains whose edges have no
common vertex has empty intersection.
\item\label{part:shallow-facets}
The following hold for the interfaces and the boundary collar.
\begin{enumerate}[label=\textup{(\alph*)}]
\item\label{part:shallow-new-facet}
For each $d_i>0$,
\[
 D_d\cap V_i
 =\operatorname{Sp}\bigl(A_Q\langle Y\rangle/(w_iY-1)\bigr)
 =Y_{Q\cap\{x_i=d_i\}}.
\]
Thus inside $X_Q=D_d$, meeting the strip $V_i$ cuts out the
regular mirror over the new $i$th edge of $Q$.
\item\label{part:shallow-partial-Q}
If every width is positive, these facet domains cover the
regular mirror over $\partial Q$.
\item\label{part:shallow-collar}
If every width is positive, the strip domains $V_i\simeq Y_{S_i}$
glue along their corner overlaps to a thin regular collar of
$\partial P$ inside $X_P$. In particular they embed
$U_\partial=\bigcup_iY_{e_i}$, the regular mirror over
the union of the edges of $P$, as an analytic domain of $X_P$.
\item\label{part:shallow-compat}
All of the identifications above commute with Floer
restriction to common refinements.
\end{enumerate}
\end{enumerate}
\end{theorem}

\subsection{Restriction to the polygon and to a strip}

Fix the testing polygons of Lemma~\ref{lm:finite-tropical-lift}.
For each edge, those that cover it contain a relative neighborhood
of the whole edge in $P$. A sufficiently thin $S_i$ lies in this
neighborhood, and each displaced vertex lies in its original vertex
cap. The cyclic order and the primitive conormals persist. Thus $Q$
is admissible, and each $S_i$ of positive width is a regular Delzant
quadrilateral in the chart along its edge. Its outer side lies on
$e_i$ and its short sides on the two adjacent edges, so $S_i\subset P$
is facewise admissible. Take $0<q<1$ bounding the errors of
Lemma~\ref{lm:finite-tropical-lift} and the products in
\eqref{eq:nonface-bound}, and require
\begin{equation}\label{eq:shallow-smallness}
             \sum_i d_i<-\log q .
\end{equation}

On one of the testing polygons meeting $e_j$, write $m_j$ for the
primitive character in its chosen marked presentation. There
$\|z_j-m_j\|\leq q$ and $|m_j|=e^{-x_j}$. Thus, at any analytic
point of that polygon and for any $\lambda>q$, each of the conditions
$|z_j|\geq\lambda$, $|z_j|\leq\lambda$ and $|z_j|=\lambda$ is
equivalent to the same condition on $|m_j|$. This follows from the
ultrametric inequality. On a testing polygon disjoint from $e_j$,
$|z_j|\leq q$. We will use this observation both for the edge lifts
after a cut and for the intersections of the resulting domains.

\begin{lemma}\label{lem:collar-retraction}
The polygon $Q$ satisfies $(T)_j$ in every degree $j$ in which
$P$ does. Each $w_i$ is a normalized lift of the corresponding edge
of $Q$. At a vertex $v'$ of $Q$, in the chart of the corresponding
cap of $P$, these lifts satisfy
\[
 \begin{aligned}
 \|w_i-T^{-d_i}m_i\|_{v'}&\leq qe^{d_i}
       &&\text{if the new $i$th edge contains $v'$},\\
 \|w_i\|_{v'}&\leq qe^{d_i} &&\text{otherwise}.
 \end{aligned}
\]
\end{lemma}
\begin{proof}
Retraction across the regular collar lifts through the torus bundle
to a deformation retraction $K_P\to K_Q$. Parallel transport identifies
the maps on ordinary cohomology to a boundary fiber, so $(T)$ for
$P$ gives $(T)$ for $Q$.

Along the new $i$th edge the character $m_i$ has absolute value
$e^{-d_i}$. Thus $T^{-d_i}m_i$ is its normalized character, and the
error on an incident testing piece is at most $qe^{d_i}$. The same
bound applies to $w_i$ near any displaced vertex where $e_i$ is
nonincident. The vertex norm formula \textup{(VN)} of
Section~\ref{sec:detection} gives $\|w_i\|_Q=1$.
\end{proof}

We will use the basis of Theorem~\ref{tmAffinoidPolygon} with
prescribed edge lifts. For a polygon $H$ satisfying that theorem,
let $u_l$ be its chosen lifts from Lemma~\ref{lm:finite-tropical-lift}.
Suppose that $v_l\in A_H$ have the same reductions as $u_l$ at
every vertex. These reductions are the normalized primitive character
at an incident vertex and zero at every other vertex. By \textup{(VN)},
$\eta=\max_l\|v_l-u_l\|_H<1$. All these lifts have norm one,
so replacing the factors one at a time gives
$\|v^{\alpha(\gamma)}-u^{\alpha(\gamma)}\|_H\leq\eta$
for every label $\gamma$. The linear operator sending the old
basis monomials to the new ones differs from the identity by norm
at most $\eta$. Its Neumann series gives an inverse of norm at
most one. Thus the prescribed lifts give an orthonormal basis as well.

\begin{proof}[Proof of Theorem~\ref{thm:shallow-localization}\,\ref{part:shallow-inward}]
By Lemma~\ref{lem:collar-retraction} the classes $w_i$ are normalized
edge lifts for $Q$. The cut keeps the cyclic edge order and the same
primitive conormals, so $Q$ has the same tropical label set
$P_{trop}(\bZ)$ as $P$. The vertex estimates of the lemma and the
basis comparison above give an orthonormal Schauder basis of $A_Q$ by the
monomials $w^{\alpha(\gamma)}$, $\gamma\in P_{trop}(\bZ)$, where
$\alpha(\gamma)$ records the corresponding powers of the edge lifts
exactly as for the basis $b_\gamma=z^{\alpha(\gamma)}$ of $A_P$.

Since $\|w_i\|_Q=1$, the map $\psi:W_d\to A_Q$ is a contraction. Sending each basis monomial $w^{\alpha(\gamma)}$ to
$X^{\alpha(\gamma)}$ gives a contractive linear right inverse $s$.
For a source basis monomial $b_\gamma=z^{\alpha(\gamma)}$ the
relations of $W_d$ give
\[
 s\psi(b_\gamma)=T^{\sum_i d_i\alpha_i(\gamma)}
                    X^{\alpha(\gamma)}=b_\gamma.
\]
The image of $A_P$ is dense in $W_d$, since it contains every
$X_i=T^{-d_i}z_i$ and hence the polynomial algebra over the image
of $A_P$. Continuity gives $s\psi=1$. So, the inverse contractions
are isometries.
\end{proof}

\begin{proof}[Proof of Theorem~\ref{thm:shallow-localization}\,\ref{part:shallow-strip}]
Fix $i$, write $d=d_i>0$ and $S=S_i(d)$, and denote the restrictions
of the $z_l$ by the same letters. On the regular pieces covering $S$
we have $|z_i|=|m_i|=e^{-x_i}\geq e^{-d}$. Thus $T^d/z_i$ defines a
function of norm at most one on the boundary of $S$. Since the torus bundle over the convex strip retracts onto a
fiber, $S$ satisfies $(T)$. By Theorem~\ref{tmAffinoidPolygon}\,\ref{part:polygon-boundary-functions}
this boundary function extends to $x\in A_S$ with $\|x\|_S\leq1$.
The identity $z_ix=T^d$ holds on the boundary and hence in $A_S$.
We obtain a contraction $\psi:L_{i,d}\to A_S$ with $X\mapsto x$.

Write $j=i-1$ and $k=i+1$ in cyclic order. View $S$ as an
admissible quadrilateral with edges in the order of $e_i$, the
short side along $e_j$, the inner cut, and the short side along
$e_k$. We claim that $z_i$, $z_j$, $x$, and $z_k$ are normalized
edge lifts for these four sides. On testing pieces meeting $e_i$
one has $z_i=m_i(1+h_i)$ with $\|h_i\|\le qe^d<1$, so along the
inner edge
\[
 x=\frac{T^d}{z_i}=\frac{T^d}{m_i}(1+h_i)^{-1}
\]
matches the normalized inverse character of $m_i$. Moreover
$\|z_i\|$ on the inner edge and $\|x\|$ on the outer edge both
equal $e^{-d}<1$. On the short sides, the estimates for $z_j$ and
$z_k$ are those of Lemma~\ref{lm:finite-tropical-lift} at the
endpoint caps. Thus at every vertex of $S$ the off-edge errors are
at most $\max(qe^d,e^{-d})<1$.
The same basis comparison gives an orthonormal Schauder basis of $A_S$ by the
monomials in these four lifts.

Lift each of their basis monomials to the same expression in
$z_i,z_j,X,z_k$ in $L_{i,d}$. This gives a contractive linear section
$s$ of $\psi$. To prove injectivity put $E=1-s\psi$, a contractive
idempotent. Before quotienting, $A_P\langle X\rangle$ has orthonormal
basis $b_\gamma X^\ell$. If $\gamma$ lies in either vertex cone on
$e_i$, cancel pairs $z_iX=T^d$. The result is a scalar multiple of
one of the chosen strip monomials, so $E$ vanishes on its image.
For every other $\gamma$, at each vertex of $P$ some factor of
$z_i b_\gamma$ is nonincident. Formula \textup{(VN)} gives
$\|z_i b_\gamma\|_P\leq q$. In the quotient,
\[
 b_\gamma X^\ell=T^{-d}(z_i b_\gamma)X^{\ell+1},
 \qquad \|b_\gamma X^\ell\|_{L_{i,d}}\leq qe^d.
\]
For the quotient map $\pi:A_P\langle X\rangle\to L_{i,d}$ this proves
$\|E\pi(f)\|\leq qe^d\|f\|$. Taking the infimum over representatives
gives $\|E\|\leq qe^d<1$. Since $E^2=E$, it follows that $E=0$.
This proves the isometry. The regular reconstruction
(Theorems~\ref{thm:local-regular-reconstruction}
and~\ref{thm:global-reconstruction}) identifies $\operatorname{Sp}A_S$
with $Y_S$ by the same restrictions on every regular subpolytope of $S$.
\end{proof}

\subsection{The intersections in the regular mirror}

\begin{proof}[Proof of Theorem~\ref{thm:shallow-localization}\,\ref{part:shallow-corners}]
A point outside $D_d$ satisfies $|z_i|>e^{-d_i}$ for some $d_i>0$,
and hence lies in $V_i$. Each of $D_d$ and the $V_i$ is affinoid by
parts~\ref{part:shallow-inward} and~\ref{part:shallow-strip}, so
these domains are a finite affinoid cover of $X_P$. If the edges
indexed by $I$ have no common vertex, a point in all the $V_i$ for
$i\in I$ would satisfy
\[
 q\ \geq\ \left|\prod_{i\in I}z_i\right|
       \ \geq\ e^{-\sum_{i\in I}d_i}\ >q
\]
by \eqref{eq:nonface-bound} and \eqref{eq:shallow-smallness}.
This also excludes a triple intersection for a triangle.

For adjacent edges, compute $V_i\cap V_j$ on $Y_{S_i}$.
The test domains covering $e_i$ can be taken to be its two endpoint caps
and polygons meeting only $e_i$. On every one except the cap at
$e_i\cap e_j$, $|z_j|\leq q<e^{-d_j}$. On that cap the inequality
$|z_j|\geq e^{-d_j}$ is exactly $x_j\leq d_j$. For small widths
the whole rectangle $C_{ij}$ lies in the cap, so this identifies
the Laurent domain with $Y_{C_{ij}}$. By uniqueness of the affinoid
domain and regular reconstruction the algebra map is restriction
\[
 A_P\langle X_i,X_j\rangle/
      (z_iX_i-T^{d_i},z_jX_j-T^{d_j})
          \ \xrightarrow{\ \sim\ }\ A_{C_{ij}}.
\]

For the claim concerning the  quotient norm, work in the marked cap chart and put
$u_l=z_l|_{C_{ij}}=m_l(1+h_l)$ for $l=i,j$. Then
$\|h_l\|\leq qe^{d_l}<1$. The rectangle has a normalized Laurent
basis obtained by taking powers of $m_l$ or $T^{d_l}/m_l$ in each
of the two directions. Replace them by $u_l$ and $T^{d_l}/u_l$.
For every integer $n$, $\|(1+h_l)^n-1\|\leq\|h_l\|$, using the
geometric series for the inverse when $n<0$. Consequently each new
monomial differs from its basis monomial by at most
$q\max(e^{d_i},e^{d_j})<1$. Lemma~\ref{cyEpsGoodGenCrit} gives an
orthonormal basis of $A_{C_{ij}}$. Every one of its elements lifts
to a monomial of quotient norm at most one. Expanding in this basis
bounds the inverse by one. Restriction itself is contractive since
the four generators are power bounded, proving the isometry.
\end{proof}

\begin{proof}[Proof of Theorem~\ref{thm:shallow-localization}\,\ref{part:shallow-facets}]
For \ref{part:shallow-new-facet}: on $D_d$ one has $|w_i|\leq1$, so
the further inequality defining $V_i$ is $|w_i|=1$. This is the
displayed Laurent algebra. On the strip $Y_{S_i}$ the same equality
is $x_i=d_i$. The remaining inequalities are $|z_j|\leq e^{-d_j}$.
On the cap at an adjacent edge of positive width this means
$x_j\geq d_j$; elsewhere on the strip it is automatic, by the same
testing estimates as in \ref{part:shallow-corners}. A zero width
imposes no condition. This cuts out precisely the new $i$th facet
of $Q$.

For \ref{part:shallow-partial-Q}: when every width is positive these
facets are the full boundary of $Q$, so their mirrors cover the
regular mirror over $\partial Q$.

For \ref{part:shallow-collar}: choose all widths positive. The
domains $V_i$ are the regular strips of a thin collar of
$\partial P$, and their overlaps are the corner rectangles of
\ref{part:shallow-corners}. In each strip the equality $|z_i|=1$
cuts out the original edge $x_i=0$. The union of the strips is
therefore a regular collar in $X_P$, and the union of these outer
faces is $U_\partial$.

For \ref{part:shallow-compat}: the maps on functions extend Floer
restriction, and the relations determine the images of the
auxiliary variables, so the identifications agree on every common
refinement. The isometries above also identify the quotient norms
on the polygon, strip and corner algebras with their spectral norms.
\end{proof}

\section{Higher-degree Floer cohomology}\label{sec:nonsmooth}

Let $P$ satisfy $(G_\partial)$, $(T)_0$ and $(T)_1$. Write
$A=A_P$, $X=X_P=\operatorname{Sp}A$ and $U=U_\partial$.
By Theorems~\ref{tmAffinoidPolygon} and~\ref{thm:shallow-localization},
$X$ is a normal affinoid surface, $U$ is a smooth analytic domain of
$X$, and restriction identifies $A$ with $\Gamma(U,\cO_U)$.
We use this identification to extend polyvector fields from $U$ to
$X$. The remaining question is which of them come from Floer
cohomology. In degree two we will construct a Floer lift of the
normalized logarithmic bivector. Its brackets with functions will
generate the tangent sheaf on the smooth locus.

On a normal analytic surface write
$\PV_X^{[q]}=(\wedge^qT_X)^{**}$ for the sheaf of reflexive
$q$-vectors in degree $q\geq0$ and set it equal to zero in negative
degrees. Its sections are determined by their restriction
to the smooth locus. With the conventions of \eqref{eq:bracket},
the \emph{bracket symbol} of a Floer class is given by
\[
 \Sigma_0(f)=f,\qquad
 \Sigma_1(v)(f)=\{v,f\},\qquad
 \Sigma_2(w)(f,g)=-\{\{w,f\},g\}.
\]
Here $f,g\in A$, $v\in\cF^1(P)$ and $w\in\cF^2(P)$. The target
is zero in all other degrees.

\begin{theorem}[Higher-degree Floer cohomology of a polygon]
\label{thm:polygon-higher-degree}
Suppose that $P$ satisfies $(G_\partial)$, $(T)_0$ and $(T)_1$.
Fix an orientation of the boundary collar.
\begin{enumerate}[label=\textup{(\roman*)},ref=(\roman*)]
\item\label{part:polygon-concentration}
$\cF^q(P)=0$ for $q\notin\{0,1,2\}$.
\item\label{part:polygon-bv}
The regular logarithmic volume form on $U$ extends uniquely to a
generator $\Omega$ of the canonical sheaf of $X$. In particular,
$X$ is Gorenstein. With divergence defined by $\Omega$ on the
smooth locus, the bracket symbol
\[
 \Sigma:\cF^*(P)\longrightarrow\Gamma(X,\PV_X^{[*]})
\]
is a homomorphism of BV algebras and is the identity in degree zero.
\item\label{part:polygon-top}
Let $\Pi$ be the reflexive bivector dual to $\Omega$ on the smooth
locus and put $Z_P=\ker\Sigma_2$. There is a split exact sequence of
$A$-modules
\begin{equation}\label{eq:top-canonical}
 0\longrightarrow Z_P\longrightarrow\cF^2(P)
 \xrightarrow{\Sigma_2}A\Pi\longrightarrow0 .
\end{equation}
Every lift $s_P$ of $\Pi$ satisfies $\Delta s_P=0$. The kernel
$Z_P$ is square-zero and $\Delta Z_P=0$.
The ideal $J_\Pi=(\Pi(df,dg)\mid f,g\in A)$ annihilates $Z_P$.
In particular $Z_P\otimes_A B=0$ on every smooth affinoid
domain $\operatorname{Sp}B\subset X$.
If $(T)_2$ also holds, then $Z_P=0$ and $\cF^2(P)=As_P$ with
$s_P$ unique.
\item\label{part:polygon-defects}
In degree one there is an exact sequence
\begin{equation}\label{eq:defect}
 0\longrightarrow\cF^1(P)\xrightarrow{\Sigma_1}
 D_A\longrightarrow Q_P\longrightarrow0,
 \qquad D_A=\Der_{\Lambda,\mathrm{cont}}(A).
\end{equation}
The module $Q_P$ has finite length and is supported on the analytic
singular locus of $X$. Thus $(\cF^1(P))^{**}\simeq D_A$.
If $X$ is smooth, $\Sigma$ is an isomorphism in all degrees.
\end{enumerate}
\end{theorem}

\subsection{Extension from the regular boundary}

Regular reconstruction identifies the restrictions of Floer classes
with polyvector fields on $U$. Write
\[
 r_q:\cF^q(P)\longrightarrow\Gamma(U,\wedge^qT_U)
\]
for these maps. The vertex norm formula \textup{(VN)} of
Section~\ref{sec:detection} gives injectivity in degrees zero and
one. By Theorem~\ref{thm:bv-reconstruction} the chosen orientation gives
a logarithmic volume form $\Omega_U$ on $U$, and $r_*$ respects
products, brackets and divergence. Write $\Pi_U$ for the bivector
dual to $\Omega_U$.

\begin{lemma}\label{lem:polyvector-extension}
Restriction induces an isomorphism
\[
 \Gamma(X,\PV_X^{[q]})\xrightarrow{\ \sim\ }
 \Gamma(U,\wedge^qT_U).
\]
Under this isomorphism the bracket symbol $\Sigma_q$ corresponds to $r_q$.
\end{lemma}

\begin{proof}
It suffices to consider $q\geq0$. Let $\Omega_A$ be the module of
analytic differentials of $A$ over $\Lambda$.
The finite module $\operatorname{Hom}_A(\wedge^q\Omega_A,A)$ is
reflexive \cite[Tag~0AV3]{Stacks}. On the smooth locus it is the
module of $q$-vectors, so its associated sheaf is $\PV_X^{[q]}$ by
normality. Choose a finite presentation
$A^m\xrightarrow{L}A^n\to\wedge^q\Omega_A\to0$. Then
\[
 \Gamma(X,\PV_X^{[q]})=\ker(L^*:A^n\to A^m).
\]
On affinoid subdomains of $U$, analytic differentials are obtained
by pullback from $X$ and affinoid localization is flat
\cite[Theorem~3.1.1, Example~3.1.4 and Exercise~3.2.10]{Conrad2008}.
Thus $\wedge^qT_U$ is the kernel of the same matrix over $\cO_U$.
Taking sections preserves kernels, and $\Gamma(U,\cO_U)=A$.

On $U$ the evaluation of $r_q(x)$ on differentials of functions
from $A$ is the signed iterated bracket defining $\Sigma_q(x)$.
These evaluations determine the extended polyvector, since
differentials of Tate generators of $A$ generate $\Omega_A$.
\end{proof}

Since $\Pi_U$ is a frame and $\Gamma(U,\cO_U)=A$, the lemma gives
$\Gamma(X,\PV_X^{[2]})=A\Pi$, where $\Pi|_U=\Pi_U$. On the affinoid
$X$ this module generates the associated reflexive sheaf, so $\Pi$ is
a frame on the smooth locus. There its dual $\Omega$ trivializes the
canonical bundle. A normal surface is
Cohen--Macaulay, and its canonical sheaf is the reflexive extension
of that bundle. Indeed, for a Tate presentation
$T_N=\Lambda\langle t_1,\ldots,t_N\rangle\twoheadrightarrow A$,
the canonical module is
$\operatorname{Ext}_{T_N}^{N-2}(A,\Omega^{N,\mathrm{an}}_{T_N/\Lambda})$.
Thus $\Omega$ trivializes $\omega_X$ and $X$ is Gorenstein
\cite[Tags~0B3D and~0DW9]{Stacks}. Uniqueness follows from
$A=\Gamma(U,\cO_U)$. Apply divergence on the smooth locus to obtain
an operator on reflexive polyvectors. Its restriction to $U$ is
divergence for $\Omega_U$. The lemma and the compatibility of $r_*$
with the Floer operations now prove \ref{part:polygon-bv}.

\subsection{Lifting the logarithmic bivector}

The topology of the polygon supplies both the orientation class we
want to lift and the vanishing in the next degree needed to lift it.

\begin{lemma}\label{lem:polygon-fiber-class}
If $(T)_1$ holds and $L$ is a fiber in the boundary collar, then
$H^2(K_P;\kk)\to H^2(L;\kk)$ is surjective. If $(T)_0$ also holds,
then $H^j(K_P;\kk)=0$ for $j\geq3$.
\end{lemma}

\begin{proof}
Take a smooth outward rounding $D$ of $P$ as in
Proposition~\ref{prpCornerSmoothing}. The region between the
boundaries is a regular collar, so $K=K_D$ retracts onto $K_P$
and we can take $L$ over a point of $\partial D$.
The compact symplectic four-manifold $K$ has boundary
$i:Z=\partial K\hookrightarrow K$, a torus bundle
$p:Z\to\partial D\simeq S^1$. The pullback $\alpha$ of a
generator of $H^1(S^1;\kk)$ is nonzero since the fibers are
connected. Its restriction to $L$ vanishes. If $\alpha=i^*\beta$
for some $\beta\in H^1(K;\kk)$, then $\beta|_L=0$, so $(T)_1$
gives $\beta=0$. Thus the connecting homomorphism for $(K,Z)$
sends $\alpha$ to a nonzero class in $H^2(K,Z;\kk)$.

The Poincar\'e dual of $\alpha$ in $Z$ is the fiber class up to
sign. By Poincar\'e--Lefschetz duality, the dual of its connecting
image in $K$ is therefore $\pm[L]\in H_2(K;\kk)$. In particular
$[L]\ne0$. Over the field $\kk$ there is a class in $H^2(K;\kk)$
pairing to one with this oriented fiber.

Under $(T)_0$, $K$ is connected. Dualizing $(T)_1$ shows that
$H_1(L;\kk)\to H_1(K;\kk)$ is surjective. This map factors through
$H_1(Z;\kk)$. The exact sequence of $(K,Z)$ therefore gives
$H_1(K,Z;\kk)=H_0(K,Z;\kk)=0$, since $Z$ is connected and nonempty.
Poincar\'e--Lefschetz duality gives the assertions in degrees three
and four. Higher degrees vanish by dimension.
\end{proof}

For $j\notin\{0,1,2\}$ ordinary restriction is now injective
because its source vanishes. Corollary~\ref{cor:concentration}
proves \ref{part:polygon-concentration}. The lemma also shows that,
under our hypotheses, $(T)_2$ is equivalent to $b_2(K_P;\kk)=1$.

\begin{proof}[Proof of Theorem~\ref{thm:polygon-higher-degree}\,\ref{part:polygon-top}]
Fix a rounding $D$ as in the lemma with $P\subset\operatorname{int}D$.
Choose $\eta\in H^2(K_D;\kk)$ restricting to the normalized
orientation class of a boundary fiber. Transport in the oriented
collar gives the same normalization at every vertex of $P$.

We will lift $\eta$ from the critical block in a short Novikov
window. For each vertex $v$, fix a sufficiently small smooth strictly convex
neighborhood $R_v\subset\operatorname{int}D$ of $v$ in a regular
affine chart. Fix finite period cutoffs below the first positive
Reeb periods of $K_D$ and the $K_{R_v}$. Take a log-norm window
$[0,\delta)$ short enough for these finitely many comparisons.
By Theorem~\ref{thm:background-indicator-input}, $\eta\otimes1$ in
the critical block gives a window class on $K_D$. Restriction to
$K_{R_v}$ sends it to $\eta|_{K_{R_v}}\otimes1$ in the critical
block \cite[Lemma~5.10]{Paper1}. Under restriction from $R_v$ to
$\{v\}$, the zero Laurent label has weight $T^0$ and its
coefficient is ordinary restriction to the fiber
\cite[Theorem~5]{Paper1}. Thus, in the oriented tropical marking,
the class at $v$ is the leading term of $\Pi_U|_{Y_v}$.
By functoriality, restricting from $D$ through $P$ gives this same
window class at every $v$.

Write $C$ for the completed Floer complex of $K_P$ and
$C_t=\{c:\|c\|<e^t\}$. The restriction of the window class to $P$
lies in $H^2(C_\delta/C_0)$. Consider the exact sequence
\[
 H^2(C_\delta)\longrightarrow H^2(C_\delta/C_0)
 \longrightarrow H^3(C_0)\longrightarrow H^3(C_\delta).
\]
Under Novikov translation the last arrow is multiplication by
$T^\delta$ on $H^3(C_0)$, which is injective by
Theorem~\ref{thm:strict-detection} in degree three, where
$(T)_3$ follows from Lemma~\ref{lem:polygon-fiber-class}. Lift the window
class to $H^2(C_\delta)$ and let $w$ denote its image in
$\cF^2(P)$. At a vertex $v$ the critical-block coefficient $1$
comes from every subwindow $[0,\epsilon)$, $0<\epsilon<\delta$. As in the
proof of Lemma~\ref{lm:finite-tropical-lift}, injectivity of the
regular open lattices lets us lift these subwindow classes with norm
less than $e^\epsilon$. Their differences from $w|_v$ come from
$H^2(C_0(v))$, so $\|w|_v\|_v\leq1$. The marked reduction is now
defined. It is the zero-label orientation class, which is also the
reduction of $\Pi_U|_{Y_v}$. Consequently their difference has
norm strictly less than one.

Write $\Sigma_2(w)=a\Pi$. Since $\Pi|_{Y_v}$ is the logarithmic
frame of norm one, the preceding comparison gives
$\|(a-1)|_{Y_v}\|_v<1$ at each vertex. There are finitely many
vertices, so \textup{(VN)} in degree zero gives $\|a-1\|_P<1$.
Thus $a$ is a unit and $s_P=a^{-1}w$ has symbol $\Pi$. Any such lift
splits \eqref{eq:top-canonical} by $f\Pi\mapsto fs_P$. For every lift,
$\Sigma_1(\Delta s_P)=\diver_\Omega\Pi=0$, since $\iota_\Pi\Omega=1$.
Injectivity of $\Sigma_1$ gives $\Delta s_P=0$.

For $z\in Z_P$ we likewise have $\Sigma_1(\Delta z)=0$ and hence
$\Delta z=0$. The product of $z$ with every positive-degree class
vanishes by concentration. Let $f,g\in A$. The bracket $\{z,g\}$ has
zero symbol and degree one, hence vanishes. Put $v=\{s_P,f\}$.
Applying the bracket with $g$ to $zv=0$ and using the derivation rule gives $z\{v,g\}=0$. By the definition of the
symbol, $\{v,g\}=-\Pi(df,dg)$. It follows that $J_\Pi Z_P=0$.
At each smooth point, differentials of Tate generators span the
cotangent space and $\Pi$ is nondegenerate. Some $\Pi(df,dg)$ is
therefore nonzero at that point. On a smooth affinoid domain
$\operatorname{Sp}B$ this gives $J_\Pi B=B$, so $Z_P\otimes_A B=0$.
If $(T)_2$ holds, \textup{(VN)} makes $r_2$ injective and
Lemma~\ref{lem:polyvector-extension} gives $Z_P=0$.
\end{proof}

\begin{proof}[Proof of Theorem~\ref{thm:polygon-higher-degree}\,\ref{part:polygon-defects}]
Let $N=\Sigma_1(\cF^1(P))\subset D_A$. For every $f\in A$,
the Hamiltonian vector field $[\Pi,f]_{\mathrm{Sch}}$ belongs to
$N$, since it is the symbol of $\{s_P,f\}$. On the smooth locus
contraction with $\Pi$ identifies one-forms with vector fields,
and differentials of Tate generators generate the cotangent
sheaf. The corresponding Hamiltonian fields therefore generate
the tangent sheaf there.

It follows that $Q_P=D_A/N$ is supported on the singular locus.
The ring $A$ is Noetherian, so $N$ and $Q_P$ are finite.
In characteristic zero the singular locus of a normal affinoid
surface is closed and has dimension at most zero, hence $Q_P$
has finite length. Thus $N$ agrees with $D_A$ in codimension one.
Since $D_A$ is reflexive, $N^{**}=D_A$ \cite[Tag~0AVB]{Stacks}.
\end{proof}

\subsection{Restriction to regular domains}

\begin{proposition}[Base change on regular overlaps]
\label{prop:regular-base-change}
Under the hypotheses of
Theorem~\ref{thm:polygon-higher-degree}, let $H\subset P$
be a regular polytope for which Floer restriction identifies
$Y_H=\operatorname{Sp}A_H$ with an affinoid domain of $X_P$.
Then restriction induces an isomorphism of $A_H$-modules
\[
 \cF^q(P)\otimes_{A_P}A_H
 \xrightarrow{\ \sim\ }\cF^q(H)
 \simeq\Gamma(Y_H,\wedge^qT_{Y_H})
\]
in every degree.

More generally, for any smooth affinoid domain
$V=\operatorname{Sp}B$ of $X_P$, the localized symbol
gives an isomorphism
\[
 \cF^q(P)\otimes_{A_P}B
 \xrightarrow{\ \sim\ }\Gamma(V,\wedge^qT_V).
\]
These identifications commute with further affinoid restriction.
\end{proposition}

\begin{proof}
Affinoid localization is flat, and analytic differentials pull
back to those of the domain. Thus the finite duals
$\operatorname{Hom}_A(\wedge^q\Omega_A,A)$ localize to polyvectors
on $V$. Tensoring \eqref{eq:defect} with $B$ kills $Q_P$, since its
support misses $V$. Tensoring \eqref{eq:top-canonical} with $B$
kills $Z_P$, since $J_\Pi B=B$. Together with degree zero and
concentration this gives the localized isomorphism in every degree.

For $V=Y_H$, the symbol of actual Floer restriction and the
localized symbol agree on differentials of functions from $A$ by
naturality of the bracket. Since these differentials generate
$\Omega^{1,\mathrm{an}}_{A_H/\Lambda}$ over $A_H$, the two maps
agree. Regular reconstruction identifies the Floer symbol over $H$
with all polyvector fields on $Y_H$. The same argument gives
compatibility with further restriction.
\end{proof}

\section{Extension over the singular value}\label{sec:radial}

Let $P$ satisfy $(G_\partial)$, $(T)_0$ and $(T)_1$ and contain the single
singular value $s$. Section~\ref{sec:localization} identifies the
regular collar of $\partial P$ in $X_P$. To extend the fibration
over $s$ we need to continue these identifications along a family
of polygons shrinking to $s$ in the base. For the nodal models,
our earlier construction with Varolg\"unes uses Weierstrass domains
over affine half-spaces
\cite[Proposition~3.15]{GromanVarolgunes2022}. For the $A_n$ smoothing,
Yuan constructs the mirror fibration by family Floer theory
\cite{YuanAn}.

Recall the full radial conicity condition $(G_{\rm rad})$ of
Definition~\ref{def:radial}.

Write $X_r^{\rm rig}=\operatorname{Sp}A_{P_r}$ for the rigid affinoid
and $X_r=\mathcal M(A_{P_r})$ for its Berkovich realization. Put
$X^{\rm rig}=X_1^{\rm rig}$ and $X=X_1$.

\begin{theorem}[Extension over the base singularity]\label{thm:radial-extension}
Assume $(G_{\rm rad})$ in addition to $(G_\partial)$, $(T)_0$ and $(T)_1$.
\begin{enumerate}[label=\textup{(\roman*)},ref=(\roman*)]
\item\label{part:radial-domains}
For $0<r<R\leq1$, Floer restriction $A_{P_R}\to A_{P_r}$
identifies $X_r^{\rm rig}$ with an affinoid domain of $X_R^{\rm rig}$.
This homomorphism is a finite composite of the Weierstrass
localizations of Theorem~\ref{thm:shallow-localization}%
\,\ref{part:shallow-inward}. Its image is dense, and $X_r$ lies in
the Berkovich analytic interior of $X_R$ over $\Lambda$.
\item\label{part:radial-maps}
The marked regular fibration extends to a continuous proper surjection
$p:X\to P$ with
\[
 p^{-1}(P_r)=X_r \qquad(0<r\leq1).
\]
Over $P\setminus\{s\}$ this is the marked regular projection.
For the corresponding rigid affinoids, the identifications
$\Gamma(X_r^{\rm rig},\cO)=A_{P_r}$ commute with Floer restriction.
\end{enumerate}
\end{theorem}

A shallow cut is available at each positive radius. To pass between
two prescribed radii we need a lower bound for the allowed step on
the intervening compact interval. Once we have this bound, the
regular strips from the successive cuts fill the punctured polygon.
The intersection of the remaining Berkovich domains will be the
fiber over $s$.

\subsection{Localization between positive radii}

Each $P_r$ inherits $(T)_j$ in every degree in which it holds
for $P$. Indeed, the annulus between
$\partial P$ and $\partial P_r$ lies in the regular torus bundle.
Lifting its radial retraction and fixing $K_{P_r}$ retracts $K_P$
onto $K_{P_r}$. Parallel transport identifies the ordinary
restrictions to their boundary fibers.

\begin{lemma}[Shallow cuts at positive radii]\label{lem:radial-uniformity}
For $0<a\leq b\leq1$ there is $\Delta>0$ such that, at each
$r\in[a,b]$, one can choose normalized edge lifts for $P_r$
for which every cut to $P_{r'}$ with $r'>0$ and $0<r-r'<\Delta$
satisfies the hypotheses of Theorem~\ref{thm:shallow-localization}.
\end{lemma}
\begin{proof}
Fix $0<r_0<1$. To compare the lifts at nearby radii on both sides of
$r_0$, we first approximate them by restrictions from a slightly
larger polygon. Put $C(r)=SC_M^*(K_{P_r};\Lambda)$. For a finite
log-norm window $I$,
\[
 \varinjlim_{R\downarrow r_0}H^*(C(R)_I)\simeq H^*(C(r_0)_I).
\]
Indeed, a Hamiltonian strictly negative on $K_{P_{r_0}}$ is strictly
negative on $K_{P_R}$ for $R>r_0$ sufficiently close to $r_0$, by
compactness. This gives cofinal Hamiltonians for the indicators.
In a finite window the completed telescope is the algebraic telescope,
as in the proof of Theorem~\ref{thm:strict-detection}.

Represent a finite collection of classes in $A_{P_{r_0}}$ by
classes in $H^0(C(r_0)_u)$ for some $u$. Given $v<u$, their images
in $I=[v,u)$ come from window classes at a common $R>r_0$.
By Theorem~\ref{thm:strict-detection} and $(T)_1$ at $P_R$,
$H^1(C(R)_0)$ is torsion-free. Novikov translation makes
$H^1(C(R)_v)\to H^1(C(R)_u)$ injective. The window exact sequence
therefore lifts the chosen window classes to $H^0(C(R)_u)$.
After restriction to $P_{r_0}$ the differences from the prescribed
classes come from $H^0(C(r_0)_v)$ and have norm less than $e^v$.
Taking $v$ sufficiently negative gives the desired approximation.

Apply this to the finitely many lifts of
Lemma~\ref{lm:finite-tropical-lift} at $r_0$, with error less
than one. Write $y_i\in A_{P_R}$ for the approximating classes.
On the testing polygons at $r_0$ their incident character errors
and their nonincident norms are bounded by a common $q_0<1$.
Choose compact boundary pieces covering the edges in the relative
interiors of these testing neighborhoods along the boundary. At an
endpoint use its vertex cap. Use marked charts extending across
these pieces inside $P_R$, and normalize each incident $m_i$ on the
supporting line at $r_0$. If the chart changes, its primitive character
has the same marked reduction on the piece. The finitely many resulting
errors still have a common bound less than one. The Gauss norms of
$y_i-m_i$ and of the nonincident $y_i$ vary continuously with the
base point, since on a compact subpolytope a convergent Laurent
series is a uniform limit of Laurent polynomials. Thus strict bounds
hold on open neighborhoods of all the chosen pieces.

Let $h_i>0$ be the support constant of $e_i$ at radius one. Its
support at radius $r$ is $rh_i$. Thus if $m_i$ is normalized at
$r_0$, the character normalized at $r$ is
$T^{(r-r_0)h_i}m_i$. For $r$ sufficiently near $r_0$ take the
lifts
\[
 T^{(r-r_0)h_i}\operatorname{res}_{P_R,P_r}(y_i)\in A_{P_r}.
\]
The absolute values of the scaling factors tend to one. Shrinking
the interval about $r_0$ therefore gives a common $q<1$ for the
scaled errors and nonincident norms. Formula \textup{(VN)} gives
norm one for each lift. The nonface products satisfy the same
bound $q$, since at each vertex some factor is nonincident.
At $r_0=1$ we use the edge lifts on $P$ itself and the same
argument on their inward neighborhoods.

Take a smaller relative interval whose closure lies in the interval
of these estimates. The moving edges are covered by the chosen
neighborhoods. On this smaller interval the strips of sufficiently
small width remain in that cover, the displaced vertices and corner
rectangles remain in their endpoint charts, and opposite endpoints
stay separated. These assertions follow by compactness, since the
edge lengths are positive and adjacent conormals remain transverse.
There is therefore a common positive geometric bound for the cuts
used in Section~\ref{sec:localization}. Finitely many smaller
intervals cover $[a,b]$. Choose $\Delta$ so that a step starting
in each smaller interval stays in its larger interval, and all
displacements $d_i=(r-r')h_i$ are below the geometric bounds.
With $q<1$ the largest of the finitely many error bounds, require
also $\Delta\sum_i h_i<-\log q$. Then \eqref{eq:shallow-smallness}
holds for every such cut.
\end{proof}

\begin{proof}[Proof of Theorem~\ref{thm:radial-extension}\,\ref{part:radial-domains}]
Subdivide $[r,R]$ into finitely many intervals of length less than
the $\Delta$ of Lemma~\ref{lem:radial-uniformity}. At each step
Theorem~\ref{thm:shallow-localization}\,\ref{part:shallow-inward}
identifies the inner polygon with a Weierstrass domain in its
predecessor. Affinoid domain inclusions compose, and the composite
on functions is $\operatorname{res}_{P_R,P_r}$ by functoriality.
This also proves independence of the subdivision. The edge lifts
may be chosen separately at each step. Each localization has dense
image, since its adjoined variables are scalar multiples of the
restricted edge lifts, and polynomials in these variables are dense.
It follows that the composite has dense image.

Consider the first step, from $P_R$ to $P_t$. By
Theorem~\ref{tmAffinoidPolygon}\,\ref{part:polygon-affinoid}, the
normalized edge lifts $z_1,\ldots,z_N$ at radius $R$ give an
admissible surjection $\Lambda\langle Z_1,\ldots,Z_N\rangle
\twoheadrightarrow A_{P_R}$. Every edge moves inward by a positive
distance $d_i$, so $|z_i|\leq e^{-d_i}<1$ on $X_t$. Thus every
coordinate in this Tate presentation has absolute value strictly
less than one there. By the analytic interior criterion
\cite[Definition~3.4.2.1]{TemkinBerkovich}, $X_t$ lies in
$\operatorname{Int}(X_R/\Lambda)$. This also places its subdomain
$X_r$ in that interior.
\end{proof}

\subsection{The regular shells and the fiber over the singular value}

\begin{lemma}[Cofinal domains and regular shells]\label{lem:cofinal-shells}
There are radii $1=r_0>r_1>\cdots\to0$ and closed regular collar
domains $S_n\subset X_{r_n}$ corresponding to the base shells
$P_{r_n}\setminus\operatorname{int}P_{r_{n+1}}$. Writing
$X_n=X_{r_n}$ in this lemma and $F_{n+1}$ for the regular domain
over $\partial P_{r_{n+1}}$, we have
\[
 X_n=X_{n+1}\cup S_n,\qquad
 S_n\cap X_{n+1}=F_{n+1},\qquad S_n\cap X_{n+2}=\varnothing .
\]
Moreover $X_{n+1}\subset\operatorname{int}_X X_n$. Floer restrictions
identify the Berkovich open $X\setminus\bigcap_n X_n$, with its
induced analytic structure, with the regular Berkovich mirror over
$P\setminus\{s\}$.
\end{lemma}
\begin{proof}
Subdivide each dyadic interval by Lemma~\ref{lem:radial-uniformity}.
This gives $r_n\downarrow0$ with every step shallow. For $S_n$ take
the union of the closed Laurent strips of this step. All facets move.
The cover and the facet identifications in
Theorem~\ref{thm:shallow-localization} give
$X_n=X_{n+1}\cup S_n$ and $S_n\cap X_{n+1}=F_{n+1}$. They identify
$S_n$ with the regular mirror over the whole closed base shell,
including its corners.

At each point of $F_{n+1}$ some edge lift of $P_{r_{n+1}}$ has absolute
value one. On $X_{n+2}$ every such lift has absolute value at most
$e^{-d_i}<1$ for the next cut. Thus $F_{n+1}$ and hence $S_n$ miss
$X_{n+2}$. The finite union of the earlier $S_0,\ldots,S_{n-1}$ is
closed in $X$ since its strips are compact Berkovich domains. Its
complement is an open neighborhood of $X_{n+1}$ contained in $X_n$.
This proves the inclusion of interiors.

Put $C=\bigcap_nX_n$. The first index at which a point outside $C$
leaves the nested domains places it in the preceding shell.
Conversely each $S_n$ misses $C$. Consecutive shells meet in
$F_{n+1}$, which lies in $S_{n+1}$ because it misses $X_{n+2}$.
All other shell intersections are empty. These are exactly the
intersections of the base shells, and their identifications on the
regular mirror are induced by the same Floer restrictions.

The closed shells form a locally finite cover off $C$, since the
complement of $X_m$ misses every $S_n$ with $n\geq m$.
The regular marked maps thus glue continuously there. To glue the
analytic identifications, take an affinoid domain $V$ in the
Berkovich open $X\setminus C$. The compact sets $V\cap X_m$ are
nested with empty intersection, so $V$ misses some $X_m$.
On $V$ the shell cover is obtained by restricting the finitely
many admissible shallow covers up to that index. The same
identifications give the regular rigid charts on each finite union
of consecutive shells.
\end{proof}

\begin{proof}[Proof of Theorem~\ref{thm:radial-extension}\,\ref{part:radial-maps}]
Each $X_{r_n}$ is nonempty and compact, since $A_{P_{r_n}}$ is a
nonzero affinoid algebra. Hence $C=\bigcap_nX_{r_n}$ is nonempty.
On $X\setminus C$ use the marked regular projection of
Lemma~\ref{lem:cofinal-shells}, and send all of $C$ to $s$.
The shell identifications give
$p^{-1}(P_{r_n})=X_{r_n}$. Given a neighborhood $O$ of
$s$, choose $P_{r_n}\subset O$. The open set
$\operatorname{int}_X X_{r_n}$ contains $X_{r_{n+1}}$ and hence $C$,
and maps into $O$. This proves continuity at $C$. Surjectivity on
the punctured base follows from regular reconstruction, and $C$
supplies the fiber over $s$.
The map is proper because $X$ is compact and $P$ is Hausdorff.

Inserting a prescribed positive radius into a finite subdivision
gives the same affinoid inclusions by
\ref{part:radial-domains}, and the regular maps agree on the refined
shells. This gives $p^{-1}(P_r)=X_r$ for every $r$.
The affinoid-domain identifications are induced by the actual
Floer restriction maps, so their function algebras have the
claimed compatibility.
\end{proof}

\begin{lemma}\label{lem:radial-chart-compatibility}
Let $P$ and $Q$ be two singular polygons containing
the same singular value $s$ and satisfying
$(G_\partial)$, $(G_{\rm rad})$, $(T)_0$ and $(T)_1$.
The analytic spaces and projections constructed in
Theorem~\ref{thm:radial-extension} agree over a
neighborhood of $s$. The identification is induced
locally by Floer restriction and agrees over the
punctured neighborhood with the regular comparison.
\end{lemma}

\begin{proof}
First suppose that $Q\subset P$, without requiring
that they belong to the same radial family.
Floer restriction induces an analytic map
\[
 f_{QP}:X_Q^{\rm rig}\longrightarrow X_P^{\rm rig}.
\]
Over the regular part of $Q$, this map agrees
with the regular analytic identification.
Indeed, on a regular polytope $H\subset Q$ the
composition of restrictions from $P$ to $Q$
and from $Q$ to $H$ is the restriction from $P$
to $H$. The radial shell identifications of
Theorem~\ref{thm:radial-extension} are induced
by these same restrictions.

The map also preserves the projection over $s$.
For any positive radius $r$ in the radial family
of $P$, choose a sufficiently small member $Q_t$
of the radial family of $Q$ contained in $P_r$.
Then the restriction from $P$ to $Q_t$ factors
through $P_r$. Hence $f_{QP}$ maps the fiber
over $s$ in $X_Q$ into every $X_{P_r}$, and
therefore into the fiber over $s$ in $X_P$.

Now let $P,Q$ be arbitrary polygons as in the
statement. Choose members of their radial
families such that
\[
 P_c\subset\operatorname{int}Q_b,\qquad
 Q_b\subset\operatorname{int}P_a,\qquad
 P_a\subset\operatorname{int}Q.
\]
Floer restriction gives maps
\[
 X_{P_c}^{\rm rig}\longrightarrow X_{Q_b}^{\rm rig}
 \longrightarrow X_{P_a}^{\rm rig}\longrightarrow X_Q^{\rm rig}.
\]
By functoriality, the two-step compositions
$X_{P_c}^{\rm rig}\to X_{P_a}^{\rm rig}$ and
$X_{Q_b}^{\rm rig}\to X_Q^{\rm rig}$
are the affinoid-domain embeddings of
Theorem~\ref{thm:radial-extension}.

Set $U=\operatorname{int}P_c$.
All these maps preserve the projection.
By the radial extension theorem the spaces
over $U$ in $X_P^{\rm rig}$, $X_{P_a}^{\rm rig}$ and
$X_{P_c}^{\rm rig}$ are identified, and likewise for
$X_Q^{\rm rig}$ and $X_{Q_b}^{\rm rig}$.
The displayed compositions therefore show
that the maps between the spaces over $U$
in $X_P^{\rm rig}$ and $X_Q^{\rm rig}$ are mutually inverse.

The identification is independent of the chosen
interleaving: after shrinking $U$, any two such
identifications can be compared by restricting
to a common smaller polygon. Functoriality of
Floer restriction gives equality there.
On the remainder of the overlap both identifications agree
with the regular comparison, so they agree on the whole overlap.
The same argument proves the cocycle condition
for three radial charts.
\end{proof}

\subsection{Global reconstruction}

We now glue the radial charts to the regular mirror. Suppose each
isolated singular value admits a polygon satisfying $(G_\partial)$,
$(T)_0$, $(T)_1$ and $(G_{\rm rad})$.
Let $\cP$ be the collection introduced in
Section~\ref{subsec:intro-singular}, consisting of regular convex polytopes
and all singular polygons satisfying
$(G_\partial)$, $(G_{\rm rad})$, $(T)_0$ and $(T)_1$.

For any two members of $\cP$, their intersection
admits a finite covering by members of $\cP$.
Indeed, if the intersection contains a singular
value, we first choose a smaller radial polygon
about that value contained in the interior of
the intersection. The remaining compact
polyhedral region lies in the regular locus
and admits a finite regular polyhedral
subdivision. The same subdivision argument
applies when the intersection is entirely regular.

\begin{theorem}[Global nonsmooth reconstruction]\label{thm:global-nonsmooth}
With these hypotheses there is a normal rigid analytic mirror $Y$
and a continuous proper surjection
\[
 p:Y^{\rm Berk}\longrightarrow B
\]
whose inverse image of $B_{\rm reg}$, with its analytic structure,
is identified with $Y_{\rm reg}^{\rm Berk}$ over $B_{\rm reg}$.
\begin{enumerate}[label=\textup{(\roman*)},ref=(\roman*)]
\item\label{part:global-functions}
For each $H\in\cP$, there is
an affinoid domain $Y_H=\operatorname{Sp}A_H$ in $Y$ with
$p^{-1}(H)=Y_H^{\rm Berk}$. Floer restriction gives an identification
\[
 \cF^0(H)\simeq\Gamma(Y_H,\cO_{Y_H})
\]
compatible with restriction. Coverings in the polyhedral site give
admissible covers of these analytic domains.
\item\label{part:global-symbol}
The surface $Y$ is Gorenstein. For a fixed orientation of $B_{\rm reg}$,
the logarithmic volume form on the regular mirror extends uniquely to
a generator of $\omega_Y$. With divergence for this generator on the
smooth locus of $Y$, the bracket symbols give homomorphisms of
graded BV algebras
\[
 \Sigma_H:\cF^*(H)\longrightarrow\Gamma(Y_H,\PV_Y^{[*]})
\]
compatible with restriction. Each is an isomorphism in degree zero.
On a singular basis polygon $H$, its degree-one map is injective
with finite-length cokernel
and its degree-two map is surjective. The degree-one cokernel and
the degree-two kernel vanish upon localization to the smooth locus
of $Y_H$. If $(T)_2$ holds at $H$, the degree-two map is
an isomorphism there. On any basis member whose affinoid is smooth,
$\Sigma_H$ is a BV isomorphism in all degrees. On regular polytopes
these are the tropically marked comparisons of
Theorem~\ref{thm:global-reconstruction}.
\item\label{part:global-stein}
The space $Y^{\rm Berk}$ is without boundary. The interiors $U$ of
the regular and radial polygons form a neighborhood basis in $B$
with Stein inverse images. Each such $U$ admits an exhaustion by
compact polygons $Q_n$ for which
\[
 \Gamma(p^{-1}(U),\cO)\simeq\varprojlim_n A_{Q_n}
\]
as topological $\Lambda$-algebras. Evaluation identifies $p^{-1}(U)$
homeomorphically with the analytic spectrum of this algebra in the
sense of \cite[Section~4.3]{koso}.
\end{enumerate}
\end{theorem}
\begin{proof}
Shrink the outer radial polygons so that they are pairwise disjoint,
and choose an inner polygon $P_\rho$ in each, with $0<\rho<1$.
Cover the complement of the interiors of the inner polygons by
regular polytopes in that complement. Along each singular chart the
overlap is then the closed collar $P_1\setminus\operatorname{int}P_\rho$.
Subdivide $[\rho,1]$ into finitely many shallow steps. The strips
from these steps cover the closed collar by affinoid domains of
$X_1^{\rm rig}$. Theorem~\ref{thm:shallow-localization} identifies
them with the corresponding domains of the regular mirror. On
overlaps these identifications come from the same Floer restrictions,
so we may glue the singular and regular charts along this finite
cover of the collar. Their projections agree there. A regular basis
polytope is compact and misses the singular values. Its intersection
with each singular chart therefore lies in finitely many regular
shells. Subdividing the polytope along these shells and the gluing
collars recovers its given regular affinoid with its original
restrictions. Lemma~\ref{lem:cofinal-shells} identifies the entire
Berkovich open over the punctured polygon with its regular mirror.

The regular charts are smooth and the singular charts are normal by
Theorem~\ref{tmAffinoidPolygon}. Their projections glue to a
continuous surjection by Theorems~\ref{thm:global-reconstruction}
and~\ref{thm:radial-extension}. These theorems identify $p^{-1}(H)$
with $Y_H^{\rm Berk}$ and its affinoid algebra with $\cF^0(H)$ for
the initial radial and regular charts.

Let $Q\in\cP$ be any other polygon containing a
singular value. By
Lemma~\ref{lem:radial-chart-compatibility},
the radial extension constructed on $X_Q^{\rm rig}$
agrees with the reconstructed mirror over a
neighborhood of the singular value.
Over the punctured polygon both are identified
with the regular mirror, and these identifications
agree on their common domain.

They therefore identify $X_Q^{\rm rig}$ with the analytic
domain of $Y$ whose Berkovich realization is
$p^{-1}(Q)$. To check the analytic-domain assertion,
one may cover $Q$ by a sufficiently small member of the
radial family of $Q$ and finitely many regular polytopes
covering the remaining compact region.
On each of these domains the identification
is an analytic-domain embedding.

For inclusions between members of $\cP$ the
resulting analytic maps are induced by Floer
restriction. This can be checked on regular
polytopes and on a sufficiently small common
radial polygon, where it follows from
functoriality. Consequently finite coverings
by members of $\cP$ pull back to admissible
affinoid coverings.

The interiors of the members of $\cP$ cover $B$, so the proof of
properness in Section~\ref{subsec:canonical-affinoid} applies to
these compact Hausdorff affinoids as well. The regular refinements
give admissible covers by Theorem~\ref{thm:global-reconstruction},
and the shallow cuts do so by Theorem~\ref{thm:shallow-localization}.

Fix the orientation and write $\Omega_{\rm reg}$ for the resulting
form on the regular mirror. Theorem~\ref{thm:polygon-higher-degree}
gives a generator $\Omega_r$ of $\omega_{X_r^{\rm rig}}$ which
agrees with $\Omega_{\rm reg}$ over $\partial P_r$. We claim that
the two agree over all of $P_r\setminus\{s\}$. This part of the
mirror is smooth and connected. Indeed, the strips covering each
shell in Lemma~\ref{lem:cofinal-shells} are connected torus domains,
adjacent strips meet at corners, and successive shells share their
common boundary. The ratio $\Omega_r/\Omega_{\rm reg}$ is therefore
an analytic function on a connected smooth space. It equals one
on the nonempty admissible domain over $\partial P_r$, hence equals
one throughout by the identity theorem
\cite[Lemma~2.1.4]{ConradIrreducible}.

For $r<R$, the generators $\Omega_R|_{X_r^{\rm rig}}$ and $\Omega_r$
thus agree over $\partial P_r$. Their ratio is a unit on $X_r^{\rm rig}$,
and restriction to $\partial P_r$ is injective on $A_{P_r}$ by
Theorem~\ref{tmAffinoidPolygon}. Hence the generators agree on
$X_r^{\rm rig}$. They glue with $\Omega_{\rm reg}$ to give the claimed
trivialization of $\omega_Y$. The same boundary injectivity gives
uniqueness. The Gorenstein property holds on each chart.

To check that the bracket symbols commute with restriction, let
$H'\subset H$ be basis members and put $A=A_H$ and $B=A_{H'}$. The differentials
of functions from $A$ generate $\Omega^{1,\mathrm{an}}_{B/\Lambda}$
over $B$, since $\operatorname{Sp}B$ is an affinoid domain of
$\operatorname{Sp}A$. On these differentials the symbol of a
restricted Floer class agrees with the restriction of its symbol,
by functoriality of the bracket. This determines the two
multiderivations on $B$. In particular a Floer lift of the dual
bivector restricts to a lift on $H'$. Two such lifts differ there by
the kernel of the degree-two symbol. The chartwise BV maps agree
with divergence for the global form and hence define the claimed
map of presheaves. Regular reconstruction and
Theorem~\ref{thm:polygon-higher-degree} give its kernel and cokernel
in each degree on each basis member.

For the Stein assertion, first take $U=\operatorname{int}P_R$ in a
radial chart. Choose $r_n\uparrow R$ with $0<r_n<R$ and put
$Q_n=P_{r_n}$. The interiors of the $Q_n$ cover $U$, so their
affinoids exhaust $p^{-1}(U)$ and every compact subset lies in one
of them. By Theorem~\ref{thm:radial-extension}\,\ref{part:radial-domains},
each affinoid lies in the analytic interior of the next over $\Lambda$,
and the restrictions on functions have dense image.

For the interior of a regular polygon, take smaller homothetic
polygons about an interior point. Laurent polynomials are defined
on the outer torus domain and are dense on each smaller one.
The inclusion $Q_n\subset\operatorname{int}Q_{n+1}$ puts the smaller
torus domain in the topological interior of the larger one. In the
torus, which is without boundary, this is also its absolute analytic
interior. Both exhaustions therefore satisfy the Stein condition of
\cite[Definition~1.1 and Remark~1.2]{MaculanPoineauStein}. They
also place every point of $Y^{\rm Berk}$ in the analytic interior
of an affinoid neighborhood. This proves that $Y^{\rm Berk}$ is
without boundary \cite[Example~4.2.4.2\textup{(i)}]{TemkinBerkovich}.

It remains to identify these domains with the spectra of their
functions. A function on $p^{-1}(U)$ is a compatible family of
functions on its affinoid exhaustion. The supremum norms on these
affinoids give the topology of uniform convergence on compact sets,
so $F=\Gamma(p^{-1}(U),\cO)$ has the stated inverse-limit topology.
Write $Q$ for the outer closed polygon. The restriction
$A_Q\to A_{Q_n}$ has dense image in both constructions and factors
through $F$.

Let $\alpha$ be a continuous multiplicative seminorm on $F$
extending the absolute value of $\Lambda$. Continuity for the
increasing supremum seminorms gives an index $n$ and a constant
$C>0$ with $\alpha(f)\leq C\|f\|_{Q_n}$ for every $f\in F$.
Applying the bound to $f^m$ and taking $m$th roots gives
$\alpha(f)\leq\|f\|_{Q_n}$. Thus $\alpha$ extends from the dense
image of $F$ to a bounded multiplicative seminorm on $A_{Q_n}$,
which defines a point of $\mathcal M(A_{Q_n})$.
Conversely every point over $U$ gives such a seminorm. The functions
from $A_Q$ distinguish these points and induce the subspace topology
from $\mathcal M(A_Q)$. Since the absolute value of every function
in $F$ is continuous on $p^{-1}(U)$, evaluation is a homeomorphism.
\end{proof}

\begin{proof}[Proof of Theorem~\ref{thm:main}]
The hypotheses of
Theorem~\ref{thm:global-nonsmooth} hold.
That theorem gives the analytic surface $Y$,
its projection, the Stein neighborhoods,
the volume form and the comparison with
reflexive polyvector fields.

By the additional condition $(T)_2$ and the
independence of the topological conditions
from the choice of polygon, the degree-two
comparison is an isomorphism on each singular
polygon. The degree-one comparison is injective
with finite-length cokernel supported on the
analytic singular locus, by
Theorem~\ref{thm:polygon-higher-degree}.
The remaining assertions follow from
Theorem~\ref{thm:global-nonsmooth} and the
regular reconstruction theorem.
\end{proof}

\begin{remark}[Rigid points over a singular value]
\label{rem:radial-rigid-points}
Over the regular base the torus coordinates give a rigid point in
each fiber. At a singular value $s$ the same conclusion follows if
$A_s=SH_M^0(\pi^{-1}(s);\Lambda)$ is a nonzero strict affinoid
algebra in its Floer topology. By the affinoid Nullstellensatz the
residue field of a maximal ideal of $A_s$ is then a finite extension
$L$ of $\Lambda$.
Composing the quotient map with Floer restriction gives compatible
bounded homomorphisms $A_{P_r}\to A_s\to L$. They define a
rigid point in the nested domains, hence in their intersection
$p^{-1}(s)$. Whether $A_s$ has the stated affinoid property under
our hypotheses is open.
\end{remark}

\bibliographystyle{alpha}
\bibliography{references}

@article{AbouzaidGromanVarolgunes24,
  title={Framed {$E_2$} structures in {Floer} theory},
  author={Abouzaid, Mohammed and Groman, Yoel and Varolgunes, Umut},
  journal={Advances in Mathematics},
  volume={450},
  pages={109755},
  year={2024},
  publisher={Elsevier}
}

@article{FHS95,
  author={Floer, Andreas and Hofer, Helmut and Salamon, Dietmar},
  title={Transversality in elliptic {Morse} theory for the
         symplectic action},
  journal={Duke Mathematical Journal},
  volume={80},
  number={1},
  pages={251--292},
  year={1995}
}

@book{Bosch14,
  title={Lectures on formal and rigid geometry},
  author={Bosch, Siegfried},
  volume={2105},
  year={2014},
  publisher={Springer}
}

@article{LutkebohmertExtension,
  author={L{\"u}tkebohmert, Werner},
  title={On extension of rigid analytic objects},
  journal={M{\"u}nster Journal of Mathematics},
  volume={15},
  pages={83--166},
  year={2022},
  doi={10.17879/23049549564},
  note={\url{https://www.uni-muenster.de/FB10/mjm/vol_15/mjm_vol_15_03.pdf}}
}

@inproceedings{Conrad2008,
  title={Several approaches to non-archimedean geometry},
  author={Conrad, Brian},
  booktitle={p-adic geometry},
  pages={9--63},
  year={2008}
}

@article{MaculanPoineauStein,
  author={Maculan, Marco and Poineau, J{\'e}r{\^o}me},
  title={Notions of {Stein} spaces in non-{Archimedean} geometry},
  journal={Journal of Algebraic Geometry},
  volume={30},
  number={2},
  pages={287--330},
  year={2021},
  doi={10.1090/jag/764},
  note={\url{https://arxiv.org/abs/1711.04008v2}}
}

@misc{TemkinBerkovich,
  author={Temkin, Michael},
  title={Introduction to {Berkovich} analytic spaces},
  year={2011},
  howpublished={\url{https://arxiv.org/abs/1010.2235v2}},
  note={Version 2, 16 December 2011; numbering in this paper refers to this version}
}

@misc{GromanBoundaryDepth,
  title={Boundary Depth and Deformations of Symplectic Cohomology},
  author={Groman, Yoel},
  year={2026},
  howpublished={\url{https://arxiv.org/abs/2510.17607v2}},
  note={Version 2, 30 September 2026; numbering in this paper refers to this version}
}

@misc{WendlLowDimensions,
  author={Wendl, Chris},
  title={Holomorphic Curves in Low Dimensions: From Symplectic Ruled Surfaces to Planar Contact Manifolds},
  year={2018},
  howpublished={\url{https://www2.mathematik.hu-berlin.de/~wendl/pub/rationalRuled.pdf}},
  note={Version 4.0, 19 March 2018; Sections 8.3.2--8.3.3 and Remark 8.52}
}

@article{SYZ96,
  author={Strominger, Andrew and Yau, Shing-Tung and Zaslow, Eric},
  title={Mirror symmetry is {$T$}-duality},
  journal={Nuclear Physics B},
  volume={479},
  number={1-2},
  pages={243--259},
  year={1996}
}

@incollection{KontsevichSoibelmanTorus,
  author={Kontsevich, Maxim and Soibelman, Yan},
  title={Homological mirror symmetry and torus fibrations},
  booktitle={Symplectic geometry and mirror symmetry
             ({S}eoul, 2000)},
  pages={203--263},
  year={2001},
  publisher={World Scientific Publishing Co., Inc.},
  address={River Edge, NJ}
}

@incollection{koso,
  title={Affine structures and non-archimedean analytic spaces},
  author={Kontsevich, Maxim and Soibelman, Yan},
  booktitle={The unity of mathematics},
  pages={321--385},
  year={2006},
  publisher={Springer}
}

@misc{AbouzaidViterbo,
  author={Abouzaid, Mohammed},
  title={Symplectic cohomology and {Viterbo}'s theorem},
  year={2014},
  howpublished={\url{https://arxiv.org/abs/1312.3354v2}}
}

@incollection{CieliebakLatschev,
  author={Cieliebak, Kai and Latschev, Janko},
  title={The role of string topology in symplectic field theory},
  booktitle={New perspectives and challenges in symplectic field theory},
  series={CRM Proceedings \& Lecture Notes},
  volume={49}, pages={113--146}, year={2009},
  publisher={American Mathematical Society}
}

@article{PascaleffLogCY,
  author={Pascaleff, James},
  title={On the symplectic cohomology of log {Calabi--Yau} surfaces},
  journal={Geometry \& Topology},
  volume={23}, number={6}, pages={2701--2792}, year={2019},
  doi={10.2140/gt.2019.23.2701}
}

@misc{GromanIndicatorV1,
  author={Groman, Yoel},
  title={The local {Floer} cohomology of indicator functions},
  year={2023},
  howpublished={\url{https://arxiv.org/abs/2307.11659v1}},
  note={Version 1}
}

@misc{Paper1,
  author={Groman, Yoel},
  title={The local {Floer} cohomology of indicator functions},
  year={2026},
  howpublished={\url{https://arxiv.org/abs/2307.11659v4}},
  note={Version 4, 8 October 2026; numbering in this paper refers to this version}
}

@misc{GromanVarolgunes2022,
  author={Groman, Yoel and Varolgunes, Umut},
  title={Closed string mirrors of symplectic cluster manifolds},
  year={2024},
  howpublished={\url{https://arxiv.org/abs/2211.07523v2}},
  note={Version 2, 5 July 2024}
}

@article{BEHWZ,
  author={Bourgeois, Fr{\'e}d{\'e}ric and Eliashberg, Yakov and
          Hofer, Helmut and Wysocki, Krzysztof and Zehnder, Eduard},
  title={Compactness results in symplectic field theory},
  journal={Geometry \& Topology},
  volume={7},
  year={2003},
  pages={799--888},
  note={\url{https://arxiv.org/abs/math/0308183}}
}

@article{FishHofer,
  author={Fish, Joel W. and Hofer, Helmut},
  title={Exhaustive {Gromov} compactness for pseudoholomorphic curves},
  journal={Ast{\'e}risque},
  volume={415},
  year={2020},
  pages={87--112},
  note={doi:10.24033/ast.1101; \url{https://arxiv.org/abs/1811.09321}}
}

@misc{McLeanDiscrepancy,
  author={McLean, Mark},
  title={{Reeb} orbits and the minimal discrepancy of an isolated singularity},
  year={2015},
  howpublished={\url{https://arxiv.org/abs/1404.1857v2}},
  note={Version 2; Lemma 8.8}
}

@misc{ConradIrreducible,
  author={Conrad, Brian},
  title={Irreducible components of rigid spaces},
  year={1999},
  howpublished={\url{https://math.stanford.edu/~conrad/papers/irredpaper.pdf}},
  note={Third version, 27 January 1999; Lemma 2.1.4}
}

@misc{Stacks,
  author={{The Stacks Project Authors}},
  title={The {Stacks Project}},
  year={2026},
  howpublished={\url{https://stacks.math.columbia.edu}},
  note={Tags \href{https://stacks.math.columbia.edu/tag/0AV3}{0AV3},
        \href{https://stacks.math.columbia.edu/tag/0AVB}{0AVB},
        \href{https://stacks.math.columbia.edu/tag/0B3D}{0B3D},
        \href{https://stacks.math.columbia.edu/tag/0CQC}{0CQC}, and
        \href{https://stacks.math.columbia.edu/tag/0DW9}{0DW9}}
}

@article{LeungSymington,
  author={Leung, Naichung Conan and Symington, Margaret},
  title={Almost toric symplectic four-manifolds},
  journal={Journal of Symplectic Geometry},
  volume={8}, number={2}, pages={143--187}, year={2010},
  doi={10.4310/JSG.2010.v8.n2.a2}
}

@article{BuzzardLu,
  author={Buzzard, Gregery T. and Lu, Steven S. Y.},
  title={Algebraic surfaces holomorphically dominable by {$\mathbb{C}^2$}},
  journal={Inventiones Mathematicae},
  volume={139}, pages={617--659}, year={2000},
  doi={10.1007/s002229900049}
}

@misc{SchuettShioda,
  author={Sch{\"u}tt, Matthias and Shioda, Tetsuji},
  title={Elliptic surfaces},
  year={2010},
  howpublished={\url{https://arxiv.org/abs/0907.0298v3}}
}

@misc{SawonIsotrivial,
  author={Sawon, Justin},
  title={Isotrivial elliptic {K3} surfaces and {Lagrangian} fibrations},
  year={2014},
  howpublished={\url{https://arxiv.org/abs/1406.1233}}
}

@misc{CheungLin,
  author={Cheung, Man-Wai Mandy and Lin, Yu-Shen},
  title={Some examples of family {Floer} mirrors},
  year={2024},
  howpublished={\url{https://arxiv.org/abs/2101.07079v4}}
}

@article{HeinInstantons,
  author={Hein, Hans-Joachim},
  title={Gravitational instantons from rational elliptic surfaces},
  journal={Journal of the American Mathematical Society},
  volume={25}, pages={355--393}, year={2012},
  doi={10.1090/S0894-0347-2011-00723-6}
}

@misc{HauselToy,
  author={Hausel, Tam{\'a}s},
  title={Geometry of the moduli space of {Higgs} bundles},
  year={2001},
  howpublished={\url{https://arxiv.org/abs/math/0107040}}
}

@misc{FMSWParabolic,
  author={Fredrickson, Laura and Mazzeo, Rafe and Swoboda, Jan and Weiss, Hartmut},
  title={Asymptotic geometry of the moduli space of parabolic
         {$\mathrm{SL}(2,\mathbb C)$}-{Higgs} bundles},
  year={2021},
  howpublished={\url{https://arxiv.org/abs/2001.03682v2}}
}

@misc{LinTropical,
  author={Lin, Yu-Shen},
  title={On the tropical discs counting on elliptic {K3} surfaces
         with general singular fibres},
  year={2019},
  howpublished={\url{https://arxiv.org/abs/1710.10625v2}}
}

@article{TuReconstruction,
  author={Tu, Junwu},
  title={On the reconstruction problem in mirror symmetry},
  journal={Advances in Mathematics},
  volume={256},
  pages={449--478},
  year={2014},
  doi={10.1016/j.aim.2014.02.005}
}

@misc{YuanAn,
  author={Yuan, Hang},
  title={Family {Floer} {SYZ} conjecture for {$A_n$} singularity},
  year={2023},
  howpublished={\url{https://arxiv.org/abs/2305.13554}}
}

@misc{YuanFamilyFloer,
  author={Yuan, Hang},
  title={Family {Floer} program and non-archimedean {SYZ} mirror construction},
  year={2025},
  howpublished={\url{https://arxiv.org/abs/2003.06106}},
  note={Version 3, 10 January 2025}
}
\end{document}